\documentclass[a4paper,oneside,11pt]{article}
 
\usepackage[a4paper]{geometry}
\usepackage{aeguill}                 
\usepackage{graphicx}
\usepackage{caption}
\usepackage{amsmath,amsfonts,amssymb,amsthm}
\usepackage{mathabx}
\usepackage{pstricks,comment} 
\usepackage{pst-plot}
\usepackage{pstricks-add}
\usepackage{multido}
\usepackage{url}
\usepackage{epstopdf,enumerate}
\usepackage{srcltx}
\usepackage[utf8]{inputenc} 
\usepackage[T1]{fontenc} 
\usepackage[english]{babel} 
\usepackage{hyperref}
\usepackage[all]{xy} 
\usepackage{titlesec}
\usepackage{mathabx}
\usepackage{tikz}
\usepackage{fancyhdr}
\usepackage[us,12hr]{datetime}
\fancypagestyle{plain}{
\fancyhf{}
\rfoot{ {\ddmmyyyydate\today}}
\lfoot{\thepage}
}

    \newtheoremstyle{TheoremNum}
        {\topsep}{\topsep}              
        {\itshape}                      
        {}                              
        {\bfseries}                     
        {.}                             
        { }                             
        {\thmname{#1}\thmnote{ \bfseries #3}}

{\theoremstyle {definition} \newtheorem {defi} {Definition}[section]}
{\theoremstyle {plain}  \newtheorem {theo} [defi] {Theorem}}
{\theoremstyle {plain}  \newtheorem {coro} [defi] {Corollary}}
{\theoremstyle {plain} \newtheorem {prop} [defi] {Proposition}}
{\theoremstyle {plain} \newtheorem {lem}[defi] {Lemma}}
{\theoremstyle {plain} \newtheorem {rmq}[defi] {Remark}}
{\theoremstyle {plain} }
{\theoremstyle{TheoremNum} }
{\theoremstyle{TheoremNum} }
{\theoremstyle{TheoremNum} }

\newcommand{\Bij}{\mathrm{Bij}}
\newcommand{\Aut}{\mathrm{Aut}}
\newcommand{\Out}{\mathrm{Out}}
\newcommand{\Inn}{\mathrm{Inn}}
\newcommand{\Mod}{\mathrm{Mod}}
\newcommand{\Comm}{\mathrm{Comm}}

\newcommand{\PO}{\mathbb{P}\mathcal{O}}
\newcommand{\Stab}{\mathrm{Stab}}

\newcommand{\rk}{\mathrm{rk}}
\newcommand{\pr}{\mathrm{prod}}
\newcommand{\id}{\mathrm{id}}
\newcommand{\ad}{\mathrm{ad}}
\newcommand{\FF}{\mathrm{FF}}
\newcommand{\ZZ}{\mathbb{Z}}
\newcommand{\QQ}{\mathbb{Q}}
\newcommand{\NN}{\mathbb{N}}
\newcommand{\RR}{\mathbb{R}}

\newcommand{\dem}{\noindent{\bf Proof. }}

\title{Commensurations of the outer automorphism group  \\ of a universal Coxeter group}
\author{Yassine Guerch}
\date{\today}
\begin{document}
\maketitle
\renewcommand*\labelenumi{(\theenumi)}

\begin{abstract}
This paper studies the rigidity properties of the abstract commensurator of the outer automorphism group of a universal Coxeter
group of rank $n$, which is the free product $W_n$ of $n$ copies of
$\ZZ/2\ZZ$. We prove that for $n\geq 5$ the natural map $\Out(W_n) \to \Comm(\Out(W_n))$ is an isomorphism and that every isomorphism between finite index subgroups of $\Out(W_n)$ is given by a conjugation by an element of $\Out(W_n)$.
\end{abstract}

\section{Introduction}\label{Section Introduction}
Given a group $G$, the \emph{abstract commensurator of $G$}, denoted by $\Comm(G)$, is the group of equivalence classes of isomorphisms between finite index subgroups of $G$. Two such isomorphisms are equivalent if they agree on some common finite index subgroup of their domain. Note that every automorphism of $G$ induces an element of $\Comm(G)$, and in particular the action of $G$ on itself by global conjugation gives a homomorphism $G \to \Comm(G)$. 

The abstract commensurator of $G$ captures a notion of symmetry for the group that is weaker than its group of automorphisms. For instance, the abstract commensurator of $\ZZ^m$ is isomorphic to $\mathrm{GL}(m,\QQ)$ while the abstract commensurator of a nonabelian free group is not finitely generated (see \cite{BartholdiBogopolski10}). However, some groups satisfy strong rigidity properties and the group $\Comm(G)$ is then not much larger than $\Aut(G)$ or $G$ itself. For instance, the Mostow-Prasad-Margulis rigidity theorem and Margulis arithmeticity theorem (see for instance~\cite{Zimmer84}) imply that if $\Gamma$ is an irreducible lattice in a connected semisimple Lie group $G$ with trivial center and no compact factor, and if $G \neq \mathrm{PSL}(2,\RR)$, then $\Gamma$ is a finite index subgroup of $\Comm(\Gamma)$ if and only if $\Gamma$ is not arithmetic, otherwise $\Comm(\Gamma)$ is dense in $G$. In the case of the extended mapping class group of a connected orientable closed surface $S_g$ of genus $g$ at least $3$, we have an even stronger result due to Ivanov \cite{Ivanov97} since the natural homomorphism $\Mod^{\pm}(S_g) \to \Comm(\Mod^{\pm}(S_g))$ is an isomorphism. This result also extends to the case of the mapping class group of a connected orientable surface with genus equal to $2$ and with at least two boundary components. In the context of the outer automorphism group of a free group $F_N$ of rank $N$, Farb and Handel (\cite{farb2007commensurations}) proved that, for $N \geq 4$, the natural map from $\Out(F_N)$ to $\Comm(\Out(F_N))$ is an isomorphism and that every isomorphism between two finite index subgroups of $\Out(F_N)$ extends to an inner automorphism of $\Out(F_N)$. This result was later extended by Horbez and Wade (\cite{HorbezWade20}) to the case $N=3$ using a more geometric approach. Their techniques also enabled them to compute the abstract commensurator of many interesting subgroups of $\Out(F_N)$, like its Torelli subgroup. These rigidity results have been extended to other groups, such as handlebody groups (\cite{hensel2018handlebody}) and big mapping class groups (\cite{BavardDowdallRafi20}).

In this article, we are interested in the outer automorphism group of a universal Coxeter group. Let $n$ be an integer greater than $1$. Let $F=\ZZ/2\ZZ$ be a cyclic group of order $2$ and $W_n= \bigast_n F$ be a universal Coxeter group of rank $n$, that is a free product of $n$ copies of $F$. We prove the following theorem.

\begin{theo}\label{Theo Comm Wn}
Let $n \geq 5$. The natural homomorphism $$\Out(W_n) \to \Comm(\Out(W_n))$$ is an isomorphism.
\end{theo}

The group $\Out(W_2)$ is finite and the group $\Out(W_3)$ is isomorphic to $\mathrm{PGL}(2,\ZZ)$. This gives an almost complete classification except for $n=4$, where our proof for $n \geq 5$ cannot be immediately adapted to this case as $\Out(W_4)$ does not contain any direct product of two nonabelian free groups. Hence the case $n=4$ remains open. Theorem~\ref{Theo Comm Wn} is a major improvement of \cite[Théorème~1.1]{Guerch2020out} which states that, for $n \geq 5$, the only automorphisms of $\Out(W_n)$ are the global conjugations. In turn, Theorem~\ref{Theo Comm Wn} implies that every isomorphism between two finite index subgroups of $\Out(W_n)$ is given by a conjugation by an element of $\Out(W_n)$. The proof of the present Theorem~\ref{Theo Comm Wn} significantly differs from the one of \cite[Théorème~1.1]{Guerch2020out} since the proof of \cite[Théorème~1.1]{Guerch2020out} is based on the study of torsion subgroups of $\Out(W_n)$, whereas $\Out(W_n)$ is virtually torsion free (see~\cite[Corollary~5.5]{Guirardel}).

We now sketch our proof of Theorem~\ref{Theo Comm Wn}. It is inspired by the proof of the similar result in the context of $\Out(F_N)$ given by Horbez and Wade (\cite{HorbezWade20}). However, their proof relies extensively on the possibility of writing a free group as an HNN extension, which is not possible in a universal Coxeter group. Instead, we use the fact that $W_n$ can be written as a free product $W_n=A \ast B$, where $B$ is a finite abelian subgroup of $W_n$. Following a strategy that dates back to Ivanov's work (\cite{Ivanov97}), we study the action of $\Out(W_n)$ on various graphs which are \emph{rigid}, that is, every graph automorphism is induced by an element of $\Out(W_n)$. These graphs include the spine $K_n$ of the Outer space of $W_n$ as defined by Guirardel and Levitt in \cite{Guirardel}, generalizing Culler and Vogtmann's Outer space of a free group (\cite{Vogtmann1986}), or the free splitting graph $\overline{K}_n$ of $W_n$ (see~\cite[Theorem~1.1 and~1.2]{Guerch2020symmetries} and Section~\ref{Section def complexes} for definitions). The proof of Theorem~\ref{Theo Comm Wn} relies on the action of $\Out(W_n)$ on a subset of the vertices of $\overline{K}_n$, called \emph{the set of $W_k$-stars}. Let \mbox{$k \in \{0,\ldots,n-1\}$}. A $W_k$-star is a free splitting $S$ of $W_n$ such that the underlying graph of the induced graph of groups $W_n \backslash S$ is a tree with $n-k$ edges, such that the degree of one of the vertices, called the \emph{center}, is equal to $n-k$, and such that the group associated with the center is isomorphic to $W_{k}$ and the groups associated with the leaves are all isomorphic to $F$. The $W_k$-stars are the analogue for $W_n$ of the roses in the Outer space of a free group. They play a significant role in the proof of other rigidity results for $\Out(W_n)$ (see~\cite{Guerch2020out,Guerch2020symmetries}). 

This allows us to introduce a graph called the \emph{graph of one-edge compatible $W_{n-2}$-stars}, and denoted by $X_n$. It is defined as follows: vertices are $W_n$-equivariant homeomorphism classes of $W_{n-2}$-stars, where two vertices $\mathcal{S}$ and $\mathcal{S}'$ are adjacent if there exist $S \in \mathcal{S}$ and $S' \in \mathcal{S}'$ such that $S$ and $S'$ have both a common refinement and a common collapse. We prove the following result. 

\begin{theo}\label{Theo rigidite Xn intro}
Let $n \geq 5$. The natural homomorphism $$\Out(W_n) \to \Aut(X_n)$$ is an isomorphism.
\end{theo}

Our proof of Theorem~\ref{Theo rigidite Xn intro} requires the rigidity of another graph, called the \emph{graph of $W_{\ast}$-stars}, and denoted by $X_n'$. It is the graph whose vertices are the $W_n$-equivariant homeomorphism classes of $W_k$-stars with $k$ varying in $\{0,\ldots,n-2\}$, where two vertices $\mathcal{S}$ and $\mathcal{S}'$ are adjacent if there exist $S \in \mathcal{S}$ and $S' \in \mathcal{S}'$ such that $S$ refines $S'$ or conversely. We first show that every graph automorphism of $X_n$ induces a graph automorphism of $X_n'$ and that the induced map $\Aut(X_n) \to \Aut(X_n')$ is injective. Using the rigidity of $X_n'$ (see Theorem~\ref{Theo rigidity Xn'}), we show that any graph automorphism of $X_n$ is induced by an element of $\Out(W_n)$.

We then show that every commensuration $f$ of $\Out(W_n)$ induces a graph automorphism of $X_n$. Once we have that result, a general argument (see Proposition~\ref{Prop isomorphism commensuration}) gives the isomorphism between $\Out(W_n)$ and $\Comm(\Out(W_n))$. In order to construct such a homomorphism $\Comm(\Out(W_n)) \to \Aut(X_n)$, we first give an algebraic characterisation of the stabilizers of equivalence classes of $W_{n-2}$-stars. The characterization relies on the examination of maximal abelian subgroups of $\Out(W_n)$ and of direct products of nonabelian free groups in $\Out(W_n)$. In particular, we prove (see Theorem~\ref{Theo product rank Wn}), using the action of $\Out(W_n)$ on a simplicial complex called the \emph{free factor complex of $W_n$}, the following result.

\begin{theo}
Let $n \geq 4$. The maximal number of factors in a direct product of nonabelian free groups contained in $\Out(W_n)$ is equal to $n-3$. 
\end{theo}

One example of such a maximal direct product of nonabelian free subgroups of $\Out(W_n)$ is the following one. Let $W_n=\left\langle x_1,\ldots,x_n \right\rangle$ be a standard generating set for $W_n$ and let $W=\left\langle x_1,x_2,x_3 \right\rangle$. For every $i \geq 4$ and every $w \in W$, let $F_{i,w}$ be the automorphism which fixes $x_j$ for every $j \neq i$ and which sends $x_i$ to $wx_iw^{-1}$. Let $[F_{i,w}]$ be the outer automorphism class of $F_{i,w}$ and let $H_i=\left\langle [F_{i,w}]_{w \in W}\right\rangle$. Then the group $\left\langle H_i \right\rangle_{i \geq 4}$ is a subgroup of $\Out(W_n)$ isomorphic to a direct product of $n-3$ nonabelian free groups. 

The complete characterisation of stabilizers of equivalence classes of $W_{n-2}$-stars being quite technical, we do not give the complete statement in the introduction (see Propositions~\ref{Prop Pwn-2 fixes Wn-1 star} and~\ref{Prop QWn-2 is stab Wn-2}). However, we remark that this characterisation relies on the following key points: the fact that stabilizers of equivalence classes of $W_{n-2}$-stars contain a maximal free abelian subgroup and the fact that it contains a direct product of $n-3$ nonabelian free groups. The characterisation also features a study of the group of twists of a $W_{n-2}$-star, which is a direct product of two virtually nonabelian free groups by a result of Levitt (\cite{levitt2005}) and such that each of which  has finite index in the centralizer in $\Out(W_n)$ of the other.

This characterisation being preserved by commensurations of $\Out(W_n)$, it induces a homomorphism from $\Comm(\Out(W_n))$ to the group $\Bij(VX_n)$ of bijections of the set of vertices of $X_n$. In order to show that this map extends to the edge set of $X_n$, we also present an algebraic characterisation of compatibility of $W_{n-2}$-stars, which is essentially based on the fact that if the intersection of stabilizers of equivalence classes of $W_{n-1}$-stars contains a maximal abelian subgroup of $\Out(W_n)$, then the $W_{n-1}$-stars are pairwise compatible (see Propositions~\ref{Prop PWn-2 compatibility} and~\ref{Prop Pcomp and compatibility}). We deduce that the map $\Comm(\Out(W_n)) \to \Bij(VX_n)$ extends to a map $\Comm(\Out(W_n)) \to \Aut(X_n)$, which completes our proof.

Finally, we prove in the appendix the rigidity of another natural graph endowed with an $\Out(W_n)$-action, called the \emph{graph of $W_{n-1}$-stars}. It is the graph whose vertices are $W_n$-equivariant homeomorphism classes of $W_{n-1}$-stars, where two vertices $\mathcal{S}$ and $\mathcal{S}'$ are adjacent if there exist $S \in \mathcal{S}$ and $S' \in \mathcal{S}'$ such that $S$ and $S'$ have a common refinement. This graph arises naturally in the study of $\Out(W_n)$ and its action on the free splitting graph $\overline{K}_n$ as it is isomorphic to the full subgraph of $\overline{K}_n$ whose vertices are the equivalence classes of $W_k$-stars, with $k$ varying in $\{0,\ldots,n-1\}$. This gives another geometric rigid model for $\Out(W_n)$ (see Theorem~\ref{Theo rigidity graph Wn-1}).

\medskip

{\small{\bf Acknowledgments. } I warmly thank my advisors, Camille Horbez and Frédéric Paulin, for their precious advices and for carefully reading the different versions of this article.}

\section{Preliminaries}\label{Section Preliminaries}

\subsection{Commensurations}

Let $G$ be a group. The \emph{abstract commensurator} of $G$, denoted by $\Comm(G)$, is the group whose elements are the equivalence classes of isomorphisms between finite index subgroups of $G$ for the following equivalence relation. Two isomorphisms between finite index subgroups $f \colon H_1 \to H_2$ and $f' \colon H_1' \to H_2'$ are equivalent if they agree on some common finite index subgroup $H$ of their domains. If $f$ is an isomorphism between finite index subgroups, we denote by $[f]$ the equivalence class of $f$. The identity of $\Comm(G)$ is the equivalence class of the identity map on $G$. Let $[f],[f'] \in \Comm(G)$, and let $f \colon H_1 \to H_2$ and $f' \colon H_1' \to H_2'$ be representatives. The composition law $[f]\cdot [f']$ is given by $[f] \cdot [f']=[f \circ f'|_{f'^{-1}(H_1) \cap H_1'}]$. Note that if $H$ is a finite index subgroup of $G$, then the natural map \mbox{$\Comm(G) \to \Comm(H)$} obtained by restriction is an isomorphism.

Two subgroups $G_1$ and $G_2$ in $G$ are commensurable if $G_1 \cap G_2$ has finite index in both $G_1$ and $G_2$. Being commensurable is an equivalence relation. If $H$ is a subgroup of $G$, we will denote by $[H]$ its commensurability class in $G$. The group $\Comm(G)$ acts on the set of all commensurability classes as follows. Let $[H]$ be the commensurability class of a subgroup $H$. Let $[f] \in \Comm(G)$ and let $f \colon H_1 \to H_2$ be a representative of $[f]$. Then we define $[f] \cdot [H]$ by setting $[f] \cdot [H]=[f(H \cap H_1)]$.

\bigskip

The next result, due to Horbez and Wade, gives a sufficient condition for $\Comm(G)$ to be rigid. It comes from ideas due to Ivanov when studying mapping class groups (see \cite{Ivanov97}). It requires the existence of a graph on which $G$ acts by graph automorphisms.

\begin{prop}\cite[Proposition~1.1]{HorbezWade20}\label{Prop isomorphism commensuration}
Let $G$ be a group. Let $X$ be a graph with no edge-loops, no multiple edges between pairs of vertices and such that $G$ acts on $X$ by graph automorphisms. Let $\Aut(X)$ be the group of graph automorphisms of $X$. Assume that:

\medskip

\noindent{$(1)$ } the natural homomorphism $G \to \Aut(X)$ is an isomorphism,

\medskip

\noindent{$(2)$ } given two distinct vertices $v$ and $w$ of $X$, the groups $\Stab_{G}(v)$ and $\Stab_{G}(w)$ are not commensurable in $G$,

\medskip

\noindent{$(3)$ } the sets $\mathcal{I}=\{[\Stab_{G}(v)] \;|\; v \in VX\}$ and $\mathcal{J}=\{([\Stab_{G}(v)],[\Stab_{G}(w)]) \;|\; vw \in EX\}$ are $\Comm(G)$-invariant (in the latter case with respect to the diagonal action).

\medskip

Then any isomorphism $f \colon H_1 \to H_2$ between finite index subgroups of $G$ is given by the conjugation by an element of $G$ and the natural map $G \to \Comm(G)$  is an isomorphism.
\hfill\qedsymbol
\end{prop}

\subsection{Free splittings and free factor systems of $W_n$}\label{Section def complexes}

Let $n$ be an integer greater than $1$. Let $F=\ZZ/2\ZZ$ be a cyclic group of order $2$ and $W_n= \bigast_n F$ be a universal Coxeter group of rank $n$. A \emph{splitting} of $W_n$ is a minimal, simplicial $W_n$-action on a simplicial tree $S$ such that:

\medskip

\noindent{$(1)$ } The finite graph $W_n \backslash S$ is not empty and not reduced to a point.

\medskip

\noindent{$(2)$ } Vertices of $S$ with trivial stabilizer have degree at least $3$.

\medskip

Here \emph{minimal} means that $W_n$ does not preserve any proper subtree of $S$. A splitting $S$ of $W_n$ is \emph{free} if all edge stabilizers are trivial. A splitting $S'$ is a \emph{blow-up}, or equivalently a \emph{refinement}, of a splitting $S$ if $S$ is obtained from $S'$ by collapsing some edge orbits in $S'$. Two splittings are \emph{compatible} if they have a common refinement. We define an equivalence class in the set of free splittings, where two splittings $S$ and $S'$ are equivalent if there exists a $W_n$-equivariant homeomorphism between them.

\bigskip

A \emph{free factor system} of $W_n$ is a set $\mathcal{F}$ of conjugacy classes of subgroups of $W_n$ which arises as the set of all conjugacy classes of nontrivial point stabilizers in some (nontrivial) free splitting of $W_n$. Equivalently, there exist $k \in \NN-\{0,1\}$ and $[A_1],\ldots,[A_k]$ conjugacy classes of nontrivial, proper subgroups of $W_n$ such that $W_n=A_1 \ast \ldots \ast A_k$ and $\mathcal{F}=\{[A_1],\ldots,[A_k]\}$. The free factor system is \emph{sporadic} if $k=2$, and \emph{nonsporadic} otherwise. The set of all free factor systems of $W_n$ has a natural partial order, where $\mathcal{F} \leq \mathcal{F}'$ if every factor of $\mathcal{F}$ is conjugate into one of the factors of $\mathcal{F}'$. Remark that if $\{x_1,\ldots,x_n\}$ is a standard generating set of $W_n$, then for every free factor system $\mathcal{F}$ of $W_n$ and every $i \in \{1,\ldots,n\}$, there exists $[A] \in \mathcal{F}$ such that $x_i$ is conjugate into $A$. In other words, the free factor system $\{[x_1],\ldots,[x_n]\}$ is a minimum for the partial order on the set of free factor systems of $W_n$.

Let $\mathcal{F}$ be a free factor system of $W_n$. We denote by $\Out(W_n,\mathcal{F})$ the subgroup of $\Out(W_n)$ consisting of all outer automorphisms that preserve all the conjugacy classes of subgroups in $\mathcal{F}$. If $\mathcal{F}=\{[A_1],\ldots,[A_k]\}$, we denote by $\Out(W_n,\mathcal{F}^{(t)})$ the subgroup of $\Out(W_n,\mathcal{F})$ consisting of all outer automorphisms which have a representative whose restriction to each $A_i$ with $i \in \{1,\ldots,k\}$ is a global conjugation by some $g_i \in W_n$.

\bigskip

A \emph{$(W_n,\mathcal{F})$-tree} is an $\RR$-tree equipped with a $W_n$-action by isometries and such that every subgroup of $W_n$ whose conjugacy class belongs to $\mathcal{F}$ is elliptic. A \emph{free splitting of $W_n$ relative to $\mathcal{F}$} is a free splitting of $W_n$ such that every free factor in $\mathcal{F}$ is elliptic. A \emph{free factor of $(W_n,\mathcal{F})$} is a subgroup of $W_n$ which arises as a point stabilizer in a free splitting of $W_n$ relative to $\mathcal{F}$. A free factor of $(W_n, \mathcal{F})$ is \emph{proper} if it is nontrivial, not equal to $W_n$ and not conjugate to an element of $\mathcal{F}$. An element $g \in W_n$ is \emph{$\mathcal{F}$-peripheral} (or simply \emph{peripheral} if there is no ambiguity) if it is conjugate into one of the subgroups of $\mathcal{F}$, and \emph{$\mathcal{F}$-nonperipheral} otherwise. In particular, for every free factor system $\mathcal{F}$ of $W_n$, and every element $x \in W_n$ appearing in a standard generating set of $W_n$, we see that $x$ is $\mathcal{F}$-peripheral.

\subsection{The Outer space of $(W_n,\mathcal{F})$}

We recall the definition of the \emph{unprojectivised Outer space of $(W_n,\mathcal{F})$}, denoted by $\mathcal{O}(W_n,\mathcal{F})$ and introduced by Guirardel and Levitt in \cite{Guirardel}. It is the set of all $(W_n,\mathcal{F})$-equivariant isometry classes $\mathcal{S}$ of metric trees with a nontrivial action of $W_n$, with trivial arc stabilizers and such that a subgroup is elliptic if and only if it is peripheral. The set $\mathcal{O}(W_n,\mathcal{F})$ is equipped with the Gromov-Hausdorff equivariant topology introduced in \cite{Paulin88}. The \emph{projectivised Outer space of $(W_n,\mathcal{F})$}, denoted by $\PO(W_n,\mathcal{F})$, is defined as the space of homothety classes of trees in $\mathcal{O}(W_n,\mathcal{F})$. The spaces $\mathcal{O}(W_n,\mathcal{F})$ and $\PO(W_n,\mathcal{F})$ come equipped with a right action of $\Out(W_n,\mathcal{F})$ given by precomposition of the actions. 

The space $\PO(W_n,\mathcal{F})$ has a natural structure of a simplicial complex with missing faces. Indeed, every element $\mathcal{S} \in \PO(W_n,\mathcal{F})$ defines an open simplex as follows. Let $S$ be a representative of $\mathcal{S}$ such that the sum of the edge lengths of $W_n \backslash S$ is equal to $1$. We associate an open simplex by varying the lengths of the edges, so that the sum of the edge lengths is still equal to $1$. A homothety class $\mathcal{S}' \in \PO(W_n,\mathcal{F})$ of a splitting $S'$ defines a codimension $1$ face of the simplex associated with $\mathcal{S}$ if we can obtain $S'$ from some representative $S$ of $\mathcal{S}$ by contracting one orbit of edges in $S$.

The closure $\overline{\mathcal{O}(W_n,\mathcal{F})}$ of Outer space in the space of all isometry classes of minimal nontrivial $W_n$-actions on $\RR$-trees, equipped with the Gromov-Hausdorff equivariant topology, was identified in \cite{Horbez17} with the space of all \emph{very small $(W_n,\mathcal{F})$-trees}, which are the $(W_n,\mathcal{F})$-trees whose arc stabilizers are either trivial, or cyclic, root-closed and nonperipheral, and whose tripod stabilizers are trivial. The space $\mathbb{P}\overline{\mathcal{O}(W_n,\mathcal{F})}$ equipped with the quotient topology is compact (see \cite[Theorem~1]{Horbez17}).

We recall the definition of a simplicial complex on which the space $\PO(W_n,\mathcal{F})$ retracts \mbox{$\Out(W_n,\mathcal{F})$-equivariantly}, called the \emph{spine of Outer space of $(W_n,\mathcal{F})$} and denoted by $K(W_n,\mathcal{F})$. It is the flag complex whose vertices are the $W_n$-equivariant homeomorphism classes $\mathcal{S}$ of free splittings relative to $\mathcal{F}$ with the property that, if $S \in \mathcal{S}$, then all elliptic subgroups in $S$ are peripheral. Two vertices $\mathcal{S}$ and $\mathcal{S}'$ in $K(W_n,\mathcal{F})$ are linked by an edge if there exist $S \in \mathcal{S}$ and $S' \in \mathcal{S}'$ such that $S$ refines $S'$ or conversely. There is an embedding $F \colon K(W_n,\mathcal{F}) \hookrightarrow \PO(W_n,\mathcal{F})$ whose image is the barycentric spine of $\PO(W_n,\mathcal{F})$. We will from now on identify $K(W_n,\mathcal{F})$ with $F(K(W_n,\mathcal{F}))$. 
If $\mathcal{F}$ consists of exactly $n$ copies of $F$, we simply write $K_n$ for $K(W_n,\mathcal{F})$. In this case the dimension of the simplicial complex $K_n$ is $n-2$. Indeed, if $\mathcal{S}$ is an equivalence class of a free splitting $S$ in $K_n$ such that the number of edges of $W_n \backslash S$ is minimal, then, the number of edges in $W_n \backslash S$ is equal to $n-1$. If $\mathcal{S}$ is an equivalence class of a free splitting $S$ in $K_n$ such that the number of edges of $W_n \backslash S$ is maximal, then $W_n \backslash S$ has $n$ leaves and every vertex of $W_n \backslash S$ that is not a leaf has degree equal to $3$. As $S$ is a tree, this shows that the number of edges in $W_n \backslash S$ is equal to $2n-3$. Since, every splitting $S$ of $K_n$ collapes onto a splitting $S'$ such that $W_n \backslash S'$ has $n-1$ edges, we see that the dimension of $K_n$ is equal to $2n-3-(n-1)=n-2$.

\begin{prop}\label{Prop vcd}
Let $n \geq 3$. The virtual cohomological dimension of $\Out(W_n)$ is equal to $n-2$. In particular, the maximal rank of a free abelian subgroup of $\Out(W_n)$ is equal to $n-2$.
\end{prop}

\dem The group $\Out(W_n)$ acts  cocompactly on $K_n$ with finite stabilizers. Since the dimension of $K_n$ is equal to $n-2$, since the Outer space $\PO(W_n)$ is contractible (see \cite[Theorem~4.2]{Guirardel}) and since $\PO(W_n)$ retracts \mbox{$\Out(W_n)$-equivariantly} on $K_n$, we see that the virtual cohomological dimension of $\Out(W_n)$ is at most equal to $n-2$.

Conversely, let $\{x_1,\ldots,x_n\}$ be a standard generating set of $W_n$. For $i \in \{3,\ldots,n\}$, let $F_i$ be the automorphism sending $x_i$ to $x_1x_2x_ix_2x_1$ and, for every $j \neq i$, fixing $x_j$, and let $[F_i]$ be its image in $\Out(W_n)$. Remark that, for distinct $i,j \in \{3,\ldots,n\}$, the automorphisms $F_i$ and $F_j$ have disjoint support, hence they commute. Since, for every $i \in \{1,\ldots,n\}$, the outer automorphism $[F_i]$ has infinite order, the group $\left\langle [F_i] \right\rangle_{i \geq 3}$ is isomorphic to $\ZZ^{n-2}$. This shows that the virtual cohomological dimension of $\Out(W_n)$ is at least $n-2$ and this concludes the proof.
\hfill\qedsymbol

\bigskip

The \emph{free splitting graph of $W_n$}, denoted by $\overline{K}_n$, is the following graph. The vertices of $\overline{K}_n$ are the $W_n$-equivariant homeomorphism classes of free splittings. Two distinct equivalence classes $\mathcal{S}$ and $\mathcal{S}'$ are joined by an edge in $\overline{K}_n$ if there exist $S \in \mathcal{S}$ and $S' \in \mathcal{S}'$ such that $S$ refines $S'$ or conversely. The free splitting graph of $W_n$ is the $1$-skeleton of the closure of $K_n$ in the space of free splittings of $W_n$.
The group $\Aut(W_n)$ acts on $\overline{K}_n$ on the right by precomposition of the action. As $\Inn(W_n)$ acts trivially on $\overline{K}_n$, the action of $\Aut(W_n)$ induces an action of $\Out(W_n)$ on $\overline{K}_n$.

\subsection{The free factor graph of $(W_n,\mathcal{F})$}

Let $\mathcal{F}$ be a free factor system of $W_n$. We now define a Gromov hyperbolic graph on which $\Out(W_n,\mathcal{F})$ acts by isometries. The \emph{free factor graph relative to $\mathcal{F}$}, denoted by $\FF(W_n,\mathcal{F})$, is the following graph. Its vertices are the $W_n$-equivariant homeomorphism classes of free splittings of $W_n$ relative to $\mathcal{F}$. Two equivalence classes $\mathcal{S}$ and $\mathcal{S}'$ are joined by an edge if there exist $S \in \mathcal{S}$ and $S' \in \mathcal{S}'$ such that $S$ and $S'$ are compatible or share a common nonperipheral elliptic element. The free factor graph is always hyperbolic (see \cite{BestvinaFeighn14,HandelMosher14,Guirardelhorbez19}). The next proposition is due to Guirardel and Horbez. Here, if $H$ is a subgroup of $\Out(W_n)$ and if $\mathcal{F}$ is a free factor system of $W_n$, we say that $\mathcal{F}$ is \emph{$H$-periodic} if there exists a finite index subgroup $H'$ of $H$ such that  $H'(\mathcal{F})=\mathcal{F}$.

\begin{prop}\cite[Theorem~5.1]{Guirardelhorbez19}\label{Prop maximal free factor system bounded orbit}
Let $n \geq 3$ and let $\mathcal{F}$ be a nonsporadic free factor system of $W_n$. Let $H$ be a subgroup of $\Out(W_n,\mathcal{F})$ which acts on $\FF(W_n,\mathcal{F})$ with bounded orbits. Then there exists an $H$-periodic free factor system $\mathcal{F}'$ such that $\mathcal{F} \leq \mathcal{F}'$ and $\mathcal{F} \neq \mathcal{F}'$.
\hfill\qedsymbol
\end{prop}

The Gromov boundary of $\FF(W_n,\mathcal{F})$ has been described in terms of \emph{relatively arational trees} (see the work of Reynolds \cite{Reynolds2012} for the definition of an arational tree in the context of a free group, the work of Bestvina-Reynolds and Hamenstädt (\cite{BestvinaReynolds15,Hamenstadt14} for the description of the boundary in the case of a free group, and the work of Guirardel-Horbez \cite{Guirardelhorbez19} in the case of a free product). A $(W_n,\mathcal{F})$-tree $T$ is \emph{arational} if no proper $(W_n,\mathcal{F})$-free factor acts elliptically on $T$ and, for every proper $(W_n,\mathcal{F})$-free factor $A$, the $A$-minimal invariant subtree of $T$ (that is the union of the axes of the loxodromic elements of $A$ for the action of $W_n$ on $T$, see~\cite[Proposition~3.1]{CulMor87}) is a simplicial $A$-tree in which every nontrivial point stabilizer can be conjugated into one of the subgroups of $\mathcal{F}$. We equip each arational $(W_n,\mathcal{F})$-tree with the \emph{observers' topology}: this is the topology on a tree $T$ such that a basis of open sets is given by the connected components of the complements of points in $T$. We equip the set of arational $(W_n,\mathcal{F})$-trees with an equivalence relation, where two arational $(W_n,\mathcal{F})$-trees are equivalent if they are $W_n$-equivariantly homeomorphic with the observers' topology.

\begin{theo}\cite[Theorem 3.4]{Guirardelhorbez19}\label{Theo Gromov boundary FF}
Let $n \geq 3$. Let $\mathcal{F}$ be a nonsporadic free factor system of $W_n$. The Gromov boundary of $\FF(W_n,\mathcal{F})$ is $\Out(W_n,\mathcal{F})$-equivariantly homeomorphic to the space of all equivalence classes of arational $(W_n,\mathcal{F})$-trees.
\hfill\qedsymbol
\end{theo}

Given $T \in \overline{\mathcal{O}(W_n,\mathcal{F})}$, let $[T]$ be the homothety class of $T$. The \emph{homothetic stabilizer} $\Stab([T])$ is the stabilizer of $[T]$ for the action of $\Out(W_n,\mathcal{F})$ on $\mathbb{P}\overline{\mathcal{O}(W_n,\mathcal{F})}$. Equivalently, $\Phi \in \Out(W_n,\mathcal{F})$ lies in $\Stab([T])$ if there exists a lift $\widetilde{\Phi} \in \Aut(W_n,\mathcal{F})$ of $\Phi$ and a homothety $I_{\widetilde{\Phi}} \colon T \to T$ such that, for all $g \in W_n$ and $x \in T$, we have $I_{\widetilde{\Phi}}(gx)=\widetilde{\Phi}(g)I_{\widetilde{\Phi}}(x)$. The scaling factor of $I_{\widetilde{\Phi}}$ does not depend on the choice of a representative of $\Phi$, and we denote it by $\lambda_T(\Phi)$. This gives a homomorphism

$$\begin{array}{ccl}
\Stab([T]) & \to & \RR_+^* \\
\Phi & \mapsto & \lambda_T(\Phi).
\end{array}
$$
The kernel of this morphism is called the \emph{isometric stabilizer of $T$} and is denoted by $\Stab^{\mathrm{is}}(T)$. It is the stabilizer of $T$ for the action of $\Out(W_n,\mathcal{F})$ on $\overline{O(W_n,\mathcal{F})}$.

\begin{lem}\cite[Lemma~6.1]{Guirardelhorbez19}\label{Lem isometric stabilizer cyclic}
Let $n \geq 3$. Let $\mathcal{F}$ be a nonsporadic free factor system of $W_n$. For every $T \in \overline{\mathcal{O}(W_n,\mathcal{F})}$, the image of the morphism $\lambda_T$ is a cyclic (maybe trivial) subgroup of $\RR_+^*$.
\hfill\qedsymbol
\end{lem}

\begin{lem}\cite[Proposition 13.5]{Guirardelhorbez19laminations}\label{Lem fixed point at infinity FF}
Let $n \geq 3$. Let $\mathcal{F}$ be a nonsporadic free factor system of $W_n$, and let $H$ be a subgroup of $\Out(W_n,\mathcal{F})$. If $H$ fixes a point in $\partial_{\infty}\FF(W_n,\mathcal{F})$, then $H$ has a finite-index subgroup that fixes the homothety class of an arational $(W_n,\mathcal{F})$-tree.

\hfill\qedsymbol
\end{lem}

Finally, we state a proposition due to Guirardel and Horbez concerning the isometric stabilizer of an arational tree.

\begin{prop}\cite[Proposition~6.5]{Guirardelhorbez19}\label{Prop isometric stabilizer free splitting fixed}
Let $n \geq 3$. Let $\mathcal{F}$ be a nonsporadic free factor system of $W_n$, and let $T$ be an arational $(W_n,\mathcal{F})$-tree. Let $H$ be a subgroup of $\Out(W_n,\mathcal{F})$ which is virtually contained in $\Stab^{\mathrm{is}}(T)$. Then $H$ has a finite index subgroup $H'$ which fixes infinitely many $(W_n,\mathcal{F})$-free splittings, and in particular $H$ fixes the conjugacy class of a proper $(W_n,\mathcal{F})$-free factor.
\hfill\qedsymbol
\end{prop}

Note that the statement of Proposition~\ref{Prop isometric stabilizer free splitting fixed} in \cite{Guirardelhorbez19} only mentions that $H'$ fixes one $(W_n,\mathcal{F})$-free splitting, but the proof uses an arbitrary free splitting of $W_n$, so that one can construct infinitely many pairwise distinct free splittings fixed by $H'$ by varying the chosen free splitting of $W_n$.

\subsection{Groups of twists}\label{Section Groups of twists}

Let $S$ be a splitting of $W_n$, let $v \in VS$, let $e$ be a half-edge incident to $v$, and let $z$ be an element of $C_{G_v}(G_e)$. We define the \emph{twist by $z$ around $e$} to be the automorphism $D_{e,z}$ of $W_n$ defined as follows (see~\cite{levitt2005}). Let $\overline{S}$ be the splitting obtained from $S$ by collapsing all the half-edges of $S$ outside of the orbit of the initial half edge of $e$. Then $\overline{S}$ is a tree. Let $\overline{e}$ be the image of $e$ in $\overline{S}$ and let $\overline{v}$ be the image of $v$ in $\overline{S}$. Let $\overline{w}$ be the endpoint of $\overline{e}$ distinct from $\overline{v}$. The automorphism $D_{e,z}$ is defined to be the unique automorphism that acts as the identity on $G_{\overline{v}}$ and as conjugation by $z$ on $G_{\overline{w}}$. The element $z$ is called the \emph{twistor of $D_{e,z}$}. It is well-defined up to composing on the right by an element of $C_{W_n}(G_{\overline{w}}) \cap C_{G_v}(G_e)$. The \emph{group of twists of $S$} is the subgroup of $\Out(W_n)$ generated by all twists around half-edges of $S$.

We now give a description of the stabilizer of a point in $\overline{K}_n$ due to Levitt. If $\mathcal{S} \in V\overline{K}_n$, we denote by $\Stab(\mathcal{S})$ the stabilizer of $\mathcal{S}$ under the action of $\Out(W_n)$. Let $S$ be a representative of $\mathcal{S}$. We denote by $\Stab^0(\mathcal{S})$ the subgroup of $\Stab(\mathcal{S})$ consisting of all elements $F \in \Out(W_n)$ such that the graph automorphism induced by $F$ on $W_n \backslash S$ is the identity.

\begin{prop}\cite[Propositions~2.2, 3.1 and 4.2]{levitt2005}\label{Prop Levitt stab}
Let $n \geq 4$ and $\mathcal{S} \in V\overline{K}_n$. Let $S$ be a representative of $\mathcal{S}$ and let $v_1,\ldots,v_k$ be the vertices of $W_n \backslash S$ with nontrivial associated groups. For $i \in \{1,\ldots,k\}$, let $G_i$ be the group associated with $v_i$. 

\noindent{$(1)$ }  The group $\Stab^0(\mathcal{S})$ fits in an exact sequence $$1 \to \mathcal{T} \to \Stab^0(\mathcal{S}) \to \prod_{i=1}^k \Out(G_i) \to 1,$$ where $\mathcal{T}$ is the group of twists of $S$.

\medskip

\noindent{$(2)$ } The group $\Stab^0(\mathcal{S})$ is isomorphic to 
$$ \prod\limits_{i=1}^k G_i^{\operatorname{deg}(v_i)-1} \rtimes \Aut(G_i),$$
where $\Aut(G_i)$ acts on $G_i^{\operatorname{deg}(v_i)-1}$ diagonally.

\noindent{$(3)$ } The group of twists $\mathcal{T}$ of $S$ is isomorphic to $$\mathcal{T} \simeq \oplus_{i=1}^k G_i^{\operatorname{deg}(v_i)}/Z(G_i),$$ where the center $Z(G_i)$ of $G_i$ is embbeded diagonally in $G_i^{\operatorname{deg}(v_i)}$.
\hfill\qedsymbol
\end{prop}

\begin{rmq}\label{Remark proposition Levitt}
In \cite[Proposition~2.2]{levitt2005}, Levitt shows that the kernel of the natural homomorphism $\Stab^0(\mathcal{S}) \to \prod_{i=1}^k \Out(G_i)$ given by the action on the vertex groups is generated by \emph{bitwists}. Since edge stabilizers are trivial, the group of bitwists is equal to the group of twists. More generally (see~\cite[Proposition~2.3]{levitt2005}), if the outer automorphism group of every edge stabilizer is finite (in particular, if edge stabilizers are isomorphic to $\ZZ$ or to $F$) then the group of twists is a finite index subgroup of the group of bitwists.

Finally, if the centralizer in $W_n$ of an edge stabilizer is trivial, then the group of bitwists about this edge is trivial. Therefore, if the edge stabilizer is not cyclic, then the group of bitwists about this edge is trivial. In all cases, we see that, for every equivalence class $\mathcal{S}$ of a splitting $S$ of $W_n$, the group of twists of $\mathcal{S}$ is a finite index subgroup of the group of bitwists of $W_n$.
\end{rmq}

We establish one last fact about twists about edges whose centralizer is cyclic (see~\cite[Lemma~5.3]{CohenLustig99} for a similar statement in the context of the outer automorphism group of a nonabelian free group).

\begin{lem}\label{Lem twists about cyclic edge central}
Let $n \geq 3$ and let $\mathcal{S}$ be the equivalence class of a splitting $S$. Suppose that there exists an edge $e$ of $S$ with cyclic stabilizer and let $D$ be the outer automorphism class of a twist about $e$. Let $H_{\mathcal{S}}$ be the subgroup of $\Stab^0(\mathcal{S})$ which induces the identity on the edge stabilizer $G_e$ of $e$. Then $D$ is central in $H_{\mathcal{S}}$

In particular, $\Stab^0(\mathcal{S})$ has a finite index subgroup $H_{\mathcal{S}}$ such that $D$ is central in $H_{\mathcal{S}}$.
\end{lem}

\dem Let $U$ be a splitting onto which $S$ collapses (or $S$ itself if $S$ does not have a nontrivial collapse), and let $\mathcal{U}$ be its equivalence class. Then $\Stab^0(\mathcal{S}) \subseteq \Stab^0(\mathcal{U})$. Thus, we may suppose, up to collapsing all orbits of edges of $S$ except the one containing $e$, that $S$ has exactly one orbit of edges. Let $v$ and $w$ be the two endpoints of $e$ and let $G_v$ and $G_w$ be their edge stabilizers. Let $f \in H_{\mathcal{S}}$ and let $F$ be a representative of $f$ such that $F(G_v)=G_v$, $F(G_w)=G_w$ and $F|_{G_e}=\mathrm{id}_{G_e}$ (this representative exists since $f \in H_{\mathcal{S}}$). Let $z \in C_{G_v}(G_e)$ be such that $D_{e,z}$ is a representative of $D$. Then, since $F(z)=z$, for every $x \in W_n$, we have $D_{e,z} \circ F \circ D_{e,z}^{-1}(x)=F(x)$. Hence $f$ and $D$ commutes and $D$ is central in $H_{\mathcal{S}}$. Since the outer automorphism group of a cyclic group is finite, we see that $H_{\mathcal{S}}$ is a finite index subgroup of $\Stab^0(\mathcal{S})$. This concludes the proof.
\hfill\qedsymbol

\section{Geometric rigidity in the graph of $W_k$-stars}\label{Section Wk stars}

We start by defining \emph{$W_k$-stars}, which are the main splittings of interest in this article. 

\begin{defi}
Let $n \geq 3$, and let $k \geq 1$ be an integer. 

\medskip

\noindent{$(1)$ } A free splitting $S$ is a \emph{$k$-edge free splitting} if $W_n \backslash S$ has exactly $k$ edges. 

\medskip

\noindent{$(2)$ } Suppose that $0 \leq k \leq n-2$. A \emph{$W_{k}$-star} is an $(n-k)$-edge free splitting such that:

\begin{itemize}
\item the underlying graph of $W_n \backslash S$ has $n-k+1$ vertices and one of them, called the \emph{center of $W_n \backslash S$}, has degree exactly $n-k$,
\item the group associated with the center of $W_n \backslash S$ is isomorphic to $W_k$ (we use the convention that $W_0=\{1\}$ and that $W_1=F$),
\item the group associated with any leaf of $W_n \backslash S$ is isomorphic to $F$.
\end{itemize}

\noindent{$(3)$ } A \emph{$W_{n-1}$-star} is a one-edge free splitting $S$ such that one of the vertex groups of $W_n \backslash S$ is isomorphic to $W_{n-1}$ while the other vertex group is isomorphic to $F$.
\end{defi}

Note that, in \cite{Guerch2020symmetries}, a $W_{n-1}$-star is called an \emph{$F$-one-edge free splitting}. Using Proposition~\ref{Prop Levitt stab}~$(2)$, we see that, if $k \in \{0,\ldots,n-2\}$, and if $\mathcal{S}$ is the equivalence class of a $W_{k}$-star, then the group $\Stab^{0}(\mathcal{S})$ is isomorphic to $F^{n-k-1} \rtimes \Aut(W_k)$.

Note that, if $S$ is a $W_k$-star with $k \in \{0,\ldots,n-2\}$ and $S'$ is a splitting on which $S$ collapses, then there exists $\ell \in \{k,\ldots,n-1\}$ such that $S'$ is a $W_{\ell}$-star. In particular, for every $k \in \{0,\ldots,n-2\}$, if $S$ is a $W_k$-star, then every one-edge free splitting on which $S$ collapses is a $W_{n-1}$-star. A similar statement is also true for refinements of $W_k$-stars (see~Lemma~\ref{Lem raffinement Wk étoile est Wl étoile}).

\subsection{Rigidity of the graph of $W_{\ast}$-stars}\label{Section graph Wn-1 stars}

We introduce in this section a graph, \emph{the graph of one-edge compatible $W_{n-2}$-stars}, on which $\Out(W_n)$ acts by simplicial automorphisms. We prove that this graph is a rigid geometric model for $\Out(W_n)$. The proof relies on the study of the rigidity of an additional graph on which $\Out(W_n)$ acts, the \emph{graph of $W_{\ast}$-stars}, to be defined after Theorem~\ref{Theo rigidity graph Wn-2}.

\begin{defi}\label{Defi graph Xn}
\noindent{$(1)$ } The \emph{graph of $W_{n-2}$-stars}, denoted by $\widetilde{X}_n$, is the graph whose vertices are the $W_n$-equivariant homeomorphism classes of $W_{n-2}$-stars, where two equivalence classes $\mathcal{S}$ and $\mathcal{S}'$ are joined by an edge if there exist $S \in \mathcal{S}$ and $S' \in \mathcal{S}'$ such that $S$ and $S'$ are compatible.

\medskip

\noindent{$(2)$ } The \emph{graph of one-edge compatible $W_{n-2}$-stars}, denoted by $X_n$, is the graph whose vertices are the $W_n$-equivariant homeomorphism classes of $W_{n-2}$-stars where two equivalence classes $\mathcal{S}$ and $\mathcal{S}'$ are joined by an edge if there exist $S \in \mathcal{S}$ and $S' \in \mathcal{S}'$ such that $S$ and $S'$ have a common refinement which is a $W_{n-3}$-star.
\end{defi}

Note that the adjacency in the graph $X_n$ is equivalent to having both a common collapse (which is a $W_{n-1}$-star) and a common refinement. The graph $X_n$ is a subgraph of $\widetilde{X}_n$. The group $\Aut(W_n)$ acts on $\widetilde{X}_n$ and $X_n$ by precomposition of the action. As $\Inn(W_n)$ acts trivially on $X_n$, the action of $\Aut(W_n)$ induces an action of $\Out(W_n)$. We denote by $\Aut(X_n)$ the group of graph automorphisms of $X_n$. In Section~\ref{Section rigidity graph Wn-2}, we prove the following theorem.

\begin{theo}\label{Theo rigidity graph Wn-2}
Let $n \geq 5$. The natural homomorphism $$\Out(W_n) \to \Aut(X_n)$$ is an isomorphism.
\end{theo}

In order to prove this theorem, we take advantage of the action of $\Out(W_n)$ on another graph, namely the \emph{graph of $W_{\ast}$-stars}, denoted by $X_n'$. The vertices of this graph are the $W_n$-equivariant homeomorphism classes of $W_k$-stars, with $k$ varying in $\{0,\ldots,n-2\}$. Two equivalence classes $\mathcal{S}$ and $\mathcal{S}'$ are joined by an edge if there exist $S \in \mathcal{S}$ and $S' \in \mathcal{S}'$ such that $S$ refines $S'$ or conversely. Note that we have a natural embedding $X_n' \hookrightarrow \overline{K}_n$. We identify from now on $X_n'$ with its image in $\overline{K}_n$. In this section, we prove the following theorem.

\begin{theo}\label{Theo rigidity Xn'}
Let $n \geq 5$. The natural homomorphism $$\Out(W_n) \to \Aut(X_n')$$ is an isomorphism.
\end{theo}

Theorem~\ref{Theo rigidity Xn'} relies on the fact that $X_n'$ contains a rigid subgraph, namely the graph of $\{0\}$-stars and $F$-stars, and denoted by $L_n$. The vertices of this graph are the $W_n$-equivariant homeomorphism classes of $\{0\}$-stars and $F$-stars. Two equivalence classes $\mathcal{S}$ and $\mathcal{S}'$ are joined by an edge if there exist $S \in \mathcal{S}$ and $S' \in \mathcal{S}'$ such that $S$ refines $S'$ or conversely.

We recall the following theorem.

\begin{theo}\cite[Theorem~3.1, Corollary~3.2]{Guerch2020symmetries}\label{Theo rigidity Ln}
Let $n \geq 4$. Let $f$ be an automorphism of $L_n$ preserving the set of $\{0\}$-stars and the set of $F$-stars. Then $f$ is induced by the action of a unique element $\gamma$ of $\Out(W_n)$. In particular, for every $n \geq 5$, the natural homomorphism $$\Out(W_n) \to \Aut(L_n)$$ is an isomorphism.
\hfill\qedsymbol
\end{theo}

The strategy in order to prove Theorem~\ref{Theo rigidity Xn'} is to show that every automorphism of $X_n'$ preserves $L_n$ and that the natural map $\Aut(X_n') \to \Aut(L_n)$ is injective. 

\begin{rmq}
Using the same techniques, we may prove that the \emph{graph of $W_{n-1}$-stars} is rigid. This is done in the appendix (see Theorem~\ref{Theo rigidity graph Wn-1}).
\end{rmq}

First we recall a theorem due to Scott and Swarup.

\begin{theo}\cite[Theorem 2.5]{ScottSwarup}\label{Theo scott swarup}
Let $n \geq 4$. Any set $\{S_1,\ldots,S_k\}$ of pairwise nonequivalent, pairwise compatible, one-edge free splittings of $W_n$ has a unique refinement $S$ such that $W_n \backslash S$ has exactly $k$ edges. Moreover, the equivalence class of $S$ only depends on the equivalence classes of $S_1,\ldots,S_k$. If $S$ is a free splitting such that $W_n \backslash S$ has exactly $k$ edges, then $S$ refines exactly $k$ pairwise nonequivalent one-edge free splittings.
\hfill\qedsymbol
\end{theo}

We also need the following lemma concerning refinements of $W_k$-stars.

\begin{lem}\label{Lem raffinement Wk étoile est Wl étoile}
Let $k, \ell \in \{0,\ldots,n-1\}$ and let $S$ and $S'$ be respectively a $W_k$-star and a $W_{\ell}$-star. If $S$ and $S'$ have a common refinement, then there exists $j \in \{0,\ldots,n-2\}$ and a $W_j$-star $S''$ which refines both $S$ and $S'$. Moreover, $S''$ can be chosen such that $S''$ is a refinement of $S$ and $S'$ with the minimal number of orbits of edges.
\end{lem}

\dem Let $S_1,\ldots,S_{n-k}$ be $n-k$ $W_{n-1}$-stars onto which $S$ collapses and let $S_1',\ldots,S_{n-\ell}'$ be $n-\ell$ $W_{n-1}$-stars onto which $S'$ collapses. Then the set $\{S_1,\ldots,S_{n-k},S_1',\ldots,S_{n-\ell}'\}$ is a set of pairwise compatible $W_{n-1}$-stars. For every $s \in \{1,\ldots,n-k\}$ and every $t \in \{1,\ldots,n-\ell\}$, let $\mathcal{S}_s$ be the equivalence class of $S_s$ and $\mathcal{S}_t'$ be the equivalence class of $S_t'$. Let $n-j=|\{\mathcal{S}_1,\ldots,\mathcal{S}_{n-k},\mathcal{S}_1',\ldots,\mathcal{S}_{n-\ell}'\}|$. By Theorem~\ref{Theo scott swarup}, there exists a free splitting $S''$ with $n-j$ edges which refines every $W_{n-1}$-star of the set $\{S_1,\ldots,S_{n-k},S_1',\ldots,S_{n-\ell}'\}$. But, as $F$ is freely indecomposable, a common refinement of two $W_{n-1}$-stars $U$ and $U'$ is obtained from $U$ by blowing-up an edge at the vertex of $W_n \backslash U$ whose associated group is isomorphic to $W_{n-1}$. Since $U'$ is also a $W_{n-1}$-star, this common refinement has two orbits of edges and the two corresponding leaves have a stabilizer isomorphic to $F$, hence it is a $W_{n-2}$-star. The same argument shows that, if $U_0$ is a $W_{n-1}$-star and if $U_1$ is a $W_k$-star with $k \in \{1,\ldots,n-1\}$ compatible with $U_0$, then a common refinement of $U_0$ and $U_1$ with a minimal number of orbits of edges is either a $W_k$-star (if the equivalence classes of $U_0$ and $U_1$ are adjacent in $\overline{K}_n$) or a $W_{k-1}$-star. Therefore, by induction on $i \in \{1,\ldots,n-\ell\}$, we see that a common refinement of $\{S_1,\ldots,S_{n-k},S_1',\ldots,S_{n-\ell}'\}$ with the minimal number of orbits of edges is a $W_j$-star. This shows that $S''$ is a $W_j$-star. This concludes the proof.
\hfill\qedsymbol

\medskip

Lemma~\ref{Lem raffinement Wk étoile est Wl étoile} implies that the set of $W_k$-stars with $k$ varying in \mbox{$\{0,\ldots,n-1\}$} is closed under taking collapse and taking refinement with a minimal number of orbits of edges.

\begin{lem}\label{Lem Morphism Xn' to Ln}
Let $n \geq 5$. For every $f \in \Aut(X_n')$, we have $f(L_n)=L_n$. Moreover, if $f|_{L_n}=\id_{L_n}$, then $f=\id_{X_n'}$.
\end{lem}

\dem Let $f \in \Aut(X_n')$. The fact that $f(L_n)=L_n$ follows from the fact that vertices of $K_n \cap X_n'$ in $X_n'$ are characterized by the fact that they are the vertices with finite valence. The proof is identical to the proof of~\cite[Proposition~5.1]{Guerch2020symmetries}.

\begin{figure}[ht]
\centering
\begin{tikzpicture}[scale=2]
\draw (0:0) -- (150:1);
\draw (0:0) -- (165:0.91);
\draw (0:0) -- (210:1);
\draw (0:0) node {$\bullet$};
\draw (165:0.91) node {$\bullet$};
\draw (210:1) node {$\bullet$};
\draw (150:1) node {$\bullet$};
\draw[dotted] (173:0.88) -- (200:0.93);
\draw (0:0) node[right, scale=0.9] {$\;\;\left\langle x_2 \right\rangle$};

\draw (165:0.91) node[left, scale=0.9] {$\left\langle x_3 \right\rangle\;$};

\draw (150:1) node[above left, scale=0.9] {$\left\langle x_1 \right\rangle$};
\draw (210:1) node[below left, scale=0.9] {$\left\langle x_n \right\rangle$};
\end{tikzpicture}
\hspace{1cm}
\begin{tikzpicture}[scale=2]
\draw (0:0) -- (150:1);
\draw (0:0) -- (173:0.88);
\draw (0:0) -- (186:0.88);
\draw (0:0) -- (210:1);
\draw (0:0) node {$\bullet$};
\draw (173:0.88) node {$\bullet$};
\draw (186:0.88) node {$\bullet$};
\draw (210:1) node {$\bullet$};
\draw (150:1) node {$\bullet$};
\draw[dotted] (155:0.97) -- (168:0.9);
\draw[dotted] (192:0.91) -- (205:0.98);
\draw (0:0) node[right, scale=0.9] {$\;\;\left\langle x_i \right\rangle$};

\draw (150:1) node[above left, scale=0.9] {$\left\langle x_1 \right\rangle$};
\draw (173:0.88) node[left, scale=0.9] {$\left\langle x_{i-1} \right\rangle$};
\draw (186:0.88) node[left, scale=0.9] {$\left\langle x_{i+1} \right\rangle$};
\draw (210:1) node[below left, scale=0.9] {$\left\langle x_{n} \right\rangle$};
\end{tikzpicture}
\caption{The $F$-stars $X$ (on the left) and $X'$ (on the right) of the proof of Lemma~\ref{Lem Morphism Xn' to Ln}.}\label{Lem Xn' proof X}
\end{figure}
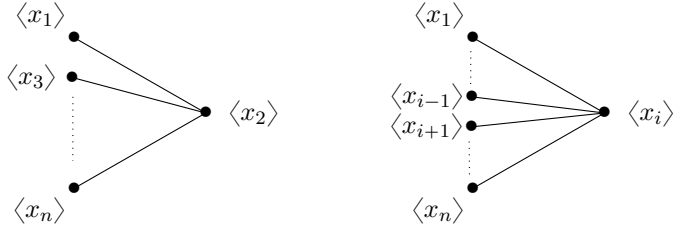

Now suppose that $f|_{L_n}=\id_{L_n}$ and let $\mathcal{S}$ be the equivalence class of a $W_{n-2}$-star $S$. Let us prove that $f(\mathcal{S})=\mathcal{S}$. Let $\{x_1,\ldots,x_n\}$ be a standard generating set of $W_n$ such that the free factor decomposition of $W_n$ induced by $S$ is $$W_n=\left\langle x_1 \right\rangle \ast \left\langle x_2,\ldots,x_{n-1} \right\rangle \ast \left\langle x_n\right\rangle.$$ Let $\mathcal{X}$ be the equivalence class of the $F$-star $X$ depicted in Figure~\ref{Lem Xn' proof X}.

We see that $\mathcal{S}$ and $\mathcal{X}$ are adjacent in $X_n'$. Therefore, as $f(\mathcal{X})=\mathcal{X}$, we see that $f(\mathcal{S})$ and $\mathcal{X}$ are adjacent in $X_n'$. 

Let $\mathcal{S}'$ be the equivalence class of a $W_{n-2}$-star adjacent to $\mathcal{X}$ and distinct from $\mathcal{S}$. Then, as $\mathcal{X}$ and $\mathcal{S}'$ are adjacent, there exist distinct $i,j \in \{1,\ldots,n\}$ and a representative $S'$ of $\mathcal{S'}$ such that the free factor decomposition of $W_n$ induced by $S'$ is $$W_n=\left\langle x_i \right\rangle \ast \left\langle x_1,\ldots,\widehat{x_i},\ldots,\widehat{x_j},\ldots,x_{n} \right\rangle \ast \left\langle x_j\right\rangle.$$ Since $\mathcal{S} \neq \mathcal{S}'$, we may suppose that $i \notin \{1,n\}$. But then $\mathcal{S}$ is adjacent to the equivalence class $\mathcal{X}'$ of the $F$-star $X'$ depicted in Figure~\ref{Lem Xn' proof X} whereas $\mathcal{S}'$ is not adjacent to $\mathcal{X}'$. Since $f(\mathcal{X}')=\mathcal{X'}$, this shows that $f(\mathcal{S}) \neq \mathcal{S}'$.

Finally, let $k \in \{2,\ldots,n-3\}$ and let $\mathcal{S}^{(2)}$ be the equivalence class of a $W_k$-star $S^{(2)}$ which is adjacent to $\mathcal{X}$. We prove that $f(\mathcal{S}) \neq \mathcal{S}^{(2)}$. Since $k \leq n-3$, the underlying graph of $W_n \backslash S^{(2)}$ has at least $3$ edges. Therefore, there exists $i \notin \{1,n\}$ and a leaf $v$ of the underlying graph of $W_n \backslash S^{(2)}$ such that the preimage by the marking of $W_n \backslash S^{(2)}$ of the generator of the group associated with $v$ is $x_i$. But then the equivalence class $\mathcal{S}^{(2)}$ is not adjacent to the equivalence class $\mathcal{X}'$ of the $F$-star $X'$ depicted in Figure~\ref{Lem Xn' proof X}. As $\mathcal{S}$ is adjacent to $\mathcal{X}'$ and as $f(\mathcal{X}')=\mathcal{X}'$, we see that $f(\mathcal{S}) \neq \mathcal{S}^{(2)}$. Therefore, $f(\mathcal{S})=\mathcal{S}$.

The above paragraphs show that $f$ fixes pointwise the set of equivalence classes of $W_{n-2}$-stars. Let $k \in \{2,\ldots,n-3\}$ and let $\mathcal{T}$ be the equivalence class of a $W_k$-star $T$. By Theorem~\ref{Theo scott swarup}, the equivalence class $\mathcal{T}$ is uniquely determined by the set of $W_{n-1}$-stars on which $T$ collapses. Since two distinct equivalence classes of $W_{n-2}$-stars are adjacent in $\overline{K}_n$ to distinct pairs of equivalence classes of $W_{n-1}$-stars, the equivalence class $\mathcal{T}$ is uniquely determined by the set of $W_{n-2}$-stars on which it collapses. Since $f$ fixes pointwise the set of equivalence classes of $W_{n-2}$-stars, we see that $f(\mathcal{T})=\mathcal{T}$. Hence $f=\id_{X_n'}$. This concludes the proof.
\hfill\qedsymbol

\bigskip

\noindent{\bf Proof of Theorem~\ref{Theo rigidity Xn'}. } Let $n \geq 5$. We first prove the injectivity. The homomorphism $\Out(W_n) \to \Aut(L_n)$ is injective by Theorem~\ref{Theo rigidity Ln}. Moreover, the homomorphism $\Out(W_n) \to \Aut(L_n)$ factors through $\Out(W_n) \to \Aut(X_n') \to \Aut(L_n)$. We therefore deduce the injectivity of $\Out(W_n) \to \Aut(X_n')$. We now prove the surjectivity. Let $f \in \Aut(X_n')$. By Lemma~\ref{Lem Morphism Xn' to Ln}, we have a homomorphism $\Phi \colon \Aut(X_n') \to \Aut(L_n)$ defined by restriction. By Theorem~\ref{Theo rigidity Ln}, the automorphism $\Phi(f)$ is induced by an element $\gamma \in \Out(W_n)$. Since the homomorphism $\Aut(X_n') \to \Aut(L_n)$ is injective by Lemma~\ref{Lem Morphism Xn' to Ln}, $f$ is induced by $\gamma$. This concludes the proof.
\hfill\qedsymbol

\subsection{Rigidity of the graph of one-edge compatible $W_{n-2}$-stars}\label{Section rigidity graph Wn-2}

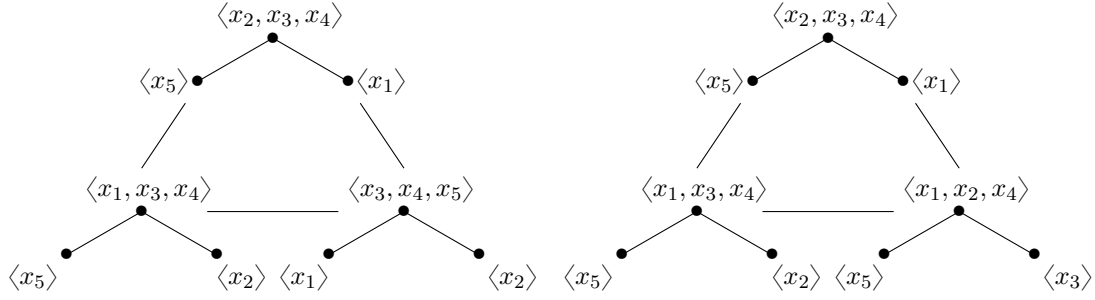
\begin{figure}
\centering
\begin{tikzpicture}[scale=1.15]
\draw (0:0) -- (0.86,-0.5);
\draw (0:0) -- (-0.86,-0.5);
\draw (0:0) node {$\bullet$};
\draw (-0.86,-0.5) node {$\bullet$};
\draw (0.86,-0.5) node {$\bullet$};
\draw (0:0) node[above, scale=0.9] {$\;\;\left\langle x_2,x_3,x_4 \right\rangle$};

\draw (0.86,-0.5) node[right, scale=0.9] {$\left\langle x_1 \right\rangle$};
\draw (-0.86,-0.5) node[left, scale=0.9] {$\left\langle x_5 \right\rangle$};

\draw (-1,-0.75)--(-1.5,-1.5);

\draw (0-1.5,0-2) -- (0.86-1.5,-0.5-2);
\draw (0-1.5,0-2) -- (-0.86-1.5,-0.5-2);
\draw (0-1.5,0-2) node {$\bullet$};
\draw (-0.86-1.5,-0.5-2) node {$\bullet$};
\draw (0.86-1.5,-0.5-2) node {$\bullet$};
\draw (0-1.5,0-2) node[above, scale=0.9] {$\;\;\left\langle x_1,x_3,x_4 \right\rangle$};

\draw (0.86-1.6,-0.5-2) node[below right, scale=0.9] {$\left\langle x_2 \right\rangle$};
\draw (-0.86-1.5,-0.5-2) node[below left, scale=0.9] {$\left\langle x_5 \right\rangle$};

\draw (1,-0.75)--(1.5,-1.5);

\draw (0+1.5,0-2) -- (0.86+1.5,-0.5-2);
\draw (0+1.5,0-2) -- (-0.86+1.5,-0.5-2);
\draw (0+1.5,0-2) node {$\bullet$};
\draw (-0.86+1.5,-0.5-2) node {$\bullet$};
\draw (0.86+1.5,-0.5-2) node {$\bullet$};
\draw (0+1.5,0-2) node[above, scale=0.9] {$\;\;\left\langle x_3,x_4,x_5 \right\rangle$};

\draw (0.86+1.5,-0.5-2) node[below right, scale=0.9] {$\left\langle x_2 \right\rangle$};
\draw (-0.86+1.6,-0.5-2) node[below left, scale=0.9] {$\left\langle x_1 \right\rangle$};

\draw (-0.75,-2)--(0.75,-2);
\end{tikzpicture}
\begin{tikzpicture}[scale=1.15]
\draw (0:0) -- (0.86,-0.5);
\draw (0:0) -- (-0.86,-0.5);
\draw (0:0) node {$\bullet$};
\draw (-0.86,-0.5) node {$\bullet$};
\draw (0.86,-0.5) node {$\bullet$};
\draw (0:0) node[above, scale=0.9] {$\;\;\left\langle x_2,x_3,x_4 \right\rangle$};

\draw (0.86,-0.5) node[right, scale=0.9] {$\left\langle x_1 \right\rangle$};
\draw (-0.86,-0.5) node[left, scale=0.9] {$\left\langle x_5 \right\rangle$};

\draw (-1,-0.75)--(-1.5,-1.5);

\draw (0-1.5,0-2) -- (0.86-1.5,-0.5-2);
\draw (0-1.5,0-2) -- (-0.86-1.5,-0.5-2);
\draw (0-1.5,0-2) node {$\bullet$};
\draw (-0.86-1.5,-0.5-2) node {$\bullet$};
\draw (0.86-1.5,-0.5-2) node {$\bullet$};
\draw (0-1.5,0-2) node[above, scale=0.9] {$\;\;\left\langle x_1,x_3,x_4 \right\rangle$};

\draw (0.86-1.6,-0.5-2) node[below right, scale=0.9] {$\left\langle x_2 \right\rangle$};
\draw (-0.86-1.5,-0.5-2) node[below left, scale=0.9] {$\left\langle x_5 \right\rangle$};

\draw (1,-0.75)--(1.5,-1.5);

\draw (0+1.5,0-2) -- (0.86+1.5,-0.5-2);
\draw (0+1.5,0-2) -- (-0.86+1.5,-0.5-2);
\draw (0+1.5,0-2) node {$\bullet$};
\draw (-0.86+1.5,-0.5-2) node {$\bullet$};
\draw (0.86+1.5,-0.5-2) node {$\bullet$};
\draw (0+1.5,0-2) node[above, scale=0.9] {$\;\;\left\langle x_1,x_2,x_4 \right\rangle$};

\draw (0.86+1.5,-0.5-2) node[below right, scale=0.9] {$\left\langle x_3 \right\rangle$};
\draw (-0.86+1.6,-0.5-2) node[below left, scale=0.9] {$\left\langle x_5 \right\rangle$};

\draw (-0.75,-2)--(0.75,-2);
\end{tikzpicture}
\caption{Two triangles in $X_n$, one corresponding to a $W_{n-3}$-star (on the left) and one corresponding to a $W_{n-4}$-star (on the right).}\label{Figure triangle in Xn}
\end{figure}

In this section, we prove Theorem~\ref{Theo rigidity graph Wn-2}. In order to do so, we construct an injective homomorphism $\Aut(X_n) \to \Aut(X_n')$. First, we need to show some technical results concerning the graph $X_n$. Indeed, let $\Delta$ be a triangle (that is, a cycle of length $3$) in $X_n$, and let $\mathcal{S}_1$, $\mathcal{S}_2$ and $\mathcal{S}_3$ be the vertices of this triangle. By Theorem~\ref{Theo scott swarup}, for every \mbox{$i \in \{1,2,3\}$}, there exists $S_i \in \mathcal{S}_i$ such that $S_1$, $S_2$ and $S_3$ have a common refinement $S$, and we suppose that $S$ has the minimal number of orbits of edges among the common refinements of $S_1$, $S_2$ and $S_3$. Since $S_1$, $S_2$ and $S_3$ are $W_{n-2}$-stars, there exists \mbox{$k \in \{0,\ldots,n-3\}$} such that $S$ is a $W_k$-star. By definition of the adjacency in $X_n$, the splitting $S$ is either a $W_{n-3}$-star or a $W_{n-4}$-star (see Figure~\ref{Figure triangle in Xn}). Our first result shows that we can distinguish these two types of triangles.

\begin{lem}\label{Lem triangle in Xn}
Let $n \geq 5$. Let $\mathcal{S}_1,\mathcal{S}_2$ and $\mathcal{S}_3$ be three equivalence classes of $W_{n-2}$-stars which are pairwise adjacent in $X_n$. Let $S_1$, $S_2$ and $S_3$ be representatives of $\mathcal{S}_1$, $\mathcal{S}_2$ and $\mathcal{S}_3$ which have a common refinement $S$. Suppose that $S$ is the refinement of $S_1$, $S_2$ and $S_3$ which has the minimal number of orbit of edges. Then $S$ is a $W_{n-4}$-star if and only if there exists an equivalence class $\mathcal{S}_4$ of a $W_{n-2}$-star $S_4$ distinct from $\mathcal{S}_1$, $\mathcal{S}_2$ and $\mathcal{S}_3$ such that, for every $i \in \{1,2,3\}$, the equivalence classes $\mathcal{S}_i$ and $\mathcal{S}_4$ are adjacent in $X_n$.
\end{lem}

\dem Suppose first that $S$ is a $W_{n-4}$-star. Let $\{x_1,\ldots,x_n\}$ be a standard generating set of $W_n$ such that the free factor decomposition of $W_n$ induced by $S$ is $$W_n=\left\langle x_1 \right\rangle \ast \left\langle x_2 \right\rangle \ast \left\langle x_3 \right\rangle \ast \left\langle x_4 \right\rangle \ast \left\langle x_5,\ldots,x_n \right\rangle.$$
Since being adjacent in $X_n$ is equivalent to having a common refinement which is a $W_{n-3}$-star and having a common collapse which is a $W_{n-1}$-star, the $W_{n-2}$-stars $S_1$ and $S_2$ share a common collapse $S'$ which is a $W_{n-1}$-star. Let $\mathcal{S}'$ be the equivalence class of $S'$. We claim that there exists an orbit of edges $E$ in $S_3$ such that the splitting obtained from $S_3$ by collapsing every orbit of edges of $S_3$ except $E$ is in $\mathcal{S}'$. Indeed, suppose towards a contradiction that this is not the case. Then, as for every $i \in \{1,2\}$, the equivalence classes $\mathcal{S}_i$ and $\mathcal{S}_3$ are adjacent in $X_n$, we see that, for every $i \in \{1,2\}$, the splittings $S_i$ and $S_3$ share a common collapse onto a $W_{n-1}$-star $S_i'$. Recall that we supposed that there does not exist an orbit of edges $E$ in $S_3$ such that the splitting obtained from $S_3$ by collapsing every orbit of edges of $S_3$ except $E$ is in $\mathcal{S}'$. This implies that for every $i \in \{1,2\}$, the equivalence class $\mathcal{S}_i'$ of $S_i'$ is distinct from $\mathcal{S}'$. Since $S_1$ and $S_2$ are $W_{n-2}$-stars, they collapse onto exactly $2$ distinct $W_{n-1}$-stars. Therefore, for every $i \in \{1,2\}$, the equivalence classes $\mathcal{S}'$ and $\mathcal{S}_i'$ are the two equivalence classes of $W_{n-1}$-stars onto which $S_i$ collapses. It follows that a common refinement of $S_1'$, $S_2'$ and $S'$ is also a common refinement of $S_1$, $S_2$ and $S_3$. But a common refinement of $S_1'$, $S_2'$ and $S_3'$ is a $W_{n-3}$-star. This contradicts the fact that $S$ has the minimal number of edges among common refinements of $S_1$, $S_2$ and $S_3$. Thus $S_3$ collapses onto a $W_{n-1}$-star in the equivalence class $\mathcal{S}'$. Let $j \in \{1,\ldots,4\}$ be such that the free factor decomposition of $W_n$ induced by $S'$ is: $$W_n=\left\langle x_j \right\rangle \ast \left\langle x_1,\ldots, \widehat{x_j},\ldots,x_n \right\rangle.$$ Let $\mathcal{S}_4$ be the equivalence class of the $W_{n-2}$-star $S_4$ whose induced free factor decomposition is: $$W_n= \left\langle x_j \right\rangle \ast \left\langle x_1,\ldots, \widehat{x_5},\ldots,\widehat{x_j},\ldots,x_n \right\rangle \ast \left\langle x_5 \right\rangle.$$ Then, for every $i \in \{1,2,3\}$, the equivalence classes $\mathcal{S}_4$ and $\mathcal{S}_i$ are adjacent in $X_n$. 

\medskip

Conversely, suppose that $S$ is a $W_{n-3}$-star. Let $\{x_1,\ldots,x_n\}$ be a standard generating set of $W_n$ such that the free factor decomposition of $W_n$ induced by $S$ is $$W_n=\left\langle x_1 \right\rangle \ast \left\langle x_2 \right\rangle \ast \left\langle x_3 \right\rangle \ast \left\langle x_4,\ldots,x_n \right\rangle.$$ Then, up to reordering, we may suppose that, for every $i \in \{1,2,3\}$ the free factor decomposition of $W_n$ induced by $S_i$ is: $$W_n=\left\langle x_i \right\rangle \ast \left\langle x_{i+1} \right\rangle \ast \left\langle x_1,\ldots,\widehat{x_i},\widehat{x_{i+1}},\ldots,x_n \right\rangle,$$ where, for $i=3$, the index $i+1$ is taken modulo $3$. Let $\mathcal{S}'$ be the equivalence class of a $W_{n-2}$-star $S'$ adjacent to $\mathcal{S}_1$ in $X_n$ and distinct from $\mathcal{S}_2$ and $\mathcal{S}_3$. Then, up to changing the representative $S'$, there exists $j \in \{1,2\}$ such that $S'$ collapses onto the $W_{n-1}$-star whose associated free factor decomposition is: $$W_n=\left\langle x_j \right\rangle \ast \left\langle x_1,\ldots,\widehat{x_j},\ldots,x_n \right\rangle.$$ If $j=1$, then, as $\mathcal{S}'$ is distinct from $\mathcal{S}_1$ and $\mathcal{S}_3$, we see that $\mathcal{S}'$ is not adjacent to $\mathcal{S}_2$ in $X_n$. If $j=2$, then, as $\mathcal{S}'$ is distinct from $\mathcal{S}_1$ and $\mathcal{S}_2$, we see that $\mathcal{S}'$ is not adjacent to $\mathcal{S}_3$ in $X_n$. In both cases, we see that there exists $i \in \{2,3\}$ such that $\mathcal{S}'$ is not adjacent to $\mathcal{S}_i$. This concludes the proof.
\hfill\qedsymbol

\begin{coro}\label{Coro simplex in Xn}
Let $n \geq 5$. Let $k \geq 4$ and let $\mathcal{S}_1,\ldots,\mathcal{S}_k$ be $k$ equivalences classes of $W_{n-2}$-stars which are pairwise adjacent in $X_n$. For $i \in \{1,\ldots,k\}$, let $S_i$ be a representative of $\mathcal{S}_i$. Let $S$ be a refinement of $S_1,\ldots,S_k$ whose number of orbits of edges is minimal. Then $S$ is a $W_{n-k-1}$-star.
\end{coro}

\dem For every distinct $i,j \in \{1,\ldots,k\}$, the equivalence classes $\mathcal{S}_i$ and $\mathcal{S}_j$ are adjacent in $X_n$. Hence, for every distinct $i,j \in \{1,\ldots,k\}$, there exists a common refinement of $S_i$ and $S_j$ which is a $W_{n-3}$-star. This implies that, for every $p \in \{1,\ldots,k\}$ and for every $i_1,\ldots,i_p \in \{1,\ldots,k\}$, a common refinement of $S_{i_1},\ldots,S_{i_p}$ is obtained from a common refinement of $S_{i_1},\ldots,S_{i_{p-1}}$ whose number of orbits of edges is minimal by adding at most one orbit of edges. We claim that a common refinement of $S_{i_1},\ldots,S_{i_p}$ whose number of orbits of edges is minimal has exactly $p+1$ orbits of edges. Indeed, otherwise there would exist $i,j,\ell \in \{1,\ldots,k\}$ pairwise distinct such that a $W_{n-3}$-star which refines both $S_i$ and $S_j$ also refines $S_{\ell}$. This is not possible by Lemma~\ref{Lem triangle in Xn} since $k \geq 4$. This proves the claim. Taking $p=k$ concludes the proof of the lemma.
\hfill\qedsymbol

\begin{prop}\label{Prop morphism Xn to Xn' well defined}
Let $n \geq 5$. There exists a $\Out(W_n)$-equivariant injective homomorphism $\widetilde{\Phi} \colon \Aut(X_n) \to \Aut(X_n')$.
\end{prop}

\begin{figure}[h]
\centering
\hspace{-0.5cm}
\begin{tikzpicture}[scale=1.5]
\draw (0:0) -- (150:1);
\draw (0:0) -- (210:1);
\draw (0:0) node {$\bullet$};
\draw (210:1) node {$\bullet$};
\draw (150:1) node {$\bullet$};
\draw[dotted] (160:0.93) -- (200:0.93);
\draw (0:0) node[right, scale=0.9] {$\left\langle x_1,\ldots,x_k \right\rangle$};

\draw (150:1) node[above left, scale=0.9] {$\left\langle x_{k+1} \right\rangle$};
\draw (210:1) node[below left, scale=0.9] {$\left\langle x_n \right\rangle$};

\draw (-0.25,-0.75) node {$S$};
\end{tikzpicture}
\hspace{0.4cm}
\begin{tikzpicture}[scale=1.5]
\draw (0:0) -- (210:1);
\draw (0:0) node {$\bullet$};
\draw (210:1) node {$\bullet$};
\draw (0:0) node[above, scale=0.9] {$\left\langle x_1,\ldots,x_{n-1} \right\rangle$};

\draw (210:1) node[below left, scale=0.9] {$\left\langle x_{n} \right\rangle$};
\draw (-0.25,-0.75) node {$S_0$};
\end{tikzpicture}
\hspace{0.4cm}
\begin{tikzpicture}[scale=1.5]
\draw (0:0) -- (150:1);
\draw (0:0) -- (210:1);
\draw (0:0) node {$\bullet$};
\draw (210:1) node {$\bullet$};
\draw (150:1) node {$\bullet$};
\draw (0:0) node[right, scale=0.9, align=left] {$\langle x_1,\ldots,x_k,\ldots,$ \\
$\; \,\widehat{x_{k+i}},\ldots,x_{n-1}\rangle$};

\draw (150:1) node[above left, scale=0.9] {$\left\langle x_{k+i} \right\rangle$};
\draw (210:1) node[below left, scale=0.9] {$\left\langle x_n \right\rangle$};
\draw (-0.25,-0.75) node {$S_i$};
\end{tikzpicture}
\caption{The construction of the map $\Aut(X_n) \to \Aut(X_n')$.}\label{Figure map aut Xn aut Xn'}
\end{figure}
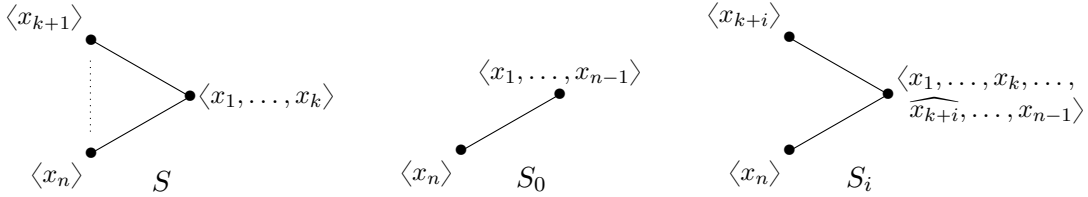

\dem We first explicit a map $\Phi \colon \Aut(X_n) \to \Bij(VX_n')$. Let $f \in \Aut(X_n)$. Let $k \in \{0,\ldots,n-2\}$ and let $\mathcal{S}$ be the equivalence class of a $W_k$-star $S$. If $k=n-2$, then we set $\Phi(f)(\mathcal{S})=f(\mathcal{S})$. If $k \leq n-3$, let $S_0$ be a $W_{n-1}$-star refined by $S$. Let $S_1,\ldots,S_{n-k-1}$ be the $W_{n-2}$-stars such that, for every $i \in \{1,\ldots,n-k-1\}$, $S$ refines $S_i$ and $S_i$ refines $S_0$ (see Figure~\ref{Figure map aut Xn aut Xn'}). For every $i \in \{1,\ldots,n-k-1\}$, let $\mathcal{S}_i$ be the equivalence class of $S_i$, and let $T_i$ be a representative of $f(\mathcal{S}_i)$. By Corollary~\ref{Coro simplex in Xn}, if $n-k-1 \geq 4$, the $W_{n-2}$-stars $T_1,\ldots,T_{n-k-1}$ are refined by a $W_k$-star $T'$. This $W_k$-star is unique up to $W_n$-equivariant homeomorphism by Theorem~\ref{Theo scott swarup}. In the case where $k=n-3$, we have $n-k-1=2$ and, since $f(\mathcal{S}_1)$ and $f(\mathcal{S}_2)$ are adjacent in $X_n$, the splittings $T_1$ and $T_2$ are refined by a $W_{n-3}$-star $T'$ and it is unique up to $W_n$-equivariant homeomorphism by Theorem~\ref{Theo scott swarup}. Finally, when $k=n-4$, Lemma~\ref{Lem triangle in Xn} implies that a common refinement of $T_1$, $T_2$ and $T_3$ with the minimal number of orbits of edges is a $W_{n-4}$-star $T'$, and it is unique up to $W_n$-equivariant homeomorphism by Theorem~\ref{Theo scott swarup}. In all cases, let $\mathcal{T}'$ be the equivalence class of $T'$. We set $\Phi(f)(\mathcal{S})=\mathcal{T}'$.

We now prove that $\Phi$ is well-defined. Let $k \in \{0,\ldots,n-2\}$ and let $\mathcal{S}$ be the equivalence class of a $W_k$-star $S$. Let $S_0$ and $S_0'$ be two distinct $W_{n-1}$-stars onto which $S$ collapses and let $\mathcal{S}_0$ and $\mathcal{S}_0'$ be their equivalence classes. Let $S_1,\ldots,S_{n-k-1}$ be the $W_{n-2}$-stars such that, for every $i \in \{1,\ldots,n-k-1\}$, $S$ refines $S_i$ and $S_i$ refines $S_0$ and let $S_1',\ldots,S_{n-k-1}'$ be the $W_{n-2}$-stars such that, for every $i \in \{1,\ldots,n-k-1\}$, $S$ refines $S_i'$ and $S_i'$ refines $S_0'$. For $i \in \{1,\ldots,n-k-1\}$, let $\mathcal{S}_i$ be the equivalence class of $S_i$ and let $\mathcal{S}_i'$ be the equivalence class of $S_i'$. For every $i \in \{1,\ldots,n-k-1\}$, let $T_i$ be a representative of $f(\mathcal{S}_i)$ and let $T_i'$ be a representative of $f(\mathcal{S}_i')$. Let $T$ be a $W_k$-star which refines $T_1,\ldots,T_{n-k-1}$ and let $T'$ be a $W_k$-star which refines $T_1',\ldots,T_{n-k-1}'$. Finally, let $\mathcal{T}$ be the equivalence class of $T$ and let $\mathcal{T}'$ be the equivalence class of $T'$. We claim that $\mathcal{T}=\mathcal{T}'$. Indeed, we first remark that there exist $i,j \in \{1,\ldots,n-k-1\}$ such that $\mathcal{S}_i=\mathcal{S}_j'$: it is the equivalence class of the $W_{n-2}$-star which refines both $S_0$ and $S_0'$. Up to reordering, we may suppose that $i=j=1$, that $S_1=S_1'$ and that $T_1=T_1'$. Therefore, both $T$ and $T'$ collapse onto $T_1$. 

Let $U_2,\ldots,U_{n-k-1}$ be the $W_{n-3}$-stars such that, for every $j \in \{2,\ldots,n-k-1\}$, the $W_{n-3}$-star $U_j$ refines $S_1$ and $U_j$ is refined by $S$. For every $j \in \{2,\ldots,n-k-1\}$ there exist $\ell,\ell' \in \{2,\ldots,n-k-1\}$ such that $S_{\ell}$ and $S_{\ell'}'$ are refined by $U_j$. Therefore, the application $g \colon \{2,\ldots,n-k-1\} \to \{2,\ldots,n-k-1\}$ sending $\ell$ to $\ell'$ is a bijection. Thus, we may suppose that $g$ is the identity, that is, we may suppose that $j=\ell=\ell'$. It follows that for every $j \in \{2,\ldots,n-k-1\}$, the equivalence class of the $W_{n-3}$-star which refines $S_1$ and $S_j$ is the same one as the equivalence class of the $W_{n-3}$-star which refines $S_1$ and $S_j'$. Therefore, for every $i \in \{2,\ldots,n-k-1\}$, the set $\{\mathcal{S}_1,\mathcal{S}_i,\mathcal{S}_i'\}$ defines a triangle in $X_n$ which corresponds to the equivalence class of a $W_{n-3}$-star. By Lemma~\ref{Lem triangle in Xn}, for every $i \in \{2,\ldots,n-k-1\}$, the set $\{f(\mathcal{S}_1),f(\mathcal{S}_i),f(\mathcal{S}_i')\}$ defines a triangle in $X_n$ which corresponds to the equivalence class of a $W_{n-3}$-star. Thus, up to changing the representative $T_i'$, for every $i \in \{1,\ldots,n-k-1\}$, the $W_{n-3}$-star which refines $T_1$ and $T_i$ is the same one as the $W_{n-3}$-star which refines $T_1$ and $T_i'$. As $\mathcal{T}$ and $\mathcal{T}'$ are characterized by the set of equivalence classes of $W_{n-3}$-stars which collapses onto $T_1$ and on which $T$ and $T'$ collapse, we see that $\mathcal{T}=\mathcal{T}'$. Therefore, the map $\Phi(f) \colon VX_n' \to VX_n'$ is well-defined. As $\Phi(f) \circ \Phi(f^{-1})=\Phi(f \circ f^{-1})=\mathrm{id}$, we see that $\Phi(f)$ is a bijection.

We now prove that the application $\Phi \colon \Aut(X_n) \to \Bij(VX_n')$ induces a monomorphism $\widetilde{\Phi} \colon \Aut(X_n) \to \Aut(X_n')$. Let $f \in \Aut(X_n)$ and let us prove that $\Phi(f)$ preserves $EX_n'$. Let $\mathcal{S},\mathcal{S}'$ be adjacent vertices in $X_n'$. Up to exchanging the roles of $\mathcal{S}$ and $\mathcal{S}'$, we may suppose that there exist $S \in \mathcal{S}$ and $S' \in \mathcal{S}'$ such that $S'$ collapses onto $S$. Let $k, \ell \in \{1,\ldots,n-2\}$ be such that $S$ is a $W_k$-star and $S'$ is a $W_{k-\ell}$-star. Let $S_0$ be a $W_{n-1}$-star such that $S$ refines $S_0$. Let $S_1,\ldots,S_{n-k-1}$ be the $W_{n-2}$-stars such that, for every $i \in \{1,\ldots,n-k-1\}$, $S$ refines $S_i$ and $S_i$ refines $S_0$. As $S'$ refines $S$, there exist $\ell$ $W_{n-2}$-stars $S_{n-k},\ldots,S_{n-k+\ell-1}$ such that the $W_{n-2}$-stars $S_1,\ldots,S_{n-k+\ell-1}$ are the $n-k+\ell-1$ $W_{n-2}$-stars which collapse onto $S_0$ and which are refined by $S'$. For every $i \in \{1,\ldots,n-k+\ell-1\}$, let $\mathcal{S}_i$ be the equivalence class of $S_i$. By definition of $\Phi(f)$, there exist a representative $T$ of $\Phi(f)(\mathcal{S})$ and representatives $T_1,\ldots,T_{n-k-1}$ of $f(\mathcal{S}_1),\ldots,f(\mathcal{S}_{n-k-1})$ such that $T$ is a common refinement of $T_1,\ldots,T_{n-k-1}$. Moreover, there exist a representative $T'$ of $\Phi(f)(\mathcal{S}')$ and representatives $T_{n-k},\ldots,T_{n-k+\ell-1}$ of $f(\mathcal{S}_{n-k}),\ldots,f(\mathcal{S}_{n-k+\ell-1})$ such that $T'$ is a common refinement of $f(\mathcal{S}_1),\ldots,f(\mathcal{S}_{n-k+\ell-1})$. As $\{f(\mathcal{S}_1),\ldots,f(\mathcal{S}_{n-k-1})\}$ is a subset of $\{f(\mathcal{S}_1),\ldots,f(\mathcal{S}_{n-k+\ell-1})\}$, we see that $f(\mathcal{S})$ and $f(\mathcal{S}')$ are adjacent. This shows that the application $\Phi(f) \colon VX_n \to VX_n'$ induces a homomorphism $\widetilde{\Phi} \colon \Aut(X_n) \to \Aut(X_n')$. Finally, the facts that $\widetilde{\Phi}$ is injective and is $\Out(W_n)$-equivariant follow from the fact that, for every equivalence class $\mathcal{S}$ of $W_{n-2}$-stars, we have $f(\mathcal{S})=\Phi(f)(\mathcal{S})$. This concludes the proof.
\hfill\qedsymbol

\bigskip

\noindent{\bf Proof of Theorem~\ref{Theo rigidity graph Wn-2}. } Let $n \geq 5$. We first prove the injectivity. The homomorphism $\Out(W_n) \to \Aut(X_n')$ is injective by Theorem~\ref{Theo rigidity Xn'}. Moreover, the homomorphism $\Out(W_n) \to \Aut(X_n')$ factors through $\Out(W_n) \to \Aut(X_n) \to \Aut(X_n')$. We therefore deduce the injectivity of $\Out(W_n) \to \Aut(X_n)$. We now prove the surjectivity. Let $f \in \Aut(X_n)$. By Proposition~\ref{Prop morphism Xn to Xn' well defined}, we have a homomorphism $\widetilde{\Phi} \colon \Aut(X_n) \to \Aut(X_n')$. By Theorem~\ref{Theo rigidity Xn'}, the automorphism $\widetilde{\Phi}(f)$ is induced by an element $\gamma \in \Out(W_n)$. Since the homomorphism $\Aut(X_n) \to \Aut(X_n')$ is injective by Proposition~\ref{Prop morphism Xn to Xn' well defined}, $f$ is induced by $\gamma$. This concludes the proof.
\hfill\qedsymbol

\section{The group of twists of a $W_{n-1}$-star}

In this section, we study the centralizers in $\Out(W_n)$ of twists about a $W_{n-1}$-star. We first need some preliminary results about stabilizers of free factors of $W_n$ isomorphic to $W_{n-1}$.

Let $\{x_1,\ldots,x_n\}$ be a standard generating set of $W_n$. For distinct $i,j \in \{1,\ldots,n\}$, let $\sigma_{j,i} \colon W_n \to W_n$ be the automorphism sending $x_j$ to $x_ix_jx_i$ and, for $k \neq j$, fixing $x_k$. For distinct $i,j \in \{1,\ldots,n\}$, let $(i\;j)$ be the automorphism of $W_n$ switching $x_i$ and $x_j$ and, for $k \neq i,j$, fixing $x_k$. The following theorem is due to Mühlherr.

\begin{theo}\cite[Theorem~B]{muhlherr1997}\label{Theo generating set Aut}
Let $n\geq 2$. The set $\{\sigma_{i,j} \, |i \neq j \} \cup \{(i\;j) \,| i \neq j\}$ is a generating set of $\Aut(W_n)$.
\end{theo}

We now introduce a finite index subgroup of $\Out(W_n)$ which will be used throughout the remainder of this paper. For every $i,j \in \{1,\ldots,n\}$ distinct, both $\sigma_{i,j}$ and $(i\;j)$ preserve the set of conjugacy classes $\{[x_1],\ldots,[x_n]\}$. Since $\{\sigma_{i,j} \, |i \neq j \} \cup \{(i\;j) \,| i \neq j\}$ generates $\Aut(W_n)$ by Theorem~\ref{Theo generating set Aut}, we see that we have a well-defined homomorphism $\Out(W_n) \to \Bij(\{[x_1],\ldots,[x_n]\})$. Let $C_n$ be the kernel of this homomorphism. The group $C_n$ has finite index in $\Out(W_n)$. We will mostly work in $C_n$ from now on because of the following lemma.

\begin{lem}\label{Lem Cn acts as the identity on graph}
Let $n \geq 3$ and let $f \in C_n$. Suppose that $f$ fixes the equivalence class $\mathcal{S}$ of a free splitting $S$. Then the graph automorphism of the underlying graph of $W_n \backslash S$ induced by $f$ is the identity. Therefore we have $\Stab_{C_n}(\mathcal{S})=\Stab_{C_n}^0(\mathcal{S})$.
\end{lem}

\dem The underlying graph $\overline{W_n \backslash S}$ of $W_n \backslash S$ is a tree. Moreover, since $S$ is a free splitting, if $L$ is the set of leaves of $\overline{W_n \backslash S}$, then the set $\{[G_v]\}_{v \in L}$ is a free factor system of $W_n$. Note that, as $\{[x_1],\ldots,[x_n]\}$ is a free factor system of $W_n$ which is minimal for inclusion, for every $i \in \{1,\ldots,n\}$, there exists one $v \in VS$ such that $x_i \in G_v$. Since $S$ is a free splitting, for every $i \in \{1,\ldots,n\}$, the element $x_i$ is contained in a unique vertex group. Moreover, for every $v \in L$, there exist $k \in \{0,\ldots,n-1\}$ and $\{i_1,\ldots,i_k\} \subseteq \{1,\ldots,n\}$ such that $G_v$ is isomorphic to $W_k$ and $\{[x_{i_1}] \cap G_v,\ldots,[x_{i_k}]\cap G_v\}$ is a free factor system of $G_v$. As $f \in C_n$, and as $f$ fixes $\mathcal{S}$, it follows that, for every $v \in L$, we have $f([G_v])=[G_v]$. Hence the graph automorphism $\widehat{f}$ of $\overline{W_n \backslash S}$ induced by $f$ acts as the identity on $L$. As any graph automorphism of a finite tree is determined by its action on the set of leaves, it follows that $\widehat{f}=\id$. This concludes the proof.
\hfill\qedsymbol

\begin{rmq}
The subgroup $C_n$ of $\Out(W_n)$ is our (weak) analogue of the subgroup $\mathrm{IA}_N(\ZZ/3\ZZ)$ of $\Out(F_N)$, which is defined as the kernel of the natural homomorphism $\Out(F_N) \to \mathrm{GL}(N,\ZZ/3\ZZ)$. Indeed, the group $\mathrm{IA}_N(\ZZ/3\ZZ)$ satisfies a statement similar to Lemma~\ref{Lem Cn acts as the identity on graph}, but it has the additional property that if $\phi \in \mathrm{IA}_N(\ZZ/3\ZZ)$ has a periodic orbit in the free splitting graph of $F_N$, then the cardinality of this orbit is equal to $1$. In the context of $C_n$, we do not know if $C_n$ contains a torsion free finite index subgroup which satisfies this property.
\end{rmq}

The next lemma relates the stabilizer of a free factor of $W_n$ isomorphic to $W_{n-1}$ and the stabilizer of a $W_{n-1}$-star.

\begin{lem}\label{Lem stab corank one Wn-1 star}
Let $n \geq 3$. Let $A$ be a free factor of $W_n$ isomorphic to $W_{n-1}$. Then, up to $W_n$-equivariant homeomorphism, there exists a unique free splitting $S$ in which $A$ is elliptic. In particular, if $f \in \Out(W_n)$ is such that $f([A])=[A]$, then $f$ fixes the equivalence class of $S$.
\end{lem}

\dem By definition of a free factor, there exists a free splitting $S$ of $W_n$ such that $A$ is elliptic in $S$. This proves the existence. We now prove the uniqueness statement. We may assume that $\{x_1,\ldots,x_{n-1}\}$ is a standard generating set of $A$ and $x_n \in W_n$ is such that $$W_n=A \ast \left\langle x_n\right\rangle.$$ Then, the free factor system $\mathcal{F}=\{[A],[\left\langle x_n \right\rangle] \}$ is a sporadic free factor system which contains $[A]$. Let $\mathcal{F}'$ be a free factor system of $W_n$ which contains $[A]$. Since the free factor system $\{[\left\langle x_1 \right\rangle],\ldots,[\left\langle x_n \right\rangle]\}$ is the minimal element of the set of free factor systems of $W_n$, we see that there exists $[B] \in \mathcal{F}'$ such that $x_n \in B$. As $\mathcal{F}'$ contains $[A]$ and as $W_n=A \ast \left\langle x_n \right\rangle$, it follows that $W_n=A \ast B$ and that $B \subseteq \left\langle x_n \right\rangle$. Therefore $[B]=[\left\langle x_n \right\rangle]$ and $\mathcal{F}'=\{[A],[\left\langle x_n\right\rangle]\}$. We deduce that $\mathcal{F}$ is the unique nontrivial free factor system which contains $[A]$. But the spine $K(W_n,\mathcal{F})$ of the Outer space relative to $\mathcal{F}$ is reduced to a point, i.e. it is reduced to a unique equivalence class of free splittings. This proves the uniqueness statement.
\hfill\qedsymbol

\begin{rmq}
In the context of $\Out(F_N)$, the analogue of the splitting given by Lemma~\ref{Lem stab corank one Wn-1 star} is the following one. Let $[A]$ be the conjugacy class of a free factor of $F_N$ isomorphic to $F_{N-1}$. Then the canonical splitting associated with $A$ is the splitting corresponding to the HNN extension $F_N=A \ast$ over the trivial group. However, there does not exist a natural choice (up to conjugacy) of an element $g \in F_N$ such that $\{[A],[g]\}$ is a free factor system of $F_N$.
\end{rmq}

Let $S$ be a splitting with exactly one orbit of edges, whose stabilizer is root-closed and isomorphic to $\ZZ$. Then the group of twists of $S$ is isomorphic to $\ZZ$ by a result of Levitt (see~\cite[Proposition~3.1]{levitt2005}). The next proposition is similar to a result in the case of the outer automorphism group of a free group (see~\cite{CohenLustig95} and \cite[Lemma~2.7]{HorbezWade20}). Recall that an element $w \in W_n$ is \emph{root-closed} if there does not exist $w_0 \in W_n$ and an integer $n \geq 2$ such that $w=w_0^n$.

\begin{lem}\label{Lem compatibility splittings twists}
Let $n \geq 3$. Let $A$ be a free factor of $W_n$ isomorphic to $W_{n-1}$ and let $w \in A$ be a root-closed element of infinite order. Let $x \in W_n$ be such that $W_n=A \ast \left\langle x \right\rangle$. Let $\mathcal{S}$ be the equivalence class of a splitting $S$ whose associated amalgamated decomposition of $W_n$ is the following: $$W_n= A \ast_{\left\langle w \right\rangle} \left(\left\langle w \right\rangle \ast \left\langle x \right\rangle \right).$$ Let $D$ be a nontrivial twist about $S$. Let $\mathcal{R}$ be the equivalence class of a free splitting $R$ of $W_n$ such that $D(\mathcal{R})=\mathcal{R}$. Let $R'$ and $S'$ be  metric representatives of $R$ and $S$, let $\mathcal{R}'$ and $\mathcal{S}'$ be their $W_n$-equivariant isometry classes and let $[\mathcal{R}']$ and $[\mathcal{S}']$ be their homothety classes.

\medskip

\noindent{$(1)$ } In $\mathbb{P}\overline{\mathcal{O}(W_n)}$, there exists an increasing function $\psi \colon \NN \to \NN$ such that $$\lim_{n \to \infty} D^{\psi(n)}([\mathcal{R}'])=[\mathcal{S}'].$$

\noindent{$(2)$ } The splittings $S$ and $R$ are compatible.
\end{lem}

\dem We prove the first part. As $\mathbb{P}\overline{\mathcal{O}(W_n)}$ is compact, up to passing to a subsequence, there exists a sequence $(\lambda_n)_{n \in \NN} \in (\RR_+^*)^{\NN}$ and an $W_n$-equivariant isometry class $\mathcal{T}$ of an $\RR$-tree $T$ such that $$\lim_{n \to \infty}\lambda_nD^n(\mathcal{R}')=\mathcal{T}.$$ Since translation length functions are continuous for the Gromov-Hausdorff topology (see~\cite{Paulin88}), for every $g \in W_n$, we have: $$\lim_{n \to \infty} \lambda_n\left\lVert g \right\rVert_{D^n(\mathcal{R}')}= \left\lVert g \right\rVert_{\mathcal{T}},$$ where $\left\lVert g \right\rVert_{\mathcal{T}}$ is the translation length of $g$ in $T$. Hence, for every $g \in W_n$, the limit $\lim_{n \to \infty} \lambda_n\left\lVert g \right\rVert_{D^n(\mathcal{R}')}$ is finite. But as $D$ has infinite order, we have $\lim_{n \to \infty} \lambda_n=0$. As there exists a representative $\phi \in \Aut(W_n)$ of $D$ such that $\phi_A=\id_A$, for every $g\in A$, we have: $$\lim_{n \to \infty} \lambda_n\left\lVert g \right\rVert_{D^n(\mathcal{R}')}=\lim_{n \to \infty} \lambda_n\left\lVert g \right\rVert_{\mathcal{R}'}=0.$$ Hence every element of $A$ fixes a point in $T$. As $A$ is finitely generated, this implies that $A$ fixes a point in $T$ (see~\cite[I.6.5~Corollary~2]{Serre83}). Similarly, we see that $\left\langle w \right\rangle \ast \left\langle x \right\rangle$ fixes a point in $T$. As $W_n=A \ast \left\langle x \right\rangle$, we see that $A$ and $\left\langle w \right\rangle \ast \left\langle x \right\rangle$ cannot fix the same point in $T$. Let $U$ be the free splitting of $W_n$ associated with the free factor decomposition $W_n=A \ast \left\langle x \right\rangle$. Let $v_0$ be the vertex of $U$ fixed by $A$, let $v_1$ be the vertex fixed by $x$ and let $v_2$ be the vertex fixed by $wxw^{-1}$. Let $e_1$ be the edge between $v_0$ and $v_1$ and $e_2$ be the edge between $v_0$ and $v_2$. The arguments above show that we have a canonical $W_n$-equivariant morphism from $U$ to $T$. This morphism is obtained by a fold of the edges $e_1$ and $e_2$ of $U$ and this fold is extended $W_n$-equivariantly. Since $w$ is root-closed, there is no other edge of $U$ that can be folded as otherwise the stabilizer of an edge of $T$ would not be cyclic. Therefore the $\RR$-tree $T$ is simplicial and the decomposition of $W_n$ associated with $W_n \backslash T$ is $$W_n= A \ast_{\left\langle w \right\rangle} \left(\left\langle w \right\rangle \ast \left\langle x \right\rangle \right).$$ Hence $\mathcal{T}=\mathcal{S}'$ and the first statement follows.

Let us prove the second statement. For every $n \in \NN$, the equivalence classes $\lambda_nD^n(\mathcal{R})$ and $\mathcal{R}$ have compatible representatives. But as $\lim_{n \to \infty} \lambda_n D^n(\mathcal{R})=\mathcal{S}$, it follows from \cite[Corollary~A.12]{guirardel2016jsj} that, in the limit, the splittings $S$ and $R$ are compatible. 
\hfill\qedsymbol

\begin{lem}\label{Lem compatibility twists Wn-1}
Let $n \geq 3$ and let $\mathcal{S}$ be the equivalence class of a $W_{n-1}$-star $S$. Let $T$ be the group of twists of $\mathcal{S}$ and let $f \in T$ be an element of infinite order. Let $\mathcal{R}$ be the equivalence class of a $W_{n-1}$-star $R$ such that $f(\mathcal{R})=\mathcal{R}$. Then $S$ and $R$ are compatible.
\end{lem}

\dem Let $$W_n=A \ast \left\langle x_n\right\rangle$$ be a free factor decomposition of $W_n$ associated with $S$ and let $z_f \in A$ be the twistor of $f$. Let $z$ be a root-closed element of $A$ such that there exists $m \geq 1$ with $z^m=z_f$. Let $h \in T$ be the twist about $z$. We see that $h^m=f$. Let $S'$ be the splitting associated with the following amalgamated decomposition of $W_n$: $$W_n= A \ast_{\left\langle z \right\rangle} \left( \left\langle x_n \right\rangle \ast \left\langle z \right\rangle \right).$$ Let $\mathcal{S}'$ be the equivalence class of $S'$. Let $T'$ be the group of twists of $\mathcal{S}'$. Since $A$ is isomorphic to $W_{n-1}$ and since $z$ is root-closed, we see that $C_A(z)=\left\langle z \right\rangle$. Therefore $T'$ is isomorphic to $\ZZ$ and a generator of $T'$ is $h$. As $f(\mathcal{R})=\mathcal{R}$, Lemma~\ref{Lem compatibility splittings twists} implies that $S'$ and $R$ are compatible. Let $U$ be a common refinement of $S'$ and $R$ whose number of orbits of edges is minimal. Since both $S'$ and $R$ are one-edge splittings and are different, the splitting $U$ has $2$ orbits of edges. It follows that $W_n \backslash U$ is obtained from $W_n \backslash S'$ by blowing-up an edge at one of the two vertices of $W_n \backslash S'$. Let $\widetilde{v}$ be the vertex of $S'$ whose stabilizer is $A$ and let $v$ be its image in $W_n \backslash S'$. Let $\widetilde{w}$ be the vertex of $S'$ fixed by $\left\langle x_n \right\rangle \ast \left\langle z \right\rangle$ and let $w$ be its image in $W_n \backslash S'$.

\medskip

\noindent{\bf Claim. } Either $\mathcal{S}=\mathcal{R}$ or the splitting $W_n \backslash U$ is obtained from $W_n \backslash S'$ by blowing-up an edge at $v$.

\medskip

\dem Suppose that $W_n \backslash U$ is obtained from $W_n \backslash S'$ by blowing-up an edge at $w$. Then, since the group $G_w$ associated with $w$ is $\left\langle x_n \right\rangle \ast \left\langle z \right\rangle$ and since $z$ must fix an edge of $U$, we see that a free splitting of $G_w$ such that $z$ fixes a vertex is a $(G_w,\{\left\langle z \right\rangle,\left\langle x_n \right\rangle\})$-free splitting. But $(G_w,\{\left\langle z \right\rangle,\left\langle x_n \right\rangle\})$ has exactly one such equivalence class of one-edge free splitting: the one with vertex stabilizers conjugated with $\left\langle z\right\rangle$ and $\left\langle x_n \right\rangle$. This implies that $\mathcal{R}=\mathcal{S}$. The claim follows.
\hfill\qedsymbol

\medskip

Suppose that $\mathcal{R} \neq \mathcal{S}$. The claim implies that the amalgamated decomposition of $W_n$ associated with $U$ is $$W_n=B \ast C \ast_{\left\langle z \right\rangle} \left(\left\langle z \right\rangle \ast \left\langle x_n \right\rangle\right),$$ where $B$ and $C$ are free factors of $W_n$ such that $A=B \ast C$ and $z \in C$. Let $U'$ be a refinement of $U$ whose associated amalgamated decomposition of $W_n$ is: $$W_n=B \ast C \ast_{\left\langle z \right\rangle} \left\langle z \right\rangle \ast \left\langle x_n \right\rangle,$$ that is, $z$ and $x_n$ fix distinct points in $U'$. Then, since $A=B \ast C$, the splitting $U'$ is a refinement of $S$. This concludes the proof.
\hfill\qedsymbol

\begin{prop}\label{Prop commutator twist splitting fixed}
Let $n \geq 3$. Let $S$ be a $W_{n-1}$-star and let $f \in \Out(W_n)$ be a twist about the unique edge of $W_n \backslash S$. Let $g \in C_n$ be such that $g \in C_{Out(W_n)}(f)$. Then $g(\mathcal{S})=\mathcal{S}$.
\end{prop}  

\dem  Let $$W_n= \left\langle x_1,\ldots,x_{n-1} \right\rangle \ast \left\langle x_{n} \right\rangle$$ be the free factor decomposition associated with $S$ and let $\mathcal{S}$ be the equivalence class of $S$. By Lemma~\ref{Lem stab corank one Wn-1 star}, in order to prove that $g(\mathcal{S})=\mathcal{S}$, it suffices to show that $g$ preserves the conjugacy class of $A=\left\langle x_1,\ldots,x_{n-1} \right\rangle$. Let $\widetilde{f}$ be a representative of $f$ such that $\widetilde{f}|_{A}=\mathrm{id}_{A}$. Let $\widetilde{g}$ be a representative of $g$. Suppose towards a contradiction that $\widetilde{g}$ does not preserve the conjugacy class of $A$. By hypothesis, there exists $I \in \Inn(W_n)$ such that $\widetilde{f} \circ \widetilde{g}=I \circ \widetilde{g} \circ \widetilde{f}$. Thus, $$\widetilde{f} \circ \widetilde{g} (A)=I \circ \widetilde{g} \circ \widetilde{f}(A)=I \circ \widetilde{g}(A).$$
Therefore, $f$ preserves the conjugacy class of $\widetilde{g}(A)$. By Lemma~\ref{Lem stab corank one Wn-1 star}, $f$ fixes the unique equivalence class $\mathcal{R}$ of the \mbox{$W_{n-1}$-star} $R$ associated with $\widetilde{g}(A)$. By Lemma~\ref{Lem compatibility twists Wn-1}, the splittings $S$ and $R$ are compatible. Since we suppose that $\widetilde{g}(A) \notin [A]$, there exists a common refinement $S'$ of $S$ and $R$ which is a $W_{n-2}$-star.  Thus, there exists $y_n \in W_n$ such that the free factor decomposition associated with $S'$ is $$W_n=\left\langle x_n \right\rangle \ast B \ast \left\langle y_n \right\rangle,$$ where $B$ is such that $A=B \ast \left\langle y_n \right\rangle$ and $B \ast \left\langle x_n \right\rangle$ is a conjugate of $\widetilde{g}(A)$. Up to changing the representative $\widetilde{g}(A)$, we may suppose that $\widetilde{g}(A)=B \ast \left\langle x_n \right\rangle$. This implies that $x_n \in \widetilde{g}(A)$, that is $\widetilde{g}^{-1}(x_n) \in A$. But, since $A=\left\langle x_1,\ldots,x_{n-1} \right\rangle$, we see that \mbox{$[\widetilde{g}^{-1}(x_n)] \in \{[x_1],\ldots,[x_{n-1}]\}$}. This contradicts the fact that $g \in C_n$.
\hfill\qedsymbol

\bigskip

Combining Lemma~\ref{Lem compatibility twists Wn-1} and Proposition~\ref{Prop commutator twist splitting fixed}, we have the following corollary.

\begin{coro}\label{Coro compatibility commuting twists}
Let $n \geq 3$. Let $\mathcal{S}$ and $\mathcal{R}$ be two distinct $W_n$-equivariant homeomorphism classes of two $W_{n-1}$-stars $S$ and $R$. Let $f$ and $g$ be twists about respectively $S$ and $R$ such that $f$ and $g$ commute. Then $S$ and $R$ are compatible.
\end{coro}

\dem Let $k \geq 1$ be such that $g^k \in C_n$. By Proposition~\ref{Prop commutator twist splitting fixed}, since $g^k$ and $f$ commute, we have $g^k(\mathcal{S})=\mathcal{S}$. Since $g^k$ is a twist about $\mathcal{R}$, by Lemma~\ref{Lem compatibility twists Wn-1}, we have that $S$ and $R$ are compatible.
\hfill\qedsymbol

\bigskip

Let $\mathcal{S}$ be the equivalence class of a $W_{n-1}$-star $S$ and let $$W_n= \left\langle x_1,\ldots,x_{n-1} \right\rangle \ast \left\langle x_n \right\rangle$$ be the free factor decomposition of $W_n$ associated with $S$. Let $A=\left\langle x_1,\ldots,x_{n-1} \right\rangle$. Let $f \in \Stab_{\Out(W_n)}(\mathcal{S})$. Then any representative of $f$ sends $A$ to a conjugate of itself. Let $\widetilde{f}'$ be a representative of $f$ such that $\widetilde{f}'(A)=A$. Since the vertices in $S$ fixed by $A$ and $x_n$ are adjacent, and since the stabilizer of every vertex in $S$ adjacent to the vertex fixed by $A$ is a conjugate of $\left\langle x_n \right\rangle$ by an element of $A$, we see that $\widetilde{f}'(x_n)=xx_nx^{-1}$ with $x \in A$. Therefore, there exists a representative $\widetilde{f}$ of $f$ such that $\widetilde{f}(A)=A$ and $\widetilde{f}(x_n)=x_n$. The automorphism $\widetilde{f}$ is the unique representative of $f$ such that $\widetilde{f}(A)=A$ and $\widetilde{f}(x_n)=x_n$. 

We have a similar result for $W_{n-2}$-stars. Indeed, let $\mathcal{S}'$ be the equivalence class of a $W_{n-2}$-star $S'$ and let $$W_n= \left\langle x_1\right\rangle \ast \left\langle x_2,\ldots,x_{n-1} \right\rangle \ast \left\langle x_n \right\rangle$$ be the free factor decomposition of $W_n$ associated with $S'$ and let $B=\left\langle x_2,\ldots,x_{n-1} \right\rangle$. Let $f \in \Stab_{C_n}(\mathcal{S}')$. A similar argument as in the case of a $W_{n-1}$-star shows that there exists a representative $\widetilde{f}$ of $f$ such that $\widetilde{f}(B)=B$ and $\widetilde{f}(x_n)=x_n$.

\begin{lem}\label{Lem twists Wn-1 star in stab Wn-2}
Let $n \geq 4$. Let $\mathcal{S}$ be the $W_n$-equivariant homeomorphism class of a \mbox{$W_{n-1}$-star} $S$. Let $T$ be the group of twists of $\mathcal{S}$. Let $\mathcal{S}'$ be the $W_n$-equivariant homeomorphism class of a $W_{n-2}$-star $S'$ which refines $S$. Let $e$ be the edge of $W_n \backslash S'$ such that a representative of $\mathcal{S}$ is obtained from $W_n \backslash S'$ by collapsing the edge distinct from $e$. Let $T'$ be the group of twists of $S'$ about the edge $e$. Then $T \cap \Stab_{C_n}(\mathcal{S}') \subseteq T'$.
\end{lem}

\dem Let $$W_n=\left\langle x_1 \right\rangle \ast \left\langle x_2,\ldots,x_{n-1} \right\rangle \ast \left\langle x_n \right\rangle$$ be the free factor decomposition of $W_n$ induced by $S'$ and let $A=\left\langle x_2,\ldots,x_{n-1} \right\rangle$. Let $$W_n=B \ast \left\langle y_n \right\rangle$$ be the free factor decomposition associated with $S$. Up to changing the representative $S$, we may suppose that $B=\left\langle x_1,\ldots,x_{n-1} \right\rangle$ and that $y_n=x_n$. Let $f \in T \cap  \Stab_{C_n}(\mathcal{S}')$. Let $\widetilde{f}$ be the representative of $f$ such that $\widetilde{f}(B)=B$ and $\widetilde{f}(x_n)=x_n$ which exists since $f \in \Stab_{\Out(W_n)}(\mathcal{S})$. Since $f \in T$, there exists $g \in B$ such that $\widetilde{f}|_B$ is the global conjugation by $g$. Let $\widetilde{f}'$ be a representative of $f$ such that $\widetilde{f}'(A)=A$ and $\widetilde{f}'(x_n)=x_n$, which exists since $f \in \Stab_{C_n}(\mathcal{S}')$. Since the centralizer in $W_n$ of $x_n$ is $\left\langle x_n \right\rangle$ and since $A$ is malnormal in $W_n$, we see that $\widetilde{f}=\widetilde{f}'$. Hence $\widetilde{f}(A)=A$, and, since $A$ is malnormal, we see that $g \in A$. Therefore, $f \in T'$, which concludes the proof.
\hfill\qedsymbol

\bigskip

\begin{lem}\label{Lem twists commuting element twistor}
Let $n \geq 3$. Let $\mathcal{S}$ be the equivalence class of a $W_{n-1}$-star $S$ and let $$W_n= \left\langle x_1,\ldots,x_{n-1} \right\rangle \ast \left\langle x_n \right\rangle$$ be the free factor decomposition associated with $S$. Let $A=\left\langle x_1,\ldots,x_{n-1} \right\rangle$. Let $T$ be the group of twist of $S$. For $f \in T$, let $z_f \in A$ be the twistor of $f$. Let $g \in \Stab(\mathcal{S})$ and let $\widetilde{g}$ be a representative of $g$ such that $\widetilde{g}(A)=A$ and $\widetilde{g}(x_n)=x_n$. Then $g \in C_{\Out(W_n)}(\left\langle f \right\rangle)$ if and only if $\widetilde{g}(z_f)=z_f$.
\end{lem}

\dem By Proposition~\ref{Prop Levitt stab}~$(2)$, the group $\Stab(\mathcal{S})$ is isomorphic to $\Aut(A)$. The isomorphism $\Stab(\mathcal{S}) \to \Aut(A)$ is defined by sending $f \in \Stab(\mathcal{S})$ to its representative $\widetilde{f}$ such that $\widetilde{f}(A)=A$ and $\widetilde{f}(x_n)=x_n$. In particular, for every $h_1,h_2 \in \Out(W_n) \cap \Stab(\mathcal{S})$, we see that $h_1$ and $h_2$ commute if and only if there exist representatives $\widetilde{h}_1$ and $\widetilde{h}_2$ of $h_1$ and $h_2$ respectively such that $\widetilde{h}_1(A)=A$, $\widetilde{h}_2(A)=A$, $\widetilde{h}_1(x_n)=\widetilde{h}_2(x_n)=x_n$ and $\widetilde{h}_1 \circ \widetilde{h}_2=\widetilde{h}_2 \circ \widetilde{h}_1$. Moreover, Proposition~\ref{Prop Levitt stab}~$(2)$ identifies the group of twists $T$ with the group $\Inn(A)$. For $a \in A$, let $\ad_a$ be the inner autmorphism of $A$ induced by $a$. Since, for every $h \in \Aut(A)$ and every $a \in A$, we have $h \, \ad_a\,h^{-1}=\ad_{h(a)}$, we see that $h$ commutes with $\ad_a$ if and only if $h(a)=a$. Hence $g \in C_{\Out(W_n)}(\left\langle f \right\rangle)$ if and only if $\widetilde{g}(z_f)=z_f$.
\hfill\qedsymbol

\section{Direct products of nonabelian free groups in $\Out(W_n)$}\label{Section properties commuting normal subgroup}

Following \cite[Section~6]{HorbezWade20}, we define the \emph{product rank} of a group $H$, denoted by $\rk_{\pr}(H)$, to be the maximal integer $k$ such that a direct product of $k$ nonabelian free groups embeds in $H$. Note that, if $H'$ is a finite index subgroup of $H$, then $\rk_{\pr}(H')=\rk_{\pr}(H)$. Moreover, if $\phi \colon H \to \ZZ$ is a homomorphism, then $\rk_{\pr}(\ker(\phi))=\rk_{\pr}(H)$. The aim of this section is to prove the following theorem.

\begin{theo}\label{Theo product rank Wn}
\noindent{$(1)$ } For every $n \geq 3$, we have $\rk_{\pr}(\Aut(W_n))=n-2$. 

\medskip

\noindent{$(2)$ } For every $n \geq 4$, we have $\rk_{\pr}(\Out(W_n))=n-3$.

\medskip

\noindent{$(3)$ } Suppose that $n \geq 5$. If $H$ is a subgroup of $\Out(W_n)$ isomorphic to a direct product of $n-3$ nonabelian free groups, then $H$ has a subgroup $H'$ isomorphic to a direct product of $n-3$ nonabelian free groups which virtually fixes the $W_n$-equivariant homeomorphism class of a $W_{n-1}$-star. In addition, $H$ does not virtually fix the $W_n$-equivariant homeomorphism class of any one-edge free splitting that is not a $W_{n-1}$-star. 
\end{theo}

We first recall an estimate regarding product ranks and group extensions due to Horbez and Wade.

\begin{lem}\cite[Lemma 6.3]{HorbezWade20}\label{Lem extension product rank}
Let $1 \to N \to G \to Q \to 1$ be a short exact sequence of groups. Then $\rk_{\pr}(G) \leq \rk_{\pr}(N) + \rk_{\pr}(Q)$.
\hfill\qedsymbol
\end{lem}

In order to compute the product rank of $\Out(W_n)$, we take advantage of its action on the Gromov hyperbolic free factor complex. We recall a general result concerning actions of direct products on a hyperbolic space.

\begin{lem}\cite[Proposition~4.2, Lemma~4.4]{HorbezWade20}\label{Lem hyperbolic space and commuting}
Let $X$ be a Gromov hyperbolic space, and let $H$ be a group acting by isometries on $X$. Assume that $H$ contains a normal subgroup $K$ isomorphic to a direct product $K=\prod_{i=1}^k K_i$. 

If there exists $j \in \{1,\ldots,k\}$ such that $K_j$ contains a loxodromic element, then $\prod_{i \neq j} K_i$ has a finite orbit in $\partial_{\infty} X$.

If there exist two distinct $i,j \in \{1,\ldots,k\}$ such that both $K_i$ and $K_j$ contain a loxodromic element, then $H$ has a finite orbit in $\partial_{\infty} X$.

If, for every $j \in \{1,\ldots,k\}$, the group $K_j$ does not contain a loxodromic element, then either $K$ has a finite orbit in $\partial_{\infty} X$ or $H$ has bounded orbits in $X$.
\hfill\qedsymbol
\end{lem}

We will also use a theorem due to Guirardel and Horbez which assigns to every nonempty collection of free splittings whose elementwise stabilizer is infinite a canonical (not necessarily free) splitting.

\begin{theo}\cite{Guirardelhorbez}\label{Theo invariant splitting constructed}
Let $n \geq 3$. There exists an $\Out(W_n)$-equivariant map which assigns to every nonempty collection $\mathcal{C}$ of free splittings of $W_n$ whose elementwise $\Out(W_n)$-stabilizer is infinite, a nontrivial splitting $U_{\mathcal{C}}$ of $W_n$ whose set of vertices $VU_{\mathcal{C}}$ has a $W_n$-invariant partition $VU_{\mathcal{C}}=V_1 \amalg V_2$ with the following properties:
\begin{enumerate}
\item For every vertex $v \in V_1$, the following holds:
\begin{enumerate}
\item either some edge incident on $v$ has trivial stabilizer, or the set of stabilizers of edges incident on $v$ induces a nontrivial free factor system of the vertex stabilizer $G_v$,
\item there exists a finite index subgroup $H_0$ of the elementwise stabilizer of the collection $\mathcal{C}$ such that every outer automorphism in $H_0$ has a representative in $\Aut(W_n)$ which restricts to the identity on $G_v$.
\end{enumerate}
\item The collection of all conjugacy classes of stabilizers of vertices in $V_2$ is a free factor system of $W_n$.
\hfill\qedsymbol
\end{enumerate}
\end{theo}

\noindent{\bf Proof of Theorem~\ref{Theo product rank Wn}. } The proof is inspired by \cite[Theorem~6.1]{HorbezWade20} due to Horbez and Wade and \cite[Theorem~4.3]{HenselHorbezWade19} due to Hensel, Horbez and Wade. 

We first prove that if $n \geq 4$, then $\rk_{\pr}(\Out(W_n)) \geq n-3$ and that, if $n \geq 3$, then $\rk_{\pr}(\Aut(W_n)) \geq n-2$. Pick a standard generating set $\{x_1,\ldots,x_n\}$ of $W_n$. Then the group $H$ generated by $\{x_1x_2,x_2x_3\}$ is a nonabelian free group (see~\cite[Theorem~A]{muhlherr1997}). 

Suppose first that $n \geq 4$. For $i \in \{4,\ldots,n\}$ and $h \in H$, let $F_{i,h}$ be the automorphism sending $x_i$ to $hx_ih^{-1}$ and, for $j \neq i$, fixing $x_j$. Then, for every distinct $i,j \in \{4,\ldots,n\}$ and for every $g,h \in H$, the automorphisms $F_{i,g}$ and $F_{j,h}$ commute, giving a direct product of $n-3$ nonabelian free groups in $\Out(W_n)$. Moreover, for every $g,h \in H$, and every $i \in \{4,\ldots,n\}$, the inner automorphism $\ad_g$ commutes with $F_{i,h}$, which yields a direct product of $n-2$ nonabelian free groups in $\Aut(W_n)$. In the case where $n=3$, the group $\Aut(W_3)$ contains the subgroup  $\left\langle \ad_h \right\rangle_{h \in  H}$, which is a nonabelian free group.

We now prove that, if $n \geq 3$, then $\rk_{\pr}(\Aut(W_n)) \leq n-2$, if $n=3$, then $\rk_{\pr}(\Out(W_n))=1$ and if $n \geq 4$, then $\rk_{\pr}(\Out(W_n)) \leq n-3$. The proof is by induction on $n$. The base case where $n=3$ follows from the fact that the group $\Aut(W_3)$ is isomorphic to $\Aut(F_2)$ (see \cite[Lemma~2.3]{varghese2019}) and the fact that the group $\Aut(F_2)$ does not contain a direct product of two nonabelian free groups (see \cite[Lemma~6.2]{HorbezWade20}). Moreover, by \cite[Proposition~2.2]{Guerch2020out}, the group $\Out(W_3)$ is isomorphic to $\mathrm{PGL}(2,\ZZ)$ which is virtually free.

Let $k \geq \max\{n-3,2\}$ and let $H=H_1 \times H_1 \times \ldots \times H_k$ be a subgroup of $\Out(W_n)$ isomorphic to a direct product of $k$ nonabelian free groups. Note that $k=n-3$ if $n \geq 5$ and $k=2$ if $n=4$. We prove that there exists a subgroup $K$ of $H$ isomorphic to a direct product of $k$ nonabelian free groups which virtually fixes a one-edge free splitting of $W_n$. Let $\mathcal{F}$ be a maximal $H$-periodic free factor system. If $\mathcal{F}$ is sporadic, then $H$ virtually fixes a one-edge free splitting, so we are done. Therefore, we may suppose that $\mathcal{F}$ is nonsporadic. As $\mathcal{F}$ is supposed to be maximal, by Proposition~\ref{Prop maximal free factor system bounded orbit}, the group $H$ acts on $\FF(W_n,\mathcal{F})$ with unbounded orbits. Lemma~\ref{Lem hyperbolic space and commuting} implies that, after possibly reordering the factors, the group $H'=H_1 \times H_2 \times \ldots \times H_{k-1}$ has a finite orbit in $\partial_{\infty} \FF(W_n, \mathcal{F})$. By Lemma~\ref{Lem fixed point at infinity FF}, the group $H'$ virtually fixes the homothety class $[T]$ of an arational $(W_n,\mathcal{F})$-tree $T$.

Let $H_0$ be a normal subgroup of finite index in $H'$ that is contained in $\Stab([T])$. 

\medskip

\noindent{\bf Claim. } The group $H$ contains a subgroup isomorphic to a direct product of $k$ nonabelian free groups, which fixes the equivalence class of a one-edge free splitting. 

\medskip

\dem By Lemma~\ref{Lem isometric stabilizer cyclic}, the morphism $\lambda_T|_{H_0}$ from $H_0$ to $\RR_+^*$ given by the scaling factor has cyclic image. 
As $H_0$ contains a direct product of $k-1$ nonabelian free groups, so does $P=\ker(\lambda_T|_{H_0})$ (see the beginning of Section~\ref{Section properties commuting normal subgroup}). As $P$ is contained in the isometric stabilizer of $T$, Proposition~\ref{Prop isometric stabilizer free splitting fixed} implies that $P$ contains a finite index subgroup $P_0$ which fixes infinitely many $(W_n,\mathcal{F})$-free splittings.

Let $\mathcal{C}$ be the (nonempty) collection of all $(W_n,\mathcal{F})$-free splittings fixed by the infinite group $P_0$, let $U_{\mathcal{C}}$ be the splitting provided by Theorem~\ref{Theo invariant splitting constructed}, and let $\mathcal{U}_{\mathcal{C}}$ be its equivalence class. Since $P_0$ commutes with $H_k$, the equivalence class $\mathcal{U}_{\mathcal{C}}$ is $\left(P_0 \times H_k\right)$-invariant. 

Suppose first that the splitting $U_{\mathcal{C}}$ contains an edge $e \in EU_{\mathcal{C}}$ with trivial stabilizer. Let $U'$ be the splitting obtained from $U_{\mathcal{C}}$ by collapsing every edge of $U_{\mathcal{C}}$ that is not contained in the orbit of $e$, and let $\mathcal{U}'$ be its equivalence class. Then $\mathcal{U}'$ is the equivalence class of a one-edge free splitting virtually fixed by $P_0 \times H_k$. Since $P_0$ contains a direct product of $k-1$ nonabelian free groups, the claim follows.

Thus, we can suppose that all edge stabilizers of $U_{\mathcal{C}}$ are nontrivial. We show that this leads to a contradiction. Let $VU_{\mathcal{C}}=V_1 \amalg V_2$ be the partition of $VU_{\mathcal{C}}$ given by Theorem~\ref{Theo invariant splitting constructed}. Let $P'$ be a finite index subgroup of $P_0$ which acts trivially on the quotient $W_n \backslash U_{\mathcal{C}}$. We claim that the intersection of $P'$ with the group of twists of $U_{\mathcal{C}}$ is trivial. Indeed, let $e$ be an half-edge of $U_{\mathcal{C}}$. As every subgroup of $W_n$ with nontrivial centralizer is cyclic, if the edge stabilizer $G_e$ of $e$ is not cyclic, the group of twists around this half-edge is trivial. Thus, half-edges with nontrivial group of twists have cyclic stabilizers. But twists about edges with cyclic stabilizers are central in a finite index subgroup of $\Stab^0(U_{\mathcal{C}})$ by Lemma~\ref{Lem twists about cyclic edge central}. As the center of every finite index subgroup of $P'$ is trivial, we see that the intersection of $P'$ with the group of twists is trivial. By Remark~\ref{Remark proposition Levitt}, up to passing to a further finite index subgroup of $P'$, we may suppose that the intersection of $P'$ with the group of bitwists is trivial.

By Proposition~\ref{Prop Levitt stab}~$(1)$ and Remark~\ref{Remark proposition Levitt}, the fact that the intersection of $P'$ with the group of bitwists is trivial implies that we have an injective homomorphism $$P' \to \prod_{v \in W_n \backslash VU_{\mathcal{C}}} \Out(G_v).$$ By Theorem~\ref{Theo invariant splitting constructed}~$(1)(b)$, for every vertex $v \in V_1$, the homomorphism $P' \to \Out(G_v)$ has finite image. Therefore, up to passing to a finite index subgroup of $P'$, we have an injective map $$P' \to \prod_{v \in W_n \backslash V_2} \Out(G_v).$$ By Theorem~\ref{Theo invariant splitting constructed}~$(2)$, for every $v \in V_2$, the vertex stabilizer $G_v$ is an element of a free factor system of $W_n$. Therefore, there exists $k$ such that $G_v$ is isomorphic to $W_k$. By Lemma~\ref{Lem extension product rank}, we have: $$n-4 \leq k-1=\rk_{\pr}(P') \leq \sum_{v \in W_n \backslash V_2} \rk_{\pr} (\Out(G_v)).$$ By induction,  we see that, if $|W_n \backslash V_2| \geq 2$, then $$\sum_{v \in W_n \backslash V_2} \rk_{\pr} (\Out(G_v)) \leq n-6,$$ which leads to a contradiction. Thus $|W_n \backslash V_2|=1$. Let $v \in W_n \backslash V_2$. Then there exists $\ell \in \{1,\ldots,n-1\}$ such that $G_v$ is isomorphic to $W_{\ell}$. If $\ell \leq n-2$, then $$\rk_{\pr} (\Out(G_v)) \leq n-5,$$ which leads to a contradiction. If $\ell=n-1$, then the free factor system $\mathcal{F}$ contains a free factor isomorphic to $W_{n-1}$ and is therefore a sporadic free factor system, which leads to a contradiction.
\hfill\qedsymbol

\medskip

Therefore, we see that there exists a subgroup $K$ of $H$ isomorphic to a direct product of $k$ nonabelian free groups such that $K$ fixes the $W_n$-equivariant homeomorphism class of a one-edge-free splitting $\mathcal{S}$. We now prove that $\mathcal{S}$ is the equivalence class of a $W_{n-1}$-star. Let $S$ be a representative of $\mathcal{S}$, let $v_1$ and $v_2$ be the vertices of the underlying graph of $W_n \backslash S$ and, for $i \in \{1,2\}$, let $k_i$ be such that $W_{k_i}$ is isomorphic to $G_{v_i}$. Let $K_0$ be the finite index subgroup of $K$ which acts as the identity on $W_n \backslash S$. Then $K_0 \subseteq \Stab^0(\mathcal{S})$. By Proposition~\ref{Prop Levitt stab}~$(2)$, the group $\Stab^0(\mathcal{S})$ is isomorphic to $\Aut(W_{k_1}) \times \Aut(W_{k_2})$. Suppose towards a contradiction that, for every $i \in \{1,2\}$, we have that $k_i \neq 1$. Suppose first that, for every $i \in \{1,2\}$, we have $k_i \geq 3$. Then, by Lemma~\ref{Lem extension product rank}, we see that: $$k=\rk_{\pr}(K_0) \leq \rk_{\pr}(\Aut(W_{k_1}))+\rk_{\pr}(\Aut(W_{k_2})) \leq k_1-2 +k_2-2=n-4,$$ where the second inequality comes from the induction hypothesis. If there exists $i \in \{1,2\}$ such that $k_i=2$, then, as $\Aut(W_2)$ is virtually cyclic (it is isomorphic to $W_2$ by \cite[Lemma~1.4.2]{thomasautotower}), we see that: $$k=\rk_{\pr}(K_0) \leq \rk_{\pr}(\Aut(W_{k_1}))+\rk_{\pr}(\Aut(W_{k_2})) \leq k_1-2\leq n-4.$$ In both cases, we have a contradiction as $k \geq n-3$ when $k \geq 5$ and $k=n-2$ when $n=4$. Thus, there exists $i \in \{1,2\}$ such that $k_i=1$. This shows that $S$ is a $W_{n-1}$-star. In particular, when $k=n-3$, that is, when $n \geq 5$, this proves Theorem~\ref{Theo product rank Wn}~$(3)$.

Since $K_0 \subseteq \Stab^0(\mathcal{S})$, Proposition~\ref{Prop Levitt stab}~$(2)$ implies that $$k=\rk_{\pr}(K_0) \leq \rk_{\pr}(\Aut(W_{n-1}))= n-1-2= n-3.$$ When $n=4$, then $k=2=n-2$. Therefore, we have a contradiction in this case. This shows that, for all $n \geq 4$, the product rank of $\Out(W_n)$ is equal to $n-3$. This concludes the proof of Theorem~\ref{Theo product rank Wn}~$(2)$.

It remains to prove that, if $n \geq 4$, we have $\rk_{\pr}(\Aut(W_n)) \leq n-2$. We have the following short exact sequence $$1 \to W_n \to \Aut(W_n) \to \Out(W_n) \to 1.$$ By Lemma~\ref{Lem extension product rank}, as $W_n$ is virtually free, we see that $$\rk_{\pr}(\Aut(W_n)) \leq \rk_{\pr}(W_n) + \rk_{\pr}(\Out(W_n))=1+n-3=n-2.$$ This concludes the proof of Theorem~\ref{Theo product rank Wn}~$(1)$. 
\hfill\qedsymbol

\bigskip

\section{Subgroups of stabilizers of $W_{n-1}$-stars}\label{Section subgroups Wn-1}

In the next two sections, we prove an algebraic characterisation of  stabilizers of equivalence classes of $W_{n-2}$-stars. In this section, we take advantage of properties satisfied by stabilizers of equivalence classes of $W_{n-2}$-stars which are sufficiently rigid to show that a subgroup $H$ of $\Out(W_n)$ which satisfies these properties virtually fixes a $W_{n-1}$-star. In the next section, we will take advantage of the fact that stabilizers of equivalence classes of compatible $W_{n-2}$-stars have large intersections to give a characterisation of stabilizers of equivalence classes of $W_{n-2}$-stars.

\medskip

Let $\Gamma$ be a finite index subgroup of the group $C_n$ (defined after Theorem~\ref{Theo generating set Aut}). We introduce the following algebraic property for a subgroup $H \subseteq \Gamma$. 

\bigskip

\noindent $(P_{W_{n-2}})$ The group $H$ satisfies the following three properties:

\begin{enumerate}
\item The group $H$ contains a normal subgroup isomorphic to a direct product $K_1 \times K_2$ of two normal subgroups such that each one contains a nonabelian finitely generated normal free subgroups of finite index and such that for every $i \in \{1,2\}$, for every nontrivial normal subgroup $P$ of a finite index subgroup $K_i'$ of $K_i$, and for every finite index subgroup $P'$ of $P$, the group $C_{C_n}(P')$ contains $K_{i+1}$ as a finite index subgroup (where indices are taken modulo $2$).

\item The group $H$ contains a direct product of $n-3$ nonabelian free groups.

\item The group $H$ contains a subgroup isomorphic to $\ZZ^{n-2}$.

\end{enumerate}

\begin{rmq}\label{Rmq centralizer PWn-2}

\noindent{$(1)$ } Notice that property~$(P_{W_{n-2}})$ is closed under taking finite index subgroups.

\noindent{$(2)$ } Hypothesis $(P_{W_{n-2}})~(1)$ implies that, if for every $i \in \{1,2\}$, the group $P_i$ is a finite index subgroup of a nontrivial normal subgroup of a finite index subgroup of $K_i$, the centralizer in $C_n$ of $P_1 \times P_2$ is finite.

\end{rmq}

We first prove that the stabilizer in $\Gamma$ of the equivalence class of a $W_{n-2}$-star satisfies $(P_{W_{n-2}})$. We then show that a group satisfying $(P_{W_{n-2}})$ virtually fixes the equivalence class of a $W_{n-1}$-star.

\subsection{Properties of $\mathcal{Z}_{RC}$-factors}

In order to prove that the stabilizer in $\Gamma$ of the equivalence class of a $W_{n-2}$-star satisfies $(P_{W_{n-2}})$, we first need some background concerning $\mathcal{Z}_{RC}$-splittings. Let $G$ be a finitely generated group. A \emph{$\mathcal{Z}_{RC}$-splitting} of $G$ is a splitting of $G$ such that every edge stabilizer is either trivial or isomorphic to $\ZZ$ and root-closed. A \emph{$\mathcal{Z}_{RC}$-factor} of $G$ is a subgroup of $G$ which arises as a vertex stabilizer of a $\mathcal{Z}_{RC}$-splitting of $G$. Note that since edge stabilizers are root-closed, so are the vertex stabilizers. We outline here some properties of {$\mathcal{Z}_{RC}$-factors (see e.g.~\cite[Proposition~7.3]{HorbezWade20}).

\begin{prop}\label{Prop Zrc factors}
Let $n \geq 3$. The $\mathcal{Z}_{RC}$-factors of $W_n$ satisfy the following properties.

\medskip

\noindent{$(1)$ } Let $H$ be a finitely generated subgroup of $W_n$ which is not virtually cyclic. There exists $g \in H$ which is not contained in any proper $\mathcal{Z}_{RC}$-factor of $H$.

\medskip

\noindent{$(2)$ } There exists $C \in \NN^*$ such that, for every strictly ascending chain $G_1 \subsetneq \ldots \subsetneq G_k$ of $\mathcal{Z}_{RC}$-factors of $W_n$, one has $k \leq C$. 

\medskip

\noindent{$(3)$ } If a subgroup $K$ of $W_n$ is not contained in any proper $\mathcal{Z}_{RC}$-factor of $W_n$ and if $P$ is either a finite index subgroup of $K$ or a nontrivial normal subgroup of $K$, then $P$ is not contained in any proper $\mathcal{Z}_{RC}$-factor of $W_n$.

\medskip

\noindent{$(4)$ } A subgroup $K$ of $W_n$ is contained in a proper $\mathcal{Z}_{RC}$-factor of $W_n$ if and only if every element of $K$ is contained in a proper $\mathcal{Z}_{RC}$-factor of $W_n$.
\end{prop} 

\dem The first assertion is a consequence of \cite[Lemma~4.3]{GenevoisHorbez2020} due to Genevois and Horbez. 

For the second assertion, let $G_1 \subsetneq \ldots \subsetneq G_k$ be a sequence of strictly ascending $\mathcal{Z}_{RC}$-factors. Then, since $\mathcal{Z}_{RC}$-factors are root-closed, for every $i \geq 3$ the group $G_i$ is not cyclic. Thus, as we want an upper bound on the number of subgroups of such a sequence, we may suppose that for every $i \in \{1,\ldots,n\}$, the group $G_i$ is not cyclic. We claim that, for every $i \in \{1,\ldots,k\}$, there exists $\phi_i \in \Aut(W_n)$ such that $\mathrm{Fix}(\phi_i)=G_i$. Indeed, let $S_i$ be a $\mathcal{Z}_{RC}$-splitting of $W_n$ such that there exists $v \in VS_i$ whose stabilizer is equal to $G_i$. Up to collapsing edges, we may suppose that every vertex of $S_i$ has nontrivial stabilizer. Let $e_1,\ldots,e_{\ell}$ be the edges with origin $v$. Let $F \subseteq \{e_1,\ldots,e_{\ell}\}$ be the subset made of all edges with nontrivial stabilizer. By the definition of a $\mathcal{Z}_{RC}$-splitting, for every $e_s \in F$, the group $G_{e_s}$ is cyclic. For every $e_s \in F$, let $z_s$ be a generator of $G_{e_s}$. For every $e_{s'} \in \{e_1,\ldots,e_{\ell}\}-F$, let $z_{s'} \in G_i$ be such that, if $w_{s'}$ is the endpoint of $e_{s'}$ distinct from $v$, we have $z_{s'}G_{w_{s'}}z_{s'}^{-1} \neq G_{w_{s'}}$. Let $\phi_i=D_{e_1,z_1} \circ \ldots \circ D_{e_{\ell},z_{\ell}}$ be a multitwist about every edge with origin $v$. Then, as the centralizer of an infinite cyclic subgroup of $W_n$ is infinite cyclic, we have $\mathrm{Fix}(\phi_i)=G_i$. Therefore, in order to prove the second assertion, it suffices to prove that there exists $C \in \NN^*$ such that for every strictly ascending chain $\mathrm{Fix}(\phi_1) \subsetneq \ldots \subsetneq \mathrm{Fix}(\phi_k)$ of fixed points sets of automorphisms of $W_n$, one has $k \leq C$.

Let $\{x_1,\ldots,x_n\}$ be a standard generating set of $W_n$. By \cite[Theorem~A]{muhlherr1997} the kernel $K'$ of the homomorphism $W_n \to F$ which, for every $i \in \{1,\ldots,n\}$, sends $x_i$ to the nontrivial element of $F$ is a nonabelian free group on $n-1$ generators. Remark that $K'$ does not depend on the choice of the basis since, for every element $x$ of order $2$, there exists $i \in \{1,\ldots,n\}$ and $g \in W_n$ such that $x=gx_ig^{-1}$. Moreover, $K'$ is a characteristic subgroup of index $2$ of $W_n$ and the natural homomorphism $\Phi \colon \Aut(W_n) \to \Aut(K')$ given by restriction is injective. Then $\mathrm{Fix}(\Phi(\phi_1)) \subseteq \ldots \subseteq \mathrm{Fix}(\Phi(\phi_k))$ is an ascending chain of fixed points sets. 

\medskip

\noindent{\bf Claim. } For every $i \in \{2,\ldots,k-1\}$, the set $\{\mathrm{Fix}(\Phi(\phi_{i-1})),\mathrm{Fix}(\Phi(\phi_i)),\mathrm{Fix}(\Phi(\phi_{i+1}))\}$ contains at least $2$ elements. 

\medskip

\dem Suppose towards a contradiction that  $$|\{\mathrm{Fix}(\Phi(\phi_{i-1})),\mathrm{Fix}(\Phi(\phi_i)),\mathrm{Fix}(\Phi(\phi_{i+1}))\}|=1.$$ As $\mathrm{Fix}(\phi_{i-1}) \subsetneq \mathrm{Fix}(\phi_i)$ and $\mathrm{Fix}(\Phi(\phi_{i-1}))=\mathrm{Fix}(\Phi(\phi_i))$, there exists $a \in W_n-K'$ such that $\phi_i(a)=a$ and $\phi_{i-1}(a) \neq a$. Since the index of $K'$ is equal to $2$, we see that $\phi_{i-1}(a^2)=a^2$. Therefore, $\phi_{i-1}(a)^2=a^2$ and $\phi_{i-1}(a)$ is a square root of $a^2$. If $a^2$ has infinite order, its only square root is $a$. This implies that $\phi_{i-1}(a)=a$, a contradiction. Thus we can assume that $a$ has order $2$ and, up to changing the basis $\{x_1,\ldots,x_n\}$, we may suppose that $a=x_1$. 

As the index of $K'$ is equal to $2$, we have $W_n=K' \amalg x_1K'$. Let $x \in \mathrm{Fix}(\phi_{i+1})-K'$. Then there exists $y \in K'$ such that $x=x_1y$. As $x_1 \in \mathrm{Fix}(\phi_i)$ and $\mathrm{Fix}(\phi_i) \subsetneq \mathrm{Fix}(\phi_{i+1})$, we have that $\phi_{i+1}(x_1)=x_1$. Hence $\phi_{i+1}(y)=y$. As $y \in K'$ and $\mathrm{Fix}(\Phi(\phi_i))=\mathrm{Fix}(\Phi(\phi_{i+1}))$, we see that $\phi_i(y)=y$ and $\phi_i(x)=\phi_i(x_1y)=x_1y=x$. Therefore we have that $\mathrm{Fix}(\phi_i)=\mathrm{Fix}(\phi_{i+1})$, which is a contradiction. The claim follows.
\hfill\qedsymbol

\medskip

The claim implies that the length of the strictly ascending chain associated with $\mathrm{Fix}(\Phi(\phi_1)) \subseteq \ldots \subseteq \mathrm{Fix}(\Phi(\phi_k))$ is at least equal to $\frac{k}{2}$. But any strictly ascending chain of fixed subgroups in a free group on $n-1$ generators as length at most $2(n-1)$ (see~\cite[Theorem~4.1]{MartinoVentura2004}). Therefore, there exists $C$ which depends only on $n$ such that $k \leq C$. The second assertion of Proposition~\ref{Prop Zrc factors} follows.

We now prove the third assertion. Let $P$ and $K$ be as in Proposition~\ref{Prop Zrc factors}~$(3)$. If $K$ is a virtually infinite cyclic group, then $K$ is either isomorphic to $\ZZ$ or to $W_2$. Let $a$ be a generator of the subgroup of $K$ isomorphic to $\ZZ$ and root-closed in $K$. Since $\left\langle a\right\rangle$ is a finite index subgroup of $K$ and since $K$ is not contained in any proper $\mathcal{Z}_{RC}$-factor of $W_n$, then neither is $a$. Remark that any nontrivial normal subgroup of $K$ intersects the subgroup $\left\langle a \right\rangle$ non trivially. Therefore, if $P$ is contained in a proper $\mathcal{Z}_{RC}$-factor of $W_n$, then $a$ is elliptic in a $\mathcal{Z}_{RC}$-splitting. This contradicts the fact that $a$ is not contained in any proper $\mathcal{Z}_{RC}$-factor of $W_n$. 

So we can assume that $K$ is not virtually cyclic. As every finite index subgroup contains a nontrivial normal subgroup of $K$, we may assume that $P$ is a nontrivial normal subgroup of $K$. Notice that $P$ is necessarily noncyclic. Suppose towards a contradiction that $P$ is contained in a $\mathcal{Z}_{RC}$-factor. Then there exists a $\mathcal{Z}_{RC}$-splitting $S$ of $W_n$ such that $P$ is elliptic in $S$. Since edge stabilizers are cyclic, the group $P$ fixes a unique vertex $x$ of $S$. But, as $P$ is normal in $K$, for every $k \in K$, we have that $kx$ is also fixed by $P$, hence we have $kx=x$. Therefore, $x$ is fixed by $K$, which contradicts the fact that $K$ is not contained in any proper $\mathcal{Z}_{RC}$-factor.

We finally prove Proposition~\ref{Prop Zrc factors}~$(4)$. Suppose that $K$ is contained in a proper $\mathcal{Z}_{RC}$-factor. Then it is clear that every element of $K$ is contained in a proper $\mathcal{Z}_{RC}$-factor.

Conversely, assume that $K$ is not contained in any proper $\mathcal{Z}_{RC}$-factor of $W_n$. Let us prove that there exists $g \in K$ such that $g$ is not contained in any proper $\mathcal{Z}_{RC}$-factor. By Proposition~\ref{Prop Zrc factors}~$(2)$, there exists a bound on the length of an increasing chain of $\mathcal{Z}_{RC}$-factors of $W_n$. Therefore, the group $K$ contains a finitely generated subgroup $K'$ which is not contained in any proper $\mathcal{Z}_{RC}$-factor. By Proposition~\ref{Prop Zrc factors}~$(1)$, there exists $g \in K'$ such that $g$ is not contained in a proper $\mathcal{Z}_{RC}$-factor of $K'$. Let $S$ be a $\mathcal{Z}_{RC}$-splitting of $W_n$. As $K'$ is not contained in any proper $\mathcal{Z}_{RC}$-factor of $W_n$, the group $K'$ has a well-defined, nontrivial minimal subtree $S_{K'}$ with respect to the action of $K'$ on $S$. As $S$ is a $\mathcal{Z}_{RC}$-splitting of $W_n$, the splitting $S_{K'}$ is a $\mathcal{Z}_{RC}$-splitting of $K'$. Since $g$ is not contained in any proper $\mathcal{Z}_{RC}$-factor of $K'$, it follows that $g$ is a hyperbolic isometry of $S_{K'}$ and is not elliptic in $S$. As $S$ is arbitrary, it follows that $g$ is not contained in any $\mathcal{Z}_{RC}$-factor of $W_n$.
\hfill\qedsymbol

\bigskip

Proper $\mathcal{Z}_{RC}$-factors appear naturally when studying stabilizers of conjugacy classes of elements as shown by the following theorem. Recall that, if $\mathcal{H}=\{H_1,\ldots,H_k\}$ is a finite family of finitely generated subgorups of $W_n$, the group $\Out(W_n,\mathcal{H}^{(t)})$ is the subgroup of $\Out(W_n)$ consisting of all outer automorphisms $\phi \in \Out(W_n)$ such that, for every \mbox{$i \in \{1,\ldots,k\}$}, there exists a representative $\widetilde{\phi}_i \in \Aut(W_n)$ of $\phi$ such that $\widetilde{\phi}_i(H_i)=H_i$ and $\widetilde{\phi}_i|_{H_i}=\mathrm{id}_{H_i}$.

\begin{theo}\cite[Theorem~7.14]{GuirardelLevitt2015}\label{Theo Zrc factors conjugacy class fixed}
Let $n \geq 3$ and let $g \in W_n$. Then the subgroup $\Out(W_n,\left\langle g \right\rangle)$ of outer automorphisms which preserve $\left\langle g \right\rangle$ up to conjugacy is infinite if and only if $g$ is contained in a proper $\mathcal{Z}_{RC}$-factor of $W_n$.

More generally, if $\mathcal{H}$ is a finite family of finitely generated subgroups, then the group $\Out(W_n,\mathcal{H}^{(t)})$ is infinite if and only if there exists a nontrivial $\mathcal{Z}_{RC}$-splitting $S$ of $W_n$ such that every subgroup of $\mathcal{H}$ fixes a vertex of $S$.
\hfill\qedsymbol
\end{theo}

\subsection{Stabilizers of $W_{n-2}$-stars satisfy $(P_{W_{n-2}})$}

\begin{lem}\label{Lem normal subgroup of twist Wn-2 fixes 2 Wn-1}
Let $n \geq 5$ and let $\Gamma$ be a finite index subgroup of $C_n$. Let $\mathcal{S}$ be the equivalence class of a $W_{n-2}$-star $S$. Let $e_1$ and $e_2$ be the two edges of $W_n \backslash S$ and, for $i \in \{1,2\}$, let $T_i'$ be the group of twists about $e_i$ in $\Stab_{\Gamma}(\mathcal{S})$. Let $i \in \{1,2\}$, let $T_i$ be a finite index subgroup of $T_i'$ and let $P'$ be a finite index subgroup of a nontrivial normal subgroup of $T_i$. Then for every finite index subgroup $P_0$ of $P'$, the group $P_0$ fixes exactly one equivalence class of $W_{n-2}$-stars.
\end{lem}

\dem Let $$W_n= \left\langle x_1 \right\rangle \ast \left\langle x_3,\ldots,x_n\right\rangle \ast \left\langle x_2 \right\rangle$$ be a free factor decomposition associated with $W_n \backslash S$ and $A=\left\langle x_3,\ldots,x_n\right\rangle$. Up to exchanging the roles of $e_1$ and $e_2$, we may suppose that $P'$ is contained in the group of twists of the equivalence class of the $W_{n-1}$-star $S_1$ whose associated free factor decomposition of $W_n$ is, up to global conjugation: 
$$W_n=\left\langle x_1 \right\rangle \ast \left\langle x_2,x_3,\ldots,x_n\right\rangle.$$ Let $B=\left\langle x_2,x_3,\ldots,x_n\right\rangle$ and let $\mathcal{S}_1$ be the equivalence class of $S_1$. Finally, let $\mathcal{S}_2$ be the equivalence class of the $W_{n-1}$-star $S_2$ whose associated free factor decomposition of $W_n$ is, up to global conjugation: 
$$W_n=\left\langle x_2 \right\rangle \ast \left\langle x_1,x_3,\ldots,x_n\right\rangle.$$ Let $C=\left\langle x_1,x_3,\ldots,x_n\right\rangle= A \ast \left \langle x_1 \right\rangle$.

We claim that the only equivalence classes of $W_{n-1}$-stars fixed by any finite index subgroup of $P'$ are $\mathcal{S}_1$ and $\mathcal{S}_2$. Indeed, fix $i \in \{1,2\}$. The group $T_i$ is isomorphic to a finite index subgroup $N$ of $W_{n-2}$. By Proposition~\ref{Prop Zrc factors}~$(3)$ applied with $K=W_{n-2}$ and $P=N$, as $n \geq 5$, the group $N$ is not contained in any proper $\mathcal{Z}_{RC}$-free factor of $W_{n-2}$. By Proposition~\ref{Prop Zrc factors}~$(4)$, there exists $g \in N$ such that $W_{n-2}$ is freely indecomposable relative to $g$. Hence there exists $g \in A$ such that $A$ is freely indecomposable relative to $g$ and $P'$ contains the twist about $e_1$ whose twistor is $g$. Note that this twist can be seen as a twist about the $W_{n-1}$-star $S_1$. Let $\mathcal{S}_1'$ be the equivalence class of the one-edge cyclic splitting $S_1'$ whose associated amalgamated decomposition of $W_n$ is, up to global conjugation: 
$$W_n=\left(\left\langle x_1 \right\rangle \ast \left\langle g \right\rangle\right) \ast_{\left\langle g \right\rangle} B.$$

Let $\mathcal{S}_3$ be the equivalence class of a $W_{n-1}$-star $S_3$ fixed by some finite index subgroup of $P'$ and distinct from $\mathcal{S}_1$. Let $$W_n=\left\langle y \right\rangle \ast D$$ be the free factor decomposition associated with $S_3$. We claim that $\mathcal{S}_3=\mathcal{S}_2$. As $P'$ contains the twist about $g$, by Lemma~\ref{Lem compatibility twists Wn-1}, the splitting $S_3$ is compatible with $S_1'$. Let $U$ be a two-edge refinement of $S_1'$ and $S_3$. Then $U$ is obtained from $S_3$ by blowing-up an edge at vertices whose stabilizers are conjugate to $D$. Moreover, $U$ is obtained from $S_1'$ by blowing-up an edge at vertices whose stabilizers are conjugate to $B$ or by blowing-up an edge at the vertices whose stabilizers are conjugate to $\left\langle x_1 \right\rangle \ast \left\langle g \right\rangle$. But, the second case can only occur when $\mathcal{S}_3=\mathcal{S}_1$ (see the claim in the proof of Lemma~\ref{Lem compatibility twists Wn-1}). Therefore, we may suppose that $U$ is obtained from $S_1'$ by blowing up an edge at vertices whose stabilizers are conjugate to $B$. Thus, up to applying a global conjugation, we may assume that $\left\langle x_1 \right\rangle \ast \left\langle g \right\rangle \subseteq D$. But, as $g$ is not contained in any proper $\mathcal{Z}_{RC}$-factor of $A$ and as $A \cap D$ is a free factor of $A$, we see that $A \cap D=A$. Hence $A \ast \left \langle x_1 \right\rangle \subseteq D$, and, as $A \ast \left \langle x_1 \right\rangle$ is isomorphic to $W_{n-1}$, we have in fact $A \ast \left \langle x_1 \right\rangle = D$. It follows that $C=D$ and, by Lemma~\ref{Lem stab corank one Wn-1 star}, we see that $\mathcal{S}_2=\mathcal{S}_3$. Thus the only equivalence classes of $W_{n-1}$-stars fixed by finite index subgroups of $P'$ are $\mathcal{S}_1$ and $\mathcal{S}_2$. 

Therefore the only equivalence classes of $W_{n-2}$-stars fixed by finite index subgroups of $P'$ are the equivalence classes of the $W_{n-2}$-stars which refine $S_1$ and $S_2$. As $S_1$ and $S_2$ are refined by a unique (up to $W_n$-equivariant homeomorphism) $W_{n-2}$-star by Theorem~\ref{Theo scott swarup}, we conclude that $\mathcal{S}$ is the only equivalence class of $W_{n-2}$-star fixed by finite index subgroups of $P'$. This completes the proof.
\hfill\qedsymbol

\bigskip

\begin{prop}\label{Prop Stab Wn-2 satisfies PWn-2}
Let $n \geq 5$ and let $\Gamma$ be a finite index subgroup of $C_n$. Let $\mathcal{S}$ be the equivalence class of a $W_{n-2}$-star $S$. Then $\Stab_{\Gamma}(\mathcal{S})$ satisfies $(P_{W_{n-2}})$. Moreover, we can choose for the subgroup $K_1 \times K_2$ of Property~$(P_{W_{n-2}})$~$(1)$ the direct product of the groups of twists of $\mathcal{S}$ about the two edges of $S$. 
\end{prop}

\dem The fact that $\Stab_{\Gamma}(\mathcal{S})$ satisfies $(P_{W_{n-2}})~(2)$ follows from the fact that $\Stab_{\Gamma}(\mathcal{S})$ contains the stabilizer in $\Gamma$ of the equivalence class of a $W_3$-star obtained from $S$ by blowing-up $n-5$ edges at the center of $W_n \backslash S$. Indeed, Proposition~\ref{Prop Levitt stab}~$(3)$ ensures that the group of twists of a $W_3$-star is isomorphic to a direct product of $n-3$ copies of $W_3$. 

The fact that $\Stab_{\Gamma}(\mathcal{S})$ satisfies $(P_{W_{n-2}})~(3)$ follows from the fact that $\Stab_{\Gamma}(\mathcal{S})$ contains the stabilizer in $\Gamma$ of the equivalence class of a $W_2$-star obtained from $S$ by blowing-up $n-4$ edges at the center of $W_n \backslash S$. Indeed the group of twists of a $W_{2}$-star is isomorphic to a direct product of $n-2$ copies of $W_2$ by Proposition~\ref{Prop Levitt stab}~$(3)$.

Let us now prove that $\Stab_{\Gamma}(\mathcal{S})$ satisfies $(P_{W_{n-2}})~(1)$. Let $T'$ be the group of twists of $\mathcal{S}$ and let $T=T' \cap \Gamma$. The group $T$ is normal in $\Stab_{\Gamma}(\mathcal{S})$ since $\Gamma \subseteq C_n$. By Proposition~\ref{Prop Levitt stab}~$(3)$, the group $T'$ is isomorphic to $T_1' \times T_2'$, where, for $i \in \{1,2\}$, $T_i'$ is the group of twists in $\Out(W_n)$ about one edge of $W_n \backslash S$. For $i \in \{1,2\}$, let $T_i=T_i' \cap \Gamma$. For every $i \in \{1,2\}$, the group $T_i$ is a normal subgroup of $\Stab_{\Gamma}(\mathcal{S})$ and the group $T_1 \times T_2$ is a normal subgroup of $\Stab_{\Gamma}(\mathcal{S})$. Let $T_1^{(2)}$ be a finite index subgroup of $T_1$ and let $P'$ be a finite index subgroup of a nontrivial normal subgroup of $T_1^{(2)}$. We prove that the centralizer of $P'$ in $\Gamma$ contains $T_2$ as a finite index subgroup. This will conclude the proof of the proposition by symmetry of $T_1$ and $T_2$. By Lemma~\ref{Lem normal subgroup of twist Wn-2 fixes 2 Wn-1}, the equivalence class $\mathcal{S}$ is the only equivalence class of $W_{n-2}$-star fixed by every finite index subgroup of $P'$. Hence $C_{\Gamma}(P')$ fixes $\mathcal{S}$. 

Let $H$ be a finite index subgroup of $C_{\Gamma}(P')$ which fixes $\mathcal{S}$. Let $$W_n= \left\langle x_1 \right\rangle \ast \left\langle x_3,\ldots,x_n\right\rangle \ast \left\langle x_2 \right\rangle$$ be a free factor decomposition associated with $W_n \backslash S$ and $A=\left\langle x_3,\ldots,x_n\right\rangle$. By Proposition~\ref{Prop Levitt stab}~$(1)$, the kernel of the natural homomorphism $H \to \Out(A)$ is isomorphic to $H \cap T$.
We claim that the image of $H$ in $\Out(A)$ is finite. Indeed, as $P'$ is a finite index subgroup of a nontrivial normal subgroup of a finite index subgroup of $T_1$ and as $T_1$ is isomorphic to a finite index subgroup of $W_{n-2}$, we see that $P'$ is isomorphic to a finite index subgroup $N$ of a nontrivial normal subgroup of a finite index subgroup of $W_{n-2}$. By Proposition~\ref{Prop Zrc factors}~$(3)$, $N$ is not contained in any proper $\mathcal{Z}_{RC}$-factor of $W_{n-2}$. By Proposition~\ref{Prop Zrc factors}~$(4)$, there exists $g \in N$ such that $g$ is not contained in any proper $\mathcal{Z}_{RC}$-factor of $W_{n-2}$. Thus, there exists $g \in A$ such that $g$ is not contained in any proper $\mathcal{Z}_{RC}$-factor of $A$ and the twist about $g$ is contained in $P'$. As $H$ commutes with the twist about $g$, Lemma~\ref{Lem twists commuting element twistor} implies that $H$ preserves the conjugacy class of $g$. Hence, by Theorem~\ref{Theo Zrc factors conjugacy class fixed}, the image of $H$ in $\Out(A)$ is finite. 

Thus, $H \cap T$ has finite index in $H$ and in $C_{\Gamma}(P')$. But, as $H$ commutes with $P' \subseteq T_1$, and as $T_1$ is virtually a nonabelian free group, the intersection $H \cap T_2$ has finite index in $H \cap T$, hence has finite index in $C_{\Gamma}(P')$. This completes the proof.
\hfill\qedsymbol

\bigskip

\subsection{Groups satisfying $(P_{W_{n-2}})$ and stabilizers of $W_{n-1}$-stars}

We prove in this section that if $H$ is a subgroup of $\Out(W_n)$ which satisfies $(P_{W_{n-2}})$, then $H$ virtually fixes the equivalence class of a $W_{n-1}$-star. We first recall a general lemma.

\begin{lem}\label{Lem finite index subgroup of normal subgroup is virtually normal}
Let $G$ be a group and let $N$ be a finitely generated normal subgroup of $G$. Let $n \in \NN^*$. 

\medskip

\noindent{$(1)$ } There exist only finitely many subgroups of $N$ of index equal to $n$. 

\medskip

\noindent{$(2)$ } For every finite index subgroup $N'$ of $N$ there exists a finite index subgroup $G'$ of $G$ such that $N'$ is a normal subgroup of $G'$.
\end{lem}

\dem Assertion~$(1)$ is well known, we only prove assertion~$(2)$. Let $N'$ be a subgroup of $N$ of index $n$ and let $g \in G$. As $N$ is a normal subgroup of $G$, the automorphism $\ad_g \colon G \to G$ induces an automorphism $\ad_g|_N \colon N \to N$ by restriction. Therefore, $\ad_g$ permutes the subgroups of index $n$ in $N$. Since there exists a finite number of subgroups of index $n$ in $N$ by the first assertion, we see that there exists a finite index subgroup $G'$ of $G$ such that, for every $g \in G'$, we have $\ad_g(N')=N'$. Therefore $N'$ is a normal subgroup of $G'$. This concludes the proof.
\hfill\qedsymbol

\bigskip

\begin{lem}\label{Lem Pwn-2 alternative fixes free splitting}
Let $n \geq 5$. Let $H$ be a subgroup of $C_n$ satisfying $(P_{W_{n-2}})$. Let $K_1 \times K_2$ be a normal subgroup of $H$ given by $(P_{W_{n-2}})~(1)$. Then one of the following holds.

\medskip

\noindent{$(1)$ } For every $i \in \{1,2\}$, the group $K_i$ does not virtually fix the equivalence class of a free splitting.

\medskip

\noindent{$(2)$ } The group $H$ virtually fixes the equivalence class of a one-edge free splitting.
\end{lem}

\dem Suppose that there exists $i \in \{1,2\}$ such that $K_i$ virtually fixes the equivalence class of a free splitting. Up to reordering, we may assume that $i=1$. Let $K_1'$ be a finite index subgroup of $K_1$ which fixes the equivalence class of a free splitting, and let $\mathcal{C}$ be the set of all equivalence classes of free splittings fixed by $K_1'$. Since $K_1$ is a finitely generated normal subgroup of $H$, by Lemma~\ref{Lem finite index subgroup of normal subgroup is virtually normal}~$(2)$, there exists a finite index subgroup $H_0$ of $H$ such that $K_1'$ is a normal subgroup of $H_0$. In particular, the set $\mathcal{C}$ is preserved by $H_0$.

Suppose first that the set $\mathcal{C}$ is finite. Then the set $\mathcal{C}$ is virtually fixed pointwise by $H_0$. Hence the group $H$ virtually fixes the equivalence class of a free splitting.

So we may assume that the set $\mathcal{C}$ is infinite. Let $U_{\mathcal{C}}$ be the splitting provided by Theorem~\ref{Theo invariant splitting constructed}, and let $\mathcal{U}_{\mathcal{C}}$ be its equivalence class. By the equivariance property in Theorem~\ref{Theo invariant splitting constructed} the equivalence class $\mathcal{U}_{\mathcal{C}}$ is $H_0$-invariant. Suppose first that the splitting $U_{\mathcal{C}}$ contains an edge $e \in EU_{\mathcal{C}}$ with trivial stabilizer. Let $U'$ be the splitting obtained from $U_{\mathcal{C}}$ by collapsing every edge of $U_{\mathcal{C}}$ that are not contained in the orbit of $e$, and let $\mathcal{U}'$ be its equivalence class. Then $\mathcal{U}'$ is the equivalence class of a one-edge free splitting virtually fixed by $H$.

Thus, we may assume that all edge stabilizers of $U_{\mathcal{C}}$ are nontrivial. We show that this leads to a contradiction.  Let $H'$ be the subgroup of finite index in $H_0$ which acts trivially on $W_n \backslash U_{\mathcal{C}}$. We claim that the intersection of $H'$ with the group of twists of $U_{\mathcal{C}}$ is finite. Indeed, let $e$ be a half-edge of $U_{\mathcal{C}}$. As $W_n$ is virtually free, if the edge stabilizer $G_e$ of $e$ is not cyclic, the group of twists about this half-edge is trivial. Thus, as we suppose that all edge stabilizers are nontrivial, half-edges with nontrivial group of twists have cyclic stabilizers. But by Lemma~\ref{Lem twists about cyclic edge central} twists about edges with cyclic stabilizers are central in a finite index subgroup of $\Stab^0(U_{\mathcal{C}})$. Note that Remark~\ref{Rmq centralizer PWn-2}~$(2)$ implies that the center of every finite index subgroup of $H'$ is finite. Therefore the intersection of $H'$ with the group of twists is finite. By Remark~\ref{Remark proposition Levitt}, the intersection of $H'$ with the group of bitwists is finite. Thus, up to passing to a finite index subgroup, we may suppose that the map $$H' \to \prod_{v \in V(W_n \backslash U_{\mathcal{C}})} \Out(G_v)$$ given by the action on the vertex groups is injective. 

Let $VU_{\mathcal{C}}=V_1 \amalg V_2$ be the partition of $VU_{\mathcal{C}}$ given by Theorem~\ref{Theo invariant splitting constructed} and, for every $i \in \{1,2\}$, let $H_i$ be the subgroup of $H'$ made of all automorphisms whose image in $\prod_{v \in W_n \backslash V_i} \Out(G_v)$ is trivial. Then $H_1$ and $H_2$ centralize each other and, by Theorem~\ref{Theo invariant splitting constructed}~$(1)(b)$, the group $H_1 \cap K_1'$ is a finite index subgroup of $K_1'$. Thus $H_2$ centralizes a finite index subgroup of $K_1'$. We prove that $\rk_{\pr}(H_2) \geq 2$, which will contradict the fact that the centralizer of every finite index subgroup of $K_1'$ is virtually free. 

By Theorem~\ref{Theo invariant splitting constructed}~$(2)$, the set of all conjugacy classes of groups $G_v$ with \mbox{$v \in V_2$} is a free factor system of $W_n$. In particular, for every $v \in V_2$, there exists $k_v \in \{0,\ldots,n-1\}$ such that $G_v$ is isomorphic to $W_{k_v}$. Suppose first that $|W_n \backslash V_2| \geq 3$. In this case, by Theorem~\ref{Theo product rank Wn}~$(2)$ and since $\rk_{\pr}(\Out(W_3))=1$ and $\rk_{\pr}(\Out(W_2))=0$, for all $v \in V_2$, we have $\rk_{\pr}(\Out(W_{k_v})) \leq k_v-2$. Hence $$\rk_{\pr}\left(\prod_{v \in W_n \backslash V_2} \Out(G_v) \right) \leq n-6.$$ Since $\rk_{\pr}(H')=n-3$, using Lemma~\ref{Lem extension product rank}, we see that $\rk_{\pr}(H_2) \geq 3$. This leads to a contradiction. Suppose now that $|W_n \backslash V_2|=2$ and let $v_1,v_2 \in W_n \backslash V_2$ be distinct. Then for every $i \in \{1,2\}$ there exists $k_i \in \{1,\ldots,n-1\}$ such that $G_{v_i}$ is isomorphic to $W_{k_i}$. If $W_n=W_{k_1} \ast W_{k_2}$, then the group $H'$ virtually fixes the equivalence class of the one-edge free splitting determined by this free factor decomposition of $W_n$. So we may assume that $W_n \neq W_{k_1} \ast W_{k_2}$. This implies that $k_1 +k_2 \leq n-1$. Hence $$\rk_{\pr}\left(\prod_{v \in W_n \backslash V_2} \Out(G_v) \right) \leq n-5.$$ Since $\rk_{\pr}(H')=n-3$, using Lemma~\ref{Lem extension product rank}, we see that $\rk_{\pr}(H_2) \geq 2$. This leads to a contradiction. Suppose now that $|W_n \backslash V_2|=1$, and let $v \in W_n \backslash V_2$. Then there exists $k \in \{1,\ldots,n-1\}$ such that $G_v$ is isomorphic to $W_k$. Suppose first that $k \leq n-2$. Then by Theorem~\ref{Theo product rank Wn}~$(2)$, and since $\rk_{\pr}(\Out(W_3))=1$, $\rk_{\pr}(\Out(W_1))=0$ and $\rk_{\pr}(\Out(W_2))=0$, if $n \neq 5$, we have $$\rk_{\pr}\left(\Out(W_k) \right) \leq n-5.$$ Thus, by Lemma~\ref{Lem extension product rank}, we see that $\rk_{\pr}(H_2) \geq 2$. When $n=5$, the case where $k=3$ and $\rk_{\pr}\left(\Out(W_k) \right)=1=n-4$ can occur. But by Property~$(P_{W_{n-2}})~(3)$, the group $H'$ contains a subgroup isomorphic to $\ZZ^3$. Since $\Out(W_3)$ is virtually free, the group $H_2$ contains a subgroup isomorphic to $\ZZ^2$. This contradicts the fact that the centralizer of every finite index subgroup of $K_1'$ is virtually nonabelian free. Hence we have $k=n-1$. But then, by Lemma~\ref{Lem stab corank one Wn-1 star}, the group $H'$ (and hence the group $H$) virtually fixes the equivalence class of a $W_{n-1}$-star. This concludes the proof.
\hfill\qedsymbol

\bigskip

\begin{lem}\label{Lem Pwn-2 fixes homothety class of arational tree}
Let $n \geq 5$. Let $\mathcal{F}$ be a nonsporadic free factor system. Let $H$ be a subgroup of $C_n \cap \Out(W_n,\mathcal{F})$ containing a direct product of $n-3$ nonabelian free groups. Then $H$ cannot contain a finite index subgroup which fixes the homothety class of a $(W_n,\mathcal{F})$-arational tree.
\end{lem}

\dem Suppose towards a contradiction that $H$ has a finite index subgroup which fixes the equivalence class of a $(W_n,\mathcal{F})$-arational tree. Up to passing to a finite index subgroup, we may suppose that $H$ itself fixes the homothety class of a $(W_n,\mathcal{F})$-arational tree. By Lemma~\ref{Lem isometric stabilizer cyclic}, there exists a homomorphism from $H$ to $\ZZ$ whose kernel $K'$ is exactly the isometric stabilizer of a $(W_n,\mathcal{F})$-arational tree. Note that $K'$ contains a direct product of $n-3$ nonabelian free groups as it is the kernel of a homomorphism from $H$ to $\ZZ$. By Proposition~\ref{Prop isometric stabilizer free splitting fixed}, there exists a finite index subgroup $K$ of $K'$ such that $K$ fixes infinitely many equivalence classes of free splittings. Let $\mathcal{C}$ be the collection of all equivalence classes of free splittings fixed by $K$. 

We claim that $\mathcal{C}$ is in fact finite, which will lead to a contradiction. Since $K \subseteq C_n$, Lemma~\ref{Lem Cn acts as the identity on graph} implies that if $\mathcal{S}$ is the equivalence class of a free splitting $S$ fixed by $K$, then the group $K$ fixes the equivalence class of every one-edge free splitting onto which $S$ collapses. By Theorem~\ref{Theo scott swarup}, if $\mathcal{S}$ is the equivalence class of a free splitting $S$, then $\mathcal{S}$ is determined by the finite set of equivalence classes of one-edge free splittings onto which $S$ collapses. Therefore, it suffices to show that $K$ can only fix finitely many equivalence classes of one-edge free splittings. Let $\mathcal{S}$ be the equivalence class of a one-edge free splitting fixed by $K$. Since $K$ contains a direct product of $n-3$ nonabelian free groups, Theorem~\ref{Theo product rank Wn}~$(3)$ implies that $S$ is a $W_{n-1}$-star. Let $$W_n=\left\langle x_1,\ldots,x_{n-1} \right\rangle \ast \left\langle x_n \right\rangle$$ be a free factor decomposition associated with $S$ and let $A=\left\langle x_1,\ldots,x_{n-1} \right\rangle$. By Proposition~\ref{Prop Levitt stab}~$(1)$, the kernel of the natural homomorphism $K \to \Out(A)$ is the intersection of $K$ with the group of twists $T$ of $\mathcal{S}$. By Theorem~\ref{Theo product rank Wn}~$(2)$, the product rank of $\Out(A)$ is equal to $n-4$. Since $K$ contains a direct product of $n-3$ nonabelian free groups, we see that $K \cap T$ is infinite. Therefore, for every equivalence class $\mathcal{S}$ of a $W_{n-1}$-star $S$ fixed by $K$, the group $K$ contains an infinite twist about $\mathcal{S}$. 

Let $\mathcal{S}$ and $\mathcal{S}'$ be two distinct equivalence classes of $W_{n-1}$-stars fixed by $K$. Let $S$ be a representative of $\mathcal{S}$ and let $S'$ be a representative of $\mathcal{S}'$. We claim that $S$ and $S'$ are compatible. Indeed, by the above, there exists $f \in K$ of infinite order such that $f$ is a twist about $\mathcal{S}$. Since $f$ fixes $\mathcal{S}'$, Lemma~\ref{Lem compatibility twists Wn-1} implies that $S$ and $S'$ are compatible. Therefore, for every distinct equivalence classes $\mathcal{S}$ and $\mathcal{S}'$ of one-edge free splittings fixed by $K$, there exist $S \in \mathcal{S}$ and $S' \in \mathcal{S}'$ such that $S$ and $S'$ are compatible. By Theorem~\ref{Theo scott swarup}, this is only possible when $\mathcal{C}$ is finite. This leads to a contradiction since $K$ must fix infinitely many equivalence classes of free splittings. This concludes the proof.
\hfill\qedsymbol

\bigskip

\begin{prop}\label{Prop Pwn-2 fixes Wn-1 star}
Let $n \geq 5$. Let $H$ be a subgroup of $C_n$ satisfying $(P_{W_{n-2}})$. Then $H$ virtually fixes the equivalence class of a $W_{n-1}$-star.
\end{prop}

\dem The proof is inspired by \cite[Proposition~8.2]{HorbezWade20} and \cite[Proposition~6.5]{HenselHorbezWade19}. We prove that $H$ virtually fixes the equivalence class of a one-edge free splitting. Since $H$ contains a direct product of $n-3$ nonabelian free groups, we will then conclude by Theorem~\ref{Theo product rank Wn}~$(3)$. Suppose towards a contradiction that $H$ does not virtually fix the equivalence class of a one-edge free splitting. Let $\mathcal{F}$ be a maximal $H$-periodic free factor system. We can assume that $\mathcal{F}$ is nonsporadic otherwise $H$ virtually fixes the equivalence class of a one-edge free splitting and we are done. As $\mathcal{F}$ is maximal, by Proposition~\ref{Prop maximal free factor system bounded orbit}, the group $H$ acts with unbounded orbits on $\FF(W_n,\mathcal{F})$. 

Let $K_1 \times K_2$ be a normal subgroup of $H$ given by $(P_{W_{n-2}})~(1)$. Suppose first that $K_1 \times K_2$ does not contain a loxodromic element on $\FF(W_n, \mathcal{F})$. As $H$ has unbounded orbits on $\FF(W_n,\mathcal{F})$, Lemma~\ref{Lem hyperbolic space and commuting} implies that $K_1 \times K_2$ has a finite orbit in $\partial_{\infty}\FF(W_n,\mathcal{F})$. 

By Lemma~\ref{Lem fixed point at infinity FF}, there exists a finite index subgroup $K_1' \times K_2'$ of $K_1 \times K_2$ such that $K_1' \times K_2'$ fixes the homothety class of an arational $(W_n,\mathcal{F})$-tree $T$. Since $K_1 \times K_2$ does not contain a loxodromic element, $K_1' \times K_2'$ fixes $T$ up to isometry, not just homothety (see e.g. \cite[Proposition~6.2]{Guirardelhorbez19}). By Proposition~\ref{Prop isometric stabilizer free splitting fixed}, the group $K_1' \times K_2'$ virtually fixes infinitely many equivalence classes of $(W_n,\mathcal{F})$-free splittings. By Lemma~\ref{Lem Pwn-2 alternative fixes free splitting}, the group $H$ virtually fixes the equivalence class of a one-edge free splitting of $W_n$.

So we may suppose that there exists a loxodromic element $\Phi \in K_1 \times K_2$. First suppose that there exists a unique $i \in \{1,2\}$ such that the group $K_i$ contains a loxodromic element $\Phi_i$. We may assume, up to reordering, that only $K_2$ contains a loxodromic element $\Phi$. Therefore by Lemma~\ref{Lem hyperbolic space and commuting}, the group $K_1$ virtually fixes a point in $\partial_{\infty}\FF(W_n,\mathcal{F})$. By Lemma~\ref{Lem fixed point at infinity FF}, the group $K_1$ virtually fixes the homothety class an arational $(W_n,\mathcal{F})$-tree $T$. Let $K_1'$ be a normal subgroup of $K_1$ of finite index that is contained in $\Stab([T])$. As $K_1'$ does not contain any loxodromic element, as in the above step, $K_1'$ fixes $T$ up to isometry. By Proposition~\ref{Prop isometric stabilizer free splitting fixed}, the group $K_1'$ fixes the equivalence class of a free splitting relative to $\mathcal{F}$. By Lemma~\ref{Lem Pwn-2 alternative fixes free splitting}, the group $H$ virtually fixes the equivalence class of a one-edge free splitting of $W_n$.

Now suppose that for every $i \in \{1,2\}$, the group $K_i$ contains a loxodromic element. By Lemma~\ref{Lem hyperbolic space and commuting}, the whole group $H$ virtually fixes a point in $\partial_{\infty}\FF(W_n,\mathcal{F})$. By Lemma~\ref{Lem fixed point at infinity FF}, the group $H$ virtually fixes the homothety class of an arational tree. This contradicts Lemma~\ref{Lem Pwn-2 fixes homothety class of arational tree}.

Therefore, in all cases, the group $H$ virtually fixes the equivalence class $\mathcal{S}$ of a one-edge free splitting $S$. By Theorem~\ref{Theo product rank Wn}~$(3)$, since $H$ contains a direct product of $n-3$ nonabelian free groups, the group $H$ virtually fixes the equivalence class of a $W_{n-1}$-star.

\hfill\qedsymbol

\bigskip

We now prove a proposition which gives a sufficient condition for equivalence classes of $W_{n-1}$-stars provided by Proposition~\ref{Prop Pwn-2 fixes Wn-1 star} to be compatible.

\begin{prop}\label{Prop PWn-2 compatibility}
Let $n \geq 5$ and let $\Gamma$ be a subgroup of $C_n$ of finite index. Let $k \in \NN^*$ and let $H_1,\ldots,H_k$ be subgroups of $\Gamma$ which satisfy $(P_{W_{n-2}})$ and such that the intersection $\bigcap_{i=1}^k H_i$ contains a subgroup $H$ isomorphic to $\ZZ^{n-2}$. For $i \in \{1,\ldots,k\}$, let $\mathcal{S}_i$ be the equivalence class of a $W_{n-1}$-star $S_i$ which is virtually fixed by $H_i$. Then, for every $i,j \in \{1,\ldots,k\}$, the $W_{n-1}$-stars $S_i$ and $S_j$ are compatible. 
\end{prop}

\dem Let $i,j \in \{1,\ldots,k\}$ be distinct integers. Let $H'$ be a finite index subgroup of $H$ contained in $\Stab_{\Gamma}(\mathcal{S}_i) \cap \Stab_{\Gamma}(\mathcal{S}_j)$. Let $A_i$ and $A_j$ be the vertex groups isomorphic to $W_{n-1}$ of respectively $W_n \backslash S_i$ and $W_n \backslash S_j$ (well defined up to conjugation). By Proposition~\ref{Prop vcd}, the rank of a maximal abelian subgroup of $\Out(W_{n-1})$ is equal to $n-3$. Therefore, the kernel of the homomorphisms $H' \to \Out(A_i)$ and $H' \to \Out(A_j)$ given by the action on the vertex group contains an element of infinite order. Let $f_i \in \ker\left(H' \to \Out(A_i)\right)$ and $f_j \in \ker\left(H' \to \Out(A_i)\right)$ be infinite order elements. By Proposition~\ref{Prop Levitt stab}~$(1)$, $f_i$ and $f_j$ are twists about respectively $S_i$ and $S_j$. As $f_i$ and $f_j$ commute, by Corollary~\ref{Coro compatibility commuting twists}, $S_i$ and $S_j$ are compatible. This concludes the proof. 
\hfill\qedsymbol

\section{Algebraic characterization of stabilizers of $W_{n-2}$-stars}\label{Section characterisation stabiliser Wn-2 stars}

In this section, we give an algebraic characterization of stabilizers of $W_{n-2}$-stars. By the previous section, we know that groups which satisfy $(P_{W_{n-2}})$ virtually stabilize equivalence classes of $W_{n-1}$-stars, and we have given an algebraic criterion to show that these $W_{n-1}$-stars are compatible. In order to prove that a group $H$ which satisfies $(P_{W_{n-2}})$ virtually stabilizes the equivalence class of a $W_{n-2}$-star, we study the intersection of a normal subgroup $K_1 \times K_2$ of $H$ given by $(P_{W_{n-2}})~(1)$ with the group of twists of the equivalence class of a $W_{n-1}$-star virtually fixed by $H$.

\subsection{Groups of twists in groups satisfying $(P_{W_{n-2}})$}

We start this section with a lemma which gives a sufficient condition for a group $H$ satisfying $(P_{W_{n-2}})$ to be the stabilizer of a $W_{n-2}$-star.

\begin{lem}\label{Lem PWn-2 twists stab Wn-2}
Let $n \geq 5$ and let $\Gamma$ be a subgroup of finite index of $C_n$. Let $H$ be a subgroup of $\Gamma$ which satisfies $(P_{W_{n-2}})$ and let $K_1 \times K_2$ be a normal subgroup of $H$ given by $(P_{W_{n-2}})~(1)$. Let $\mathcal{S}_1$ be the equivalence class of a $W_{n-1}$-star $S_1$ virtually fixed by $H$ and let $T_1$ be the group of twists of $\mathcal{S}_1$. 

Suppose that $T_1 \cap K_1$ is infinite and that there exists an equivalence class $\mathcal{S}_2$ of a $W_{n-1}$-star $S_2$ such that the intersection of $K_2$ with the group of twists $T_2$ of $\mathcal{S}_2$ is infinite. Then $S_1$ and $S_2$ are compatible and $H$ virtually fixes the equivalence class $\mathcal{S}$ of the $W_{n-2}$-star which refines $S_1$ and $S_2$. Moreover, $\mathcal{S}$ is the unique equivalence class of a $W_{n-2}$-star virtually fixed by $H$. Finally, the groups $T_1 \cap \Stab_{\Gamma}(\mathcal{S})$ and $K_1$ (resp. $T_2 \cap \Stab_{\Gamma}(\mathcal{S})$ and $K_2$) are commensurable.
\end{lem}

\dem For $i \in \{1,2\}$, let $f_i \in T_i \cap K_i$ be of infinite order. First remark that, as $f_1$ and $f_2$ generate a free abelian group of order $2$, we have $T_1 \neq T_2$ because the group of twists of a $W_{n-1}$-star is virtually a nonabelian free group. Hence we have $\mathcal{S}_1 \neq \mathcal{S}_2$. As $K_1$ commutes with $f_2$, Proposition~\ref{Prop commutator twist splitting fixed} shows that $K_1$ fixes $\mathcal{S}_2$. As $K_1$ contains a twist of $\mathcal{S}_1$, Lemma~\ref{Lem compatibility twists Wn-1} shows that $S_1$ and $S_2$ are compatible.

Let $S$ be a $W_{n-2}$-star which refines $S_1$ and $S_2$, let $\mathcal{S}$ be its equivalence class and let $T$ be the group of twists of $\mathcal{S}$ in $\Gamma$. Then $T$ contains a finite index normal subgroup isomorphic to $K_1^{\mathcal{S}_1} \times K_2^{\mathcal{S}_2}$, where $K_1^{\mathcal{S}_1}$ and $K_2^{\mathcal{S}_2}$ are virtually nonabelian free groups. By Proposition~\ref{Prop Stab Wn-2 satisfies PWn-2}, we can choose $K_1^{\mathcal{S}_1} \times K_2^{\mathcal{S}_2}$ such that $K_1^{\mathcal{S}_1} \times K_2^{\mathcal{S}_2}$ is a group satisfying Property~$(P_{W_{n-2}})~(1)$. Moreover, up to reordering, $K_1^{\mathcal{S}_1} \subseteq T_1$ and $K_2^{\mathcal{S}_2} \subseteq T_2$. Since $K_1$ fixes both $\mathcal{S}_1$ and $\mathcal{S}_2$, we see that $K_1$ fixes $\mathcal{S}$. Therefore, by Proposition~\ref{Prop Levitt stab}~$(1)$, we have a homomorphism $\Phi \colon K_1 \to \Out(W_{n-2})$ whose kernel is exactly $K_1 \cap T$. By Lemma~\ref{Lem twists Wn-1 star in stab Wn-2}, we see that $T_1 \cap \Stab_{\Gamma}(\mathcal{S}) \cap K_1^{\mathcal{S}_1}$ is a finite index subgroup of $T_1 \cap \Stab_{\Gamma}(\mathcal{S})$. As $K_1 \cap T_1$ is infinite, so is $K_1 \cap K_1^{\mathcal{S}_1}$. Let $P=\ker(\Phi) \cap K_1^{\mathcal{S}_1}=K_1 \cap K_1^{\mathcal{S}_1}$. Then, since $K_1 \subseteq C_n$, the group $K_1^{\mathcal{S}_1} \cap K_1$ is a normal subgroup of $K_1$. Therefore $P$ is a nontrivial normal subgroup of $K_1$. By Property~$(P_{W_{n-2}})~(1)$, we see that $K_2$ is a finite index subgroup of $C_{\Gamma}(P)$. But $P$ is centralized by $K_2^{\mathcal{S}_2}$ since $P \subseteq K_1^{\mathcal{S}_1}$. Hence $K_2^{\mathcal{S}_2} \cap K_2$ is a finite index subgroup of $K_2^{\mathcal{S}_2}$. As $K_1^{\mathcal{S}_1}$ is a finite index subgroup of the centralizer of $K_2^{\mathcal{S}_2}$ by Property~$(P_{W_{n-2}})~(1)$, and as $K_1$ is a finite index subgroup of the centralizer of $K_2$, we see that $K_1^{\mathcal{S}_1} \cap K_1$ has finite index in $K_1$ and therefore $P$ has finite index in $K_1$. Let $$W_n = \left\langle x_1 \right\rangle \ast \left\langle x_3,\ldots,x_n\right\rangle \ast \left\langle x_2 \right\rangle$$ be the free factor decomposition of $W_n$ induced by $S$ and let $A=\left\langle x_3,\ldots,x_n\right\rangle$. Then, up to reordering, for every $f \in P$, there exists $z_f \in A$ and a representative $F$ of $f$ such that $F$ sends $x_1$ to $z_fx_1z_f^{-1}$, and, for every $i \neq 1$, fixes $x_i$.

\medskip

\noindent{\bf Claim. } The only equivalence classes of $W_{n-1}$-stars which are virtually fixed by $K_1$ are $\mathcal{S}_1$ and $\mathcal{S}_2$.

\medskip

\dem Let $\mathcal{S}_3$ be the equivalence class of a $W_{n-1}$-star $S_3$ virtually fixed by $K_1$. Suppose towards a contradiction that $\mathcal{S}_3$ is distinct from both $\mathcal{S}_1$ and $\mathcal{S}_2$. Let $K_1'=K_1 \cap \Stab_{\Gamma}(\mathcal{S}_3)$ and $P'=P \cap \Stab_{\Gamma}(\mathcal{S}_3)$. Then, as $P$ is an infinite subgroup of the group of twists of $\mathcal{S}_1$, and as $P'$ is a finite index subgroup of $P$, we see that $P'$ is an infinite subgroup of the group of twists of $\mathcal{S}_1$. By Lemma~\ref{Lem compatibility twists Wn-1}, we see that $S_1$ and $S_3$ are compatible. Let $S'$ be a $W_{n-2}$-star that refines $S_1$ and $S_3$ and let $\mathcal{S}'$ be its equivalence class. Let $$W_n = \left\langle y_1 \right\rangle \ast \left\langle y_3,\ldots,y_n\right\rangle \ast \left\langle y_2 \right\rangle$$ be the free factor decomposition of $W_n$ induced by $S'$ and let $B=\left\langle y_3,\ldots,y_n\right\rangle$. Since $S$ is a refinement of $S_1$, we may suppose that $B \ast \left\langle y_2 \right\rangle=  A \ast \left\langle x_2 \right\rangle$ and that $y_1$ is a conjugate of $x_1$ by an element of $B \ast \left\langle y_2 \right\rangle$. Up to applying a global conjugation, we may also suppose that $y_1=x_1$ and that $B \ast \left\langle y_2 \right\rangle=  A \ast \left\langle x_2 \right\rangle$.

Let $T'$ be the group of twists of $\mathcal{S}'$. Then $T'$ contains a finite index normal subgroup isomorphic  to $P_1' \times P_2'$, where both $P_1'$ and $P_2'$ are virtually nonabelian free subgroups of $T'$ which correspond to the groups of twists about the two edges of $W_n \backslash S'$. Then, as $P'$ is a group of twists of $\mathcal{S}_1$, and as $P'$ fixes $\mathcal{S}'$, by Lemma~\ref{Lem twists Wn-1 star in stab Wn-2}, up to reordering, the group $P'$ is contained in $P_1'$. 

Let $f' \in P_1'$, let $F'$ be the representative of $f'$ which acts as the identity on $ B \ast \left\langle y_2 \right\rangle$ and let $z_{f'} \in B$ be the twistor of $F'$. Then $F'$ acts as the identity on $A \ast \left\langle x_2 \right\rangle$ and $F'(x_1)=z_{f'}x_1z_{f'}^{-1}$.
Recall that for every $\psi \in P'$, there exists a unique $z_{\psi} \in A$ and a unique representative $\Psi$ of $\psi$ such that $\Psi$ sends $x_1$ to $z_{\psi}x_1z_{\psi}^{-1}$, and, for every $i \neq 1$, fixes $x_i$. Thus, a necessary condition for $f'$ to be in $P'$ is that $z_{f'} \in A \cap B$.

But as $A$ and $B$ are free factors of $W_n$, the group $A \cap B$ is a free factor of $B$. To see this, let $U$ be a free splitting of $W_n$ such that $A$ is a vertex stabilizer of $U$ and let $U_B$ be the minimal subtree of $B$ in $U$. Then, as $U$ is a free splitting of $W_n$, we see that $U_B$ is a free splitting of $B$. But then, as $A$ is a vertex stabilizer in $U$, we see that $A \cap B$ is a vertex stabilizer in $U_B$. Therefore, $A \cap B$ is a free factor of $B$. Thus one can find a $W_{n-3}$-star $S^{(2)}$ which refines $S'$ and such that, for every $f' \in P'$, the twistor $z_{f'}$ fixes a vertex of $S^{(2)}$. Indeed, one can equivariantly blow-up an edge $e$ at the vertex of $S'$ whose stabilizer is $B$ such that the stabilizer of one of the endpoints of $e$ is a subgroup $C$ isomorphic to $W_{n-3}$ with $A \cap B \subseteq C$. Therefore we may also assume that $S^{(2)}$ is a $W_{n-3}$-star. Let $\mathcal{S}^{(2)}$ be the equivalence class of $S^{(2)}$. By Proposition~\ref{Prop Levitt stab}~$(3)$, the group of twists of $\mathcal{S}^{(2)}$ is isomorphic to a direct product $W_{n-3}^3$ of three infinite groups, where each factor is a group of twists about an edge of $W_n \backslash S^{(2)}$. This implies that $P'$ is contained in exactly one of the three factors isomorphic to $W_{n-3}$. It follows that the centralizer of $P'$ contains two elements which generates a free abelian group of order $2$. This contradicts the fact that the centralizer of $P'$ is virtually a nonabelian free group by $(P_{W_{n-2}})$~$(1)$. The claim follows.
\hfill\qedsymbol

\bigskip

The claim above then implies, as $K_1$ is a normal subgroup of $H$, that $H$ virtually fixes $\mathcal{S}_2$. As $H$ virtually fixes $\mathcal{S}_1$, we see that $H$ virtually fixes the equivalence class $\mathcal{S}$. Moreover, the above claim shows that $\mathcal{S}$ is the unique equivalence class of a $W_{n-2}$-star virtually fixed by $K_1$, and hence virtually fixed by $H$. 
 
We finally prove that $K_1$ and $T_1\cap \Stab_{\Gamma}(\mathcal{S})$ (resp. $K_2$ and $T_2\cap \Stab_{\Gamma}(\mathcal{S})$) are commensurable. By Lemma~\ref{Lem twists Wn-1 star in stab Wn-2}, for every $i \in \{1,2\}$ we see that $K_i^{\mathcal{S}_i} \cap T_i\cap \Stab_{\Gamma}(\mathcal{S})$ is a finite index subgroup of $T_i\cap \Stab_{\Gamma}(\mathcal{S})$. Moreover, for every $i \in \{1,2\}$ and every $f \in K_i^{\mathcal{S}_i}$, the twist $f$ of $\mathcal{S}$ is also a twist of $\mathcal{S}_i$. Hence we have $K_i^{\mathcal{S}_i} \subseteq T_i\cap \Stab_{\Gamma}(\mathcal{S})$. Therefore, for every $i \in \{1,2\}$, the groups $K_i^{\mathcal{S}_i}$ and $T_i\cap \Stab_{\Gamma}(\mathcal{S})$ are commensurable. Hence it suffices to show that, for every $i \in \{1,2\}$, the groups $K_i$ and $K_i^{\mathcal{S}_i}$ are commensurable. 

Recall that $K_2^{\mathcal{S}_2} \cap K_2$ is a finite index subgroup of $K_2^{\mathcal{S}_2}$ and that $K_1^{\mathcal{S}_1} \cap K_1$ has finite index in $K_1$. Since $H$ virtually fixes $\mathcal{S}$, and since $K_2^{\mathcal{S}_2}$ is a normal subgroup of $\Stab_{\Gamma}(\mathcal{S})$, we see that $K_2^{\mathcal{S}_2} \cap K_2$ is a normal subgroup of a finite index subgroup of $K_2$. We know that $K_2^{\mathcal{S}_2} \cap K_2$ commutes with $K_1^{\mathcal{S}_1}$ because $K_1^{\mathcal{S}_1}$ and $K_2^{\mathcal{S}_2}$ commute with each other. Thus, by Property~$(P_{W_{n-2}})~(1)$ applied to $K_1 \times K_2$, the centralizer of $K_2^{\mathcal{S}_2} \cap K_2$ contains $K_1$ as a finite index subgroup. This shows that $K_1 \cap K_1^{\mathcal{S}_1}$ is a finite index subgroup of $K_1^{\mathcal{S}_1}$. Hence $K_1$ and $K_1^{\mathcal{S}_1}$ are commensurable. By Property~$(P_{W_{n-2}})~(1)$ applied to $K_1^{\mathcal{S}_1} \times K_2^{\mathcal{S}_2}$, the centralizer of a finite index subgroup of $K_1^{\mathcal{S}_1}$ contains $K_2^{\mathcal{S}_2}$ as a finite index subgroup. Moreover, the centralizer of a finite index subgroup of $K_1$ contains $K_2$ as a finite index subgroup. Hence the centralizer of $K_1 \cap K_1^{\mathcal{S}_1}$ contains both $K_2$ and $K_2^{\mathcal{S}_2}$ as finite index subgroups. Thus $K_2$ and $K_2^{\mathcal{S}_2}$ are commensurable. This completes the proof of Lemma~\ref{Lem PWn-2 twists stab Wn-2}.
\hfill\qedsymbol

\bigskip

Lemma~\ref{Lem PWn-2 twists stab Wn-2} suggests that in order to show that a group $H$ which satisfies $(P_{W_{n-2}})$ is in fact virtually the stabilizer of the equivalence class of a $W_{n-2}$-star, it suffices to study the intersection of $H$ with groups of twists. A first step towards such a result is the following lemma.

\begin{lem}\label{Lem PWn-2 one normal subgroup intersects twists}
Let $n \geq 5$ and let $\Gamma$ be a subgroup of $C_n$ of finite index. Let $H$ be a subgroup of $\Gamma$ satisfying $(P_{W_{n-2}})$ and let $K_1 \times K_2$ be a normal subgroup of $H$ given by $(P_{W_{n-2}})~(1)$. Let $\mathcal{S}$ be the equivalence class of a $W_{n-1}$-star $S$ virtually fixed by $H$ and let $T$ be the group of twists of $\mathcal{S}$ contained in $\Gamma$.

There exists a unique $i \in \{1,2\}$ such that $K_i \cap T$ is infinite. Moreover, $H \cap T \cap K_i$ has finite index in $H \cap T$.
\end{lem}

\dem Up to passing to a finite index subgroup of $H$, we may suppose that $H$ fixes $\mathcal{S}$. The uniqueness assertion follows from the fact that $T$ is virtually a nonabelian free group and that $K_1 \times K_2$ is a direct product. Therefore, up to reordering, we may suppose that $K_1 \cap T$ is finite. So there exists a finite index subgroup $K_1'$ of $K_1$ such that $K_1' \cap T$ is trivial. Since $K_1$ is a finitely generated normal subgroup of $H$, Lemma~\ref{Lem finite index subgroup of normal subgroup is virtually normal} implies that there exists a finite index subgroup $H'$ of $H$ such that $K_1'$ is a normal subgroup of $H'$. Therefore, we may suppose that $K_1 \cap T$ is trivial. By Proposition~\ref{Prop Levitt stab}~$(1)$, the natural homomorphism $K_1 \to \Out(W_{n-1})$ given by the action on the vertex groups is injective.

We claim that $H \cap T$ is infinite. Indeed, consider the natural homomorphism $ \Phi \colon H \to \Out(W_{n-1})$. By Proposition~\ref{Prop vcd}, the rank of a maximal free abelian subgroup of $\Out(W_{n-1})$ is equal to $n-3$. As $H$ contains a subgroup isomorphic to $\ZZ^{n-2}$ by $(P_{W_{n-2}})~(3)$, the kernel of $H \to \Out(W_{n-1})$ is infinite. But, by Proposition~\ref{Prop Levitt stab}~$(1)$, this is precisely $H \cap T$. Therefore, $H \cap T$ is infinite.

We now prove that $H \cap T \cap K_2$ has finite index in $H \cap T$. This will conclude the proof as $H \cap T$ is infinite. Let $K=\Phi^{-1}(\Phi(K_2))$. Note that $H \cap T \subseteq K$. Then, as $K_2$ is normal in $H$, we see that $K$ is a normal subgroup of $H$ which contains $H \cap T$ and $K_2$. We claim that $K \cap K_1$ is finite. Indeed, suppose towards a contradiction that there exists $f \in K \cap K_1$ of infinite order. Then, as the homomorphism $$\Phi|_{K_1} \colon K_1 \to \Out(W_{n-1})$$ is injective, the element $\Phi(f)$ has infinite order. By definition of $K$, we see that \mbox{$\Phi(f) \in \Phi(K_1) \cap \Phi(K_2)$}. But, as the homomorphism $\Phi|_{K_1} \colon K_1 \to \Out(W_{n-1})$ is injective, and as $K_1$ is virtually a nonabelian free group, there exists $g \in K_1$ of infinite order such that $\Phi(g)$ does not commute with $\Phi(f)$. Since $\Phi(f) \in \Phi(K_2)$ this contradicts the fact that $K_1$ and $K_2$ commute with each other. Hence $K \cap K_1$ is finite. 

The groups $K$ and $K_1$ are two normal subgroups of $H$ with finite intersection. Let $K_1^{(2)}$ be a finite index normal subgroup of $K_1$ such that $K \cap K_1^{(2)}=\{1\}$. Since $K_1$ is finitely generated, by Lemma~\ref{Lem finite index subgroup of normal subgroup is virtually normal}~$(2)$, there exists a finite index subgroup $H^{(2)}$ of $H$ such that $K_1^{(2)}$ is a normal subgroup of $H^{(2)}$. Hence $K_1^{(2)}$ and $K \cap H^{(2)}$ are normal subgroups of $H^{(2)}$ with trivial intersection. Therefore, $K \cap H^{(2)} \subseteq C_{\Gamma}(K_1^{(2)})$. But, Property~$(P_{W_{n-2}})~(1)$ implies that $K$ and $K_2$ are commensurable. Since $K$ contains $H \cap T$, we see that $K_2 \cap T$ and $H \cap T$ are commensurable. This concludes the proof.
\hfill\qedsymbol

\subsection{Groups satisfying $(P_{W_{n-2}})$ and stabilizers of $W_{n-2}$-stars}

In this section we prove that a subgroup of $C_n$ which satisfies $(P_{W_{n-2}})$ virtually fixes the equivalence class of a $W_{n-2}$-star. We first recall a theorem due to Guirardel and Levitt which provides a canonical splitting for a relative one-ended group (recall that a group $G$ is \emph{one-ended relative to a family of subgroups} $\mathcal{H}$ if $G$ does not have a one-edge splitting with finite edge stabilizers such that every subgroup of $\mathcal{H}$ fixes a point).

\begin{theo}\cite[Theorem~9.14]{guirardel2016jsj}\label{Theo canonical splitting for one-ended group}
Let $G$ be a hyperbolic group and let $\mathcal{H}$ be a family of subgroups such that $G$ is one-ended relative to $\mathcal{H}$. There exists a JSJ splitting $T$ such that:

\medskip

\noindent{$(1)$ } Every edge stabilizer is virtually infinite cyclic.

\medskip

\noindent{$(2)$ } For every $H \in \mathcal{H}$, the group $H$ is elliptic in $T$.

\medskip

\noindent{$(3)$ } The tree $T$ is invariant under all automorphisms of $G$ preserving $\mathcal{H}$. Moreover, $T$ is compatible with every splitting $S$ with virtually cyclic edge stabilizers and such that for every $H \in \mathcal{H}$, the group $H$ is elliptic in $S$.
\end{theo}

We also need some results about splittings over virtually cyclic groups, whose generalization to virtually free groups is due to Cashen.

\begin{theo}\cite[Theorem~1.2]{Cashen2017} \label{Theo splittings over virtually cyclic groups}
Let $G_1$ and $G_2$ be finitely generated virtually nonabelian free groups, and let $C$ be a virtually cyclic group which is a proper subgroup of both $G_1$ and $G_2$. Then $G_1 \ast_C G_2$ is virtually a nonabelian free group if and only if there exists $i \in \{1,2\}$ such that $G_i$ has a splitting with finite edge stabilizers such that $C$ is a vertex stabilizer.
\end{theo}

\begin{coro}\label{Coro splittings over W2}
Let $n \geq 3$ and let $G_1$, $G_2$ be subgroups of $W_n$ such that $W_n=G_1 \ast_C G_2$ is a (possibly trivial) amalgamated product of $W_n$, where $C$ is isomorphic to $W_2$. 

\medskip

\noindent{$(1)$ } There exists $i \in \{1,2\}$ such that $C$ is a free factor of $G_i$. Moreover, if $j \in \{1,2\}-\{i\}$, then $G_j$ is a free factor of $W_n$. 

\medskip

\noindent{$(2)$ } Suppose that $W_n=G_1 \ast_C G_2$ is a nontrivial amalgamated product of $W_n$. Then there exist $3 \leq k_1,k_2 \leq n-1$ such that $k_1 +k_2=n+2$ and, for every $i \in \{1,2\}$, the group $G_i$ is isomorphic to $W_{k_i}$. 

\medskip

\noindent{$(3)$ } If $n=3$, there exists $i \in \{1,2\}$ such that $G_i=C$.
\end{coro}

\dem $(1)$ The assertion is immediate when $G_1 \ast_C G_2$ is trivial, so suppose that $C$ is a proper subgroup of both $G_1$ and $G_2$. Since $W_n$ is virtually a nonabelian free group, since every maximal virtually cyclic subgroups of $W_n$ is isomorphic to $W_2$ and since the amalgamated product $W_n=G_1 \ast_C G_2$ is nontrivial, the groups $G_1$ and $G_2$ are virtually nonabelian free groups. Moreover, since $W_n$ and $C$ are finitely generated, so are $G_1$ and $G_2$. By Theorem~\ref{Theo splittings over virtually cyclic groups}, up to exchanging the roles of $G_1$ and $G_2$, we may suppose that $G_1$ has a splitting $S$ such that every edge stabilizer is finite and $C$ is the stabilizer of a vertex $v \in VS$. Note that, since the finite subgroups of $W_n$ are all isomorphic to $F$, every edge stabilizer of $S$ is either trivial or isomorphic to $F$. In particular, every nontrivial edge stabilizer is a free factor of both of its endpoint stabilizers. Let $V_1$ be the set of vertices of $S$ distinct from $v$ and fixed by an element of $C$ isomorphic to $F$. Therefore, for every $w \in V_1$, there exists a subgroup $A_w$ of $G_w$ and an element $x_w \in C$ of order $2$ such that $G_w=A_w \ast \left\langle x_w \right\rangle$. Let $S'$ be the splitting of $G_1$ defined as follows. The underlying tree of $G_1 \backslash S'$ is the same one as the underlying tree of $G_1 \backslash S$, the stabilizer in $G_1$ of every vertex of $S'$ not in the orbit of an element of $V_1$ is the same one as the stabilizer of the corresponding vertex of $S$, and the stabilizer in $G_1$ of an element $w$ of $V_1$ is $A_w$. Then the stabilizer in $G_1$ of every edge of $S'$ adjacent to the vertex fixed by $C$ has trivial stabilizer. Thus, $C$ is a free factor of $G_1$ and there exists $H_1 \subset G_1$ such that $G_1=H_1 \ast C$. This proves the first assertion of $(1)$. The second assertion of $(1)$ follows from the fact that $$W_n=G_1 \ast_C \ast G_2= \left(H_1 \ast C\right) \ast_C G_2=H_1 \ast G_2.$$
Hence $H_1$ and $G_2$ are free factors of $W_n$. 

\medskip

\noindent{$(2)$ } Therefore, there exist $h_1,k_2 \in \{1,\ldots,n-2\}$ with $h_1+k_2=n$ such that $H_1$ is isomorphic to $W_{h_1}$ and $G_2$ is isomorphic to $W_{k_2}$. Thus $G_1$ is isomorphic to $W_{h_1+2}$. Set $k_1=h_1+2$. Since the amalgamated product is nontrivial, we have $3 \leq k_1,k_2 \leq n-1$. This proves $(2)$. 

\medskip

\noindent{$(3)$ } Finally, the third point is a direct consequence of the inequality $3 \leq k_1,k_2 \leq n-1$, which is impossible when $n=3$. This concludes the proof. 
\hfill\qedsymbol

\bigskip

In the next lemma, we will use the notion of the \emph{abelian rank} of a group $G$. The \emph{abelian rank} of $G$ is the rank of a maximal free abelian subgroup of $G$. It is closed under taking finite index subgroups.

\begin{lem}\label{Lem rank product splitting over W2}
Let $n \geq 4$ and let $S$ be a one-edge splitting of $W_n$ whose edge stabilizers are isomorphic to $W_2$. Let $\mathcal{S}$ be its equivalence class. Let $v_1$ and $v_2$ be adjacent vertices of $S$ and let $e$ be the edge between $v_1$ and $v_2$. Suppose that the abelian rank of the image of the natural homomorphism $$ \Phi \colon\Stab(\mathcal{S}) \to \Out(G_{v_1},G_e) \times \Out(G_{v_2},G_e)$$ is equal to the abelian rank of $\Out(W_n)$. For every $i \in \{1,2\}$, let $$ \Phi_i \colon\Stab(\mathcal{S}) \to \Out(G_{v_i},G_e)$$ be the natural homomorphism given by the action on the vertex group.

\medskip

\noindent{$(1)$ } For every $i \in \{1,2\}$, the abelian rank of $\Phi_i(\Stab(\mathcal{S}))$ is equal to the abelian rank of $\Out(G_{v_i})$. 

\medskip

\noindent{$(2)$ } For every $i \in \{1,2\}$, there exists a refinement of $S$ by blowing-up a one-edge $\mathcal{Z}_{RC}$-splitting of $G_{v_i}$ at $v_i$.
\end{lem}

\dem $(1)$ By Corollary~\ref{Coro splittings over W2}~$(2)$, for every $i \in \{1,2\}$, there exists $k_i \in \NN^*$ such that $G_{v_i}$ is isomorphic to $W_{k_i}$. Moreover, we have $k_1+k_2=n+2$ and, for every $i \in \{1,2\}$, we have $k_i \geq 3$. By Proposition~\ref{Prop vcd}, for every $i \in \{1,2\}$, the abelian rank of $\Out(G_{v_i})$ is equal to $k_i-2$ and the abelian rank of $\Out(W_n)$ is equal to $n-2$. By the assumption of the lemma, the abelian rank of $\Phi(\Stab(\mathcal{S}))$ is hence equal to $n-2$. This implies that $$n-2 \leq k_1-2+k_2-2=n+2-4=n-2.$$ Therefore, we conclude that for every $i \in \{1,2\}$, the abelian rank of $\Phi_i(\Stab(\mathcal{S}))$ is equal to the abelian rank of $\Out(G_{v_i})$. This proves Assertion~$(1)$.

\medskip

\noindent{$(2)$ } Let $i \in \{1,2\}$. By the first assertion the group $\Out(G_{v_i},G_e)$ is infinite. Since $G_e$ is isomorphic to $W_2$, the group $\Out(G_e)$ is finite. Therefore, the group $\Out(G_{v_i},G_e^{(t)})$ is infinite. By Theorem~\ref{Theo Zrc factors conjugacy class fixed}, the group $G_{v_i}$ has a $\mathcal{Z}_{RC}$ splitting $U$ such that $G_e$ fixes exactly one vertex. Let $S'$ be the splitting obtained from $S$ by blowing-up $U$ at $v_i$ and by attaching $e$ to the vertex of $U$ fixed by $G_e$. Then $S'$ satisties Assertion~$(2)$. This concludes the proof.
\hfill\qedsymbol

\begin{prop}\label{Prop QWn-2 is stab Wn-2}
Let $n \geq 5$ and let $\Gamma$ be a finite index subgroup of $C_n$. Let $H$ be a subgroup of $\Gamma$ which satisfies $(P_{W_{n-2}})$. Then $H$ virtually stabilizes the equivalence class of a $W_{n-2}$-star. Moreover, this equivalence class is unique.
\end{prop}

\dem By Proposition~\ref{Prop Pwn-2 fixes Wn-1 star}, the group $H$ virtually fixes the equivalence class $\mathcal{S}$ of a $W_{n-1}$-star $S$. Let $$W_n=A \ast \left\langle x_n \right\rangle$$ be the free factor decomposition of $W_n$ induced by $S$. Up to passing to a finite index subgroup, we may suppose that $H$ fixes $\mathcal{S}$. Let $T$ be the group of twists of $\mathcal{S}$ contained in $\Gamma$. By Proposition~\ref{Prop Levitt stab}~$(2)$, the group $\Stab(\mathcal{S})$ is isomorphic to $\Aut(A)$ and the group of twists of $\mathcal{S}$ is identified with the inner automorphism group of $A$.

Let $K_1 \times K_2$ be a normal subgroup of $H$ given by Property~$(P_{W_{n-2}})~(1)$. By Lemma~\ref{Lem PWn-2 one normal subgroup intersects twists}, up to exchanging the roles of $K_1$ and $K_2$, we may assume that $K_1 \cap T$ is infinite, that $H \cap T \cap K_1$ is a finite index subgroup of $H \cap T$ and that $K_2 \cap T$ is finite. Up to passing to a finite index subgroup, we may assume that $K_2 \cap T=\{1\}$. In particular, the natural homomorphism $\phi \colon K_2 \to \Out(A)$ is injective. Let $K \subseteq A$ be the group of twistors associated with twists contained in $K_1$. Note that to every splitting $S_0$ of $A$ such that $K$ fixes a unique vertex of $S_0$, one can deduce a splitting $S_0'$ of $W_n$ such that $K$ fixes a point of $S_0'$. Indeed, by blowing-up the splitting $S_0$ at the vertex $v$ of $S$ whose associated group is $A$, and by attaching the edges of $S$ adjacent to $v$ to the vertex fixed by $K$, we obtain a splitting $S_0'$ of $W_n$ such that $K$ fixes a point of $S_0'$. Let $\mathcal{S}_0'$ be the equivalence class of $S_0'$. We claim that the group $K_1 \cap T$ fixes $\mathcal{S}_0'$. Indeed, let $e_0$ be the edge of $S_0'$ adjacent to the vertex $v_0$ fixed by $K$ and the vertex fixed by $\left\langle x_n \right\rangle$. Since the stabilizer of $e_0$ is trivial, Proposition~\ref{Prop Levitt stab} implies that the group of twists about $e_0$ at the vertex $v_0$ contains all the twists whose twistor is an element of $K$. Hence $K_1 \cap T$ fixes $\mathcal{S}_0'$.

We now construct a one-edge free splitting $S_0$ of $A$ such that $K$ fixes a vertex of $S_0$. By the above discussion, this will give a two-edge free splitting of $W_n$ such that $K$ fixes a vertex of this splitting which is not a leaf and whose equivalence class is fixed by $K_1 \cap T$. Moreover, we prove at the same time that there does not exist a free splitting of $W_n$ with at least $3$ orbits of edges and such that $K$ fixes a vertex of this splitting. We distinguish between three cases, according to whether $A$ is one-ended relative to $K$ and according to the edge stabilizers of a splitting of $A$ relative to $K$.

\medskip

\noindent{\bf Case 1. } There exists a free splitting $S_0$ of $A$ such that $K$ fixes a vertex of $S_0$. 

In particular, the corresponding splitting $S_0'$ of $W_n$ constructed above is a free splitting of $W_n$.
We claim that the splitting $S_0'$ has two orbits of edges. Indeed, suppose that $S_0'$ has $k$ orbits of edges, with $k \geq 3$. Then, $S_0'$ is obtained from $S$ by blowing-up at least two orbits of edges at $v$. Therefore, the group of twistors $K$ is contained in a free factor $B$ of $W_n$ isomorphic to $W_{n-3}$. Let $B'$ be a free factor of $W_n$ isomorphic to $W_2$ such that 
$$W_n=\left\langle x_n \right\rangle \ast B \ast B'$$ and let $R$ be the free splitting associated with this decomposition. Then the equivalence class $\mathcal{R}$ of $R$ is a free splitting of $W_n$ fixed by $K_1 \cap T$. But by Proposition~\ref{Prop Levitt stab}~$(3)$, the group of twists of $\mathcal{R}$ is isomorphic to $B \times B \times W_2$. Moreover, the group $K_1 \cap T$ is contained in one of the factors of $B \times B \times W_2$ isomorphic to $B$. Therefore, the centralizer of $K_1 \cap T$ contains a free abelian group of rank $2$. Since $K_1 \cap T$ is a normal subgroup of $K_1$, this contradicts the fact that the centralizer of $K_1 \cap T$ is virtually a nonabelian free group by Property~$(P_{W_{n-2}})~(1)$. Therefore, the splitting $S_0'$ is a two-edge free splitting.

\medskip

\noindent{\bf Case 2 } There exists a splitting $S_0$ of $A$ such that $K$ fixes a vertex of $S_0$ and such that one of the edge stabilizers of $S_0$ is finite.

Let $S_0'$ be the corresponding splitting of $W_n$ constructed in the above discussion. If $S_0$ has an edge $e'$ with trivial stabilizer, then by collapsing every orbit of edges of $S_0$ except the one containing $e'$, we obtain a splitting $S_1$ of $A$ such that $K$ fixes a vertex of $K$. Then the corresponding splitting $S_1'$ of $W_n$ is a free splitting. Thus, we can apply Case~$1$.

Therefore, we may assume that every edge stabilizer of $S_0$ is infinite or a nontrivial finite subgroup of $W_n$. By collapsing every edge of $S_0$ with infinite stabilizer and by collapsing all but one orbit of edges with finite edge stabilizer, we may suppose that $S_0$ is a one-edge splitting such that every edge stabilizer of $S_0$ is a nontrivial finite subgroup of $W_n$. Every finite subgroup of $W_n$ is isomorphic to $F$ and is in fact a free factor of $W_n$. We claim that we can construct a splitting $X_0$ of $A$ which contains an edge with trivial stabilizer and such that $K$ fixes a vertex of $X_0$. Indeed, let $x_0$ be the vertex of $S_0$ fixed by $K$, let $f_0$ be an edge adjacent to $x_0$ and let $x_1$ be the vertex of $f_0$ distinct from $v_0$. Let $G_{x_0}$ be the stabilizer of $x_0$, let $G_{x_1}$ be the stabilizer of $x_1$ and let $G_{f_0}$ be the stabilizer of $f_0$. Note that, since there does not exist HNN extensions in $W_n$, the groups $G_{x_0}$ and $G_{x_1}$ are not conjugate in $W_n$. The group $G_{f_0}$ is a free factor of both $G_{x_0}$ and $G_{x_1}$. Thus there exists a free factor $A'$ of $G_{x_1}$ such that $G_{x_1}=G_{f_0} \ast A'$. Let $U$ be the splitting of $A$ such that the underlying tree of $W_n \backslash U$ is the same one as the underlying tree of $W_n \backslash S_0$, such that the stabilizer of every vertex which is not in the orbit of $x_1$ is the same one as the stabilizer of the corresponding vertex in $S_0$ and the stabilizer of $x_1$ is $A'$. Then the edge $f_0$ has trivial stabilizer in $U$ and $K$ fixes a vertex of $U$. This proves the claim. Therefore Case~$2$ is a consequence of Case~$1$.

\medskip

\noindent{\bf Case 3 } The group $A$ is one-ended relative to $K$. 

By Theorem~\ref{Theo canonical splitting for one-ended group}, there exists a canonical splitting $S_0$ of $A$ whose edge stabilizers are virtually infinite cyclic, such that $K$ fixes a point of $S_0$ and such that every automorphism of $A$ preserving $K$ fixes the equivalence class of $S_0$. Up to collapsing some orbits of edges, we may suppose that $S_0$ has exactly one orbit of edges and that $S_0$ is fixed by a finite index subgroup of the group of automorphisms of $A$ preserving $K$. Let $S_0'$ be the corresponding splitting of $W_n$, and let $\mathcal{S}_0'$ be its equivalence class. Recall that the group $K_1 \cap T$ is a normal subgroup of $\Inn(A)$. Hence the group $H$ viewed as a subset of $\Aut(A)$ preserves $K$. Thus $H$ has a finite index subgroup $H'$ which preserves $\mathcal{S}_0'$. 

Let $v_0$ be the vertex of $S_0'$ fixed by $K$, let $e_0$ be the edge of $S_0'$ between $v_0$ and the point fixed by $\left\langle x_n \right\rangle$, let $e$ be an edge adjacent to $v_0$ which is not in the orbit of $e_0$ and let $w_0$ be the endpoint of $e$ distinct from $v_0$. By construction, the stabilizer of every edge of $S_0'$ which is not in the orbit of $e_0$ is virtually cyclic, that is it is either isomorphic to $\ZZ$ or to $W_2$. By Lemma~\ref{Lem twists about cyclic edge central}, a twist about an edge whose stabilizer is isomorphic to $\ZZ$ is central in a finite index subgroup of $\Stab_{\Out(W_n)}(\mathcal{S}_0')$. Since any finite index subgroup of $H$ has finite center by Remark~\ref{Rmq centralizer PWn-2}~$(2)$, we see that the stabilizer of every edge of $S_0'$ which is not in the orbit of $e_0$ is isomorphic to $W_2$. In particular, the stabilizer of $e$ is isomorphic to $W_2$. Therefore, Remark~\ref{Remark proposition Levitt} implies that the group of bitwists about every edge of $S_0'$ which is not in the orbit of $e_0$ is trivial. Thus, the group of bitwists $T_0$ of $\mathcal{S}_0'$ is reduced to the group of twists about $e_0$. Since the group of twists about $e_0$ is virtually a nonabelian free group, the abelian rank of $T_0$ is equal to $1$. Moreover, by Corollary~\ref{Coro splittings over W2}~$(2)$, there exist $k_{v_0},k_{w_0} \in \NN^*$ such that the groups $G_{v_0}$ and $G_{w_0}$ are isomorphic to $W_{k_{v_0}}$ and $W_{k_{w_0}}$, with $3 \leq k_{v_0},k_{w_0} \leq n-2$ and $k_{v_0}+k_{w_0}=n-1+2=n+1$.

By Proposition~\ref{Prop Levitt stab} and Remark~\ref{Remark proposition Levitt}, we have a natural homomorphism $$\Psi \colon H' \to \Out(G_{v_0}) \times \Out(G_{w_0})$$ whose image is contained in $\Out(G_{v_0},\{K,G_e\}) \times \Out(G_{w_0},G_e)$ and whose kernel is $T_0 \cap H'$. Since $H'$ contains a subgroup isomorphic to $\ZZ^{n-2}$ by the third part of Property~$(P_{W_{n-2}})$ and since $T_0$ is virtually a nonabelian free group, we see that the abelian rank of $\Out(G_{v_0},\{K,G_e\})\times \Out(G_{w_0},G_e)$ is at least equal to $n-3$. Recall that, by Proposition~\ref{Prop vcd}, for every $m \geq 3$, the abelian rank of $\Out(W_m)$ is equal to $m-2$. Since $k_{v_0}+k_{w_0}=n+1$ and since the abelian rank of a direct product is the sum of the abelian ranks of the factors, we see that the abelian rank of $\Out(G_{v_0},\{K,G_e\})\times \Out(G_{w_0},G_e)$ is at most equal to $n-3$. Therefore, the abelian rank of $\Out(G_{v_0},\{K,G_e\})\times \Out(G_{w_0},G_e)$ is equal to $n-3$. Thus the group $\Out(G_{v_0},\{K,G_e\})\times \Out(G_{w_0},G_e)$ contains a free abelian subgroup whose rank is equal to the abelian rank of $\Out(A)$. Lemma~\ref{Lem rank product splitting over W2}~$(1)$ then shows that the abelian rank of $\Out(G_{v_0},\{K,G_e\})$ is equal to the abelian rank of $\Out(G_{v_0})$. Moreover, Lemma~\ref{Lem rank product splitting over W2}~$(2)$ implies that there exists a refinement $S_1'$ of $S_0'$ by blowing-up a one-edge $\mathcal{Z}_{RC}$-splitting $U_1$ of $G_{w_0}$ at $w_0$. Let $\mathcal{S}_1'$ be the equivalence class of $S_1'$. Note that, since the group of twists about the edge $e_0$ of $S_0'$ is contained in the group of twists of $\mathcal{S}_1'$, the group $K_1 \cap T$ fixes $\mathcal{S}_1'$. Note also that the stabilizer of the vertex of $S_1'$ fixed by $K$ is equal to $G_{v_0}$. If $U_1$ is a free splitting, then $S_1'$ is a splitting of $W_n$ with one trivial edge stabilizer and such that $K$ fixes a point of $S_1'$. In this case, after collapsing every orbit of edges of $S_1'$ with nontrivial stabilizer, we can apply Case~$1$ to conclude. So we may assume that $S_1'$ has an edge stabilizer isomorphic to $\ZZ$. By Lemma~\ref{Lem twists commuting element twistor}, the twist $D$ about this edge is central in a finite index subgroup of $\Stab(\mathcal{S}_1')$. Hence $D$ commutes with a finite index subgroup of $K_1 \cap T$.

We now prove that the centralizer of the group $K_1 \cap T$ contains a free abelian group of rank $2$. This will contradict Property~$(P_{W_{n-2}})~(1)$. Since $W_{k_{v_0}}$ is a hyperbolic group, and since $G_{v_0}$ is one-ended relative to $K$ and $G_e$, one may apply Theorem~\ref{Theo canonical splitting for one-ended group} to $(G_{v_0},\{K,G_e\})$ to obtain a canonical one-edge splitting $U_2$ of $G_{v_0}$ such that both $K$ and $G_e$ fixes a vertex. The edge stabilizers of $U_2$ are all virtually infinite cyclic. Let $S_2'$ be the refinement of $S_1'$ obtained by blowing-up $U_2$ at the vertex of $S_1'$ fixed by $K$ and by attaching the edges according to their fixed points in $U_2$. Let $\mathcal{S}_2'$ be the equivalence class of $S_2'$. 

Suppose first that one of the edge stabilizer of $U_2$ is isomorphic to $\ZZ$. Note that this case always happen when $k_{v_0}=3$ by Corollary~\ref{Coro splittings over W2}~$(3)$. Then $S_2'$ contains two edges in distinct orbits whose edge stabilizers are isomorphic to $\ZZ$. By Lemma~\ref{Lem twists about cyclic edge central}, the group $\Stab(\mathcal{S}_2')$ has a finite index subgroup which contains a central subgroup isomorphic to $\ZZ^2$. In particular, the centralizer of a finite index subgroup of $K_1 \cap T$ contains a free abelian group of rank $2$. This contradicts Property~$(P_{W_{n-2}})~(1)$.

So we may suppose that $k_{v_0} \geq 4$ and that every edge stabilizer of $U_2$ is isomorphic to $W_2$. By Remark~\ref{Remark proposition Levitt}, the group of bitwists of the equivalence class of $U_2$ trivial. Since $U_2$ is canonical, $U_2$ is fixed by a finite index subgroup $H_2$ of $\Out(G_{v_0},\{K,G_e\})$. Consider the natural homomorphism $$\Psi_2 \colon H_2 \to \prod_{x \in G_{v_0} \backslash U_2} \Out(G_x)$$ given by the action on the vertex groups. By Proposition~\ref{Prop Levitt stab} and Remark~\ref{Remark proposition Levitt}, the kernel of this homomorphism is the group of bitwists, which is trivial. Recall that the abelian rank of $\Out(G_{v_0},\{K,G_e\})$ is equal to the abelian rank of $\Out(G_{v_0})$. Therefore the abelian rank of the image of $\Psi_2$ is equal to the abelian rank of $\Out(G_{v_0})$. By Lemma~\ref{Lem rank product splitting over W2}~$(2)$, the splitting $U_2$ has a refinement $U_3$ obtained by blowing-up a $\mathcal{Z}_{RC}$-splitting at a vertex of $U_2$ which is not in the orbit of the vertex fixed by $K$. This gives a refinement $S_3'$ of $S_1'$ by blowing-up $U_3$ at the vertex of $S_1'$ fixed by $K$. Let $\mathcal{S}_3'$ be the equivalence class of $S_3'$. Note that, as $K$ fixes a point in $S_3'$, and as $\left\langle x_n \right\rangle$ is adjacent to the vertex fixed by $K$, we see that $K_1 \cap T$ is contained in the group of twists of $\mathcal{S}_3'$. If $U_3$ has a trivial edge stabilizer, then after collapsing every orbit of edges in $S_3'$ with nontrivial edge stabilizer, one can apply Case~$1$ to conclude. Otherwise, $U_3$ has an edge stabilizer isomorphic to $\ZZ$. In this case, we see that $S_3'$ contains two edges in distinct orbits whose edge stabilizers are isomorphic to $\ZZ$. By Lemma~\ref{Lem twists about cyclic edge central}, the group $\Stab(\mathcal{S}_2')$ contains a finite index subgroup with a central subgroup isomorphic to $\ZZ^2$. In particular, the centralizer of a finite index subgroup of $K_1 \cap T$ contains a free abelian group of rank $2$. This contradicts Property~$(P_{W_{n-2}})~(1)$. Therefore, there exists a refinement $S_0^{(2)}$ of $S_0'$ such that $K$ fixes a vertex of $S_0^{(2)}$ and $S_0^{(2)}$ has an edge with trivial edge stabilizer. After collapsing every orbit of edges of $S_0^{(2)}$ with nontrivial stabilizer, we can apply Case~$1$ to conclude. The conclusion in Case~$3$ follows.

\medskip

Therefore, we have constructed a free splitting $S_0'$ of $W_n$ which is a two-edge free splitting fixed by $K_1 \cap T$. Moreover, the construction of the splitting is such that the vertex of the underlying graph of $W_n \backslash S_0'$ whose associated group contains $K$ is not a leaf. We now prove that $S_0'$ is a $W_{n-2}$-star. Let $C$ be the vertex stabilizer of $S_0'$ containing $K$, and let $C'$ be a vertex stabilizer of $S_0'$ which is not a conjugate of $C$ nor $\left\langle x_n \right\rangle$. Then $C'$ is the vertex group of a leaf of the underlying graph of $W_n \backslash S_0'$. By Proposition~\ref{Prop Levitt stab}~$(3)$, the group of twists of $\mathcal{S}_0'$ is isomorphic to $C \times C \times C'/Z(C')$. Since the centralizer of $K \cap T_1$ is virtually a nonabelian free group by Property~$(P_{W_{n-2}})~(1)$, we conclude that $C'/Z(C')$ is finite. Hence $C'$ is isomorphic to $F$ and $S_0'$ is a $W_{n-2}$-star. 

We now prove that $H$ virtually fixes $\mathcal{S}_0'$. By Proposition~\ref{Prop Levitt stab}~$(3)$, the group of twists of $\mathcal{S}_0'$ is isomorphic to $W_{n-2} \times W_{n-2}$. By Lemma~\ref{Lem twists Wn-1 star in stab Wn-2}, the group $K_1 \cap T$ is contained in one of the factors isomorphic to $W_{n-2}$ of the group of twists of $\mathcal{S}_0'$. Therefore, $K_1 \cap T$ is centralized by the other factor of the group of twists of $\mathcal{S}_0'$. Since the centralizer of $K_1 \cap  T$ contains $K_2$ as a finite index subgroup, the group $K_2$ contains a twist $f$ of infinite order about the edge $e$ of $S_0'$ which does not collapse onto $S$. This twist is a twist about a $W_{n-1}$-star obtained from $S_0'$ by collapsing the orbit of edges which does not contain $e$. By Lemma~\ref{Lem PWn-2 twists stab Wn-2}, the group $H$ virtually fixes $\mathcal{S}_0'$. Moreover, $K_1$ is commensurable with $T \cap \Stab(\mathcal{S}_0')$, that is $K_1$ is commensurable with the group of twists about one edge of $S_0'$. Lemma~\ref{Lem normal subgroup of twist Wn-2 fixes 2 Wn-1} then implies that $K_1$ virtually fixes a unique equivalence class of $W_{n-2}$-stars. Therefore, since $K_1$ is a normal subgroup of $H$, we see that $H$ virtually fixes a unique equivalence class of $W_{n-2}$-stars. This concludes the proof. 
\hfill\qedsymbol

\begin{prop}\label{Prop characterization Wn-2 stars}
Let $n \geq 5$ and let $\Gamma$ be a finite index subgroup of $C_n$. Let \mbox{$\Psi \in \Comm(\Gamma)$.} Then for every equivalence class $\mathcal{S}$ of $W_{n-2}$-stars, there exists a unique equivalence class $\mathcal{S}'$ of $W_{n-2}$-stars such that $\Psi([\Stab_{\Gamma}(\mathcal{S})])= [\Stab_{\Gamma}(\mathcal{S}')]$.
\end{prop}

\dem The uniqueness statement follows from Lemma~\ref{Lem normal subgroup of twist Wn-2 fixes 2 Wn-1} which shows that the stabilizer in finite index subgroups of $\Out(W_n)$ of two distinct equivalence classes of $W_{n-2}$-stars are not commensurable.

We now prove the existence statement. Let $f \colon \Gamma_1 \to \Gamma_2$ be an isomorphism between finite index subgroups of $\Gamma$ that represents $\Psi$. By Proposition~\ref{Prop Stab Wn-2 satisfies PWn-2}, the group $\Stab_{\Gamma_1}(\mathcal{S})$ satisfies $(P_{W_{n-2}})$. As $f$ is an isomorphism, we deduce that $f(\Stab_{\Gamma_1}(\mathcal{S}))$ also satisfies $(P_{W_{n-2}})$. Proposition~\ref{Prop QWn-2 is stab Wn-2} implies that there exists a unique equivalence class of $W_{n-2}$-stars $\mathcal{S}'$ such that $f(\Stab_{\Gamma_1}(\mathcal{S})) \subseteq \Stab_{\Gamma_2}(\mathcal{S}')$, where the inclusion holds up to a finite index subgroup. Applying the same argument with $f^{-1}$, we see that there exists an equivalence class $\mathcal{S}''$ of a $W_{n-2}$-star such that $$\Stab_{\Gamma_1}(\mathcal{S}) \subseteq f^{-1}(\Stab_{\Gamma_2}(\mathcal{S}')) \subseteq \Stab_{\Gamma_1}(\mathcal{S}''),$$ where the inclusion holds up to a finite index subgroup. Lemma~\ref{Lem normal subgroup of twist Wn-2 fixes 2 Wn-1} then implies that $\mathcal{S}$ is the unique equivalence class of $W_{n-2}$-stars virtually fixed by $\Stab_{\Gamma_1}(\mathcal{S})$. Therefore, we see that $\mathcal{S}=\mathcal{S}''$ and we have equality everywhere. This completes the proof.
\hfill\qedsymbol

\section{Algebraic characterization of compatibility of $W_{n-2}$-stars and conclusion}

\subsection{Algebraic characterization of compatibility of $W_{n-2}$-stars}

In this section, we give an algebraic characterization of the fact that two equivalence classes of $W_{n-2}$-stars have both a common collapse and a common refinement. This will imply that $\Comm(\Out(W_n))$ preserves the set of pairs of commensurability classes of stabilizers of adjacent pairs in the graph $X_n$ introduced in Definition~\ref{Defi graph Xn}~$(2)$.

Let $n \geq 5$ and let $\Gamma$ be a finite index subgroup of $C_n$. We consider the following properties of a pair $(H_1,H_2)$ of subgroups of $\Gamma$.

\bigskip

\noindent $(P_{\mathrm{comp}})$ The pair $(H_1,H_2)$ satisfies the following properties.

\medskip

\noindent{$(1)$ } For every $i \in \{1,2\}$, the group $H_i$ satisfies $(P_{W_{n-2}})$.

\medskip

\noindent{$(2)$ } For every normal subgroups $K_1^{(1)} \times K_2^{(1)}$ of $H_1$ and $K_1^{(2)} \times K_2^{(2)}$ of $H_2$ given by $(P_{W_{n-2}})~(1)$, there exist $i,j \in \{1,2\}$ such that $K_i^{(1)} \cap K_j^{(2)}$ is infinite.

\medskip

\noindent{$(3)$ } The group $H_1 \cap H_2$ contains a subgroup isomorphic to $\ZZ^{n-2}$.

\bigskip

\begin{prop}\label{Prop Pcomp and compatibility}
Let $n \geq 5$ and let $\Gamma$ be a finite index subgroup of $C_n$. Let $\mathcal{S}_1$ and $\mathcal{S}_2$ be two distinct equivalence classes of $W_{n-2}$-stars $S_1$ and $S_2$ and, for every $i \in \{1,2\}$, let $H_i=\Stab_{\Gamma}(\mathcal{S}_i)$. Then $S_1$ and $S_2$ have a refinement $S$ which is a $W_{n-3}$-star if and only if $(H_1,H_2)$ satisfies Property~$(P_{\mathrm{comp}})$. 
\end{prop}

\dem We first assume that $S_1$ and $S_2$ have a common refinement $S$ which is a $W_{n-3}$-star. Let $\mathcal{S}$ be the equivalence class of $S$. Let us prove that $(H_1,H_2)$ satisfies $(P_{\mathrm{comp}})$. By Proposition~\ref{Prop Stab Wn-2 satisfies PWn-2}, for every $i \in \{1,2\}$, the group $H_i$ satisfies $(P_{W_{n-2}})$. This proves that the pair $(H_1,H_2)$ satisfies $(P_{\mathrm{comp}})~(1)$. 

Let us check Property~$(P_{\mathrm{comp}})~(2)$. For every $i \in \{1,2\}$, let $T_1^{(i)} \times T_2^{(i)}$ be the group of twists of $\mathcal{S}_i$ and let $K_1^{(i)}=T_1^{(i)} \cap \Gamma$ and $K_2^{(i)}=T_2^{(i)} \cap \Gamma$. By Proposition~\ref{Prop Stab Wn-2 satisfies PWn-2}, for every $i \in \{1,2\}$, the group $K_1^{(i)} \times K_2^{(i)}$ satisfies $(P_{W_{n-2}})~(1)$ and Lemma~\ref{Lem PWn-2 twists stab Wn-2} implies that every normal subgroup of $H_i$ given by $(P_{W_{n-2}})~(1)$ is commensurable with $K_1^{(i)} \times K_2^{(i)}$. Thus it suffices to check $(P_{\mathrm{comp}})~(2)$ for $K_1^{(1)} \times K_2^{(1)}$ and $K_1^{(2)} \times K_2^{(2)}$. The group of twists of $\mathcal{S}$ is isomorphic to a direct product $A_1 \times A_2 \times A_3$ of three copies of $W_{n-3}$.  Since $n \geq 5$, we have $n-3 \geq 2$ and $W_{n-3}$ is infinite. Since $S$ is a common refinement of $S_1$ and $S_2$ and since $S$ has $3$ orbits of edges there exists a $W_{n-1}$-star $S_0$ which is a common collapse of $S_1$ and $S_2$. Moreover, there exists $k \in \{1,2,3\}$ such that $A_k$ is contained in the group of twists of $S_0$. Therefore, for every $i \in \{1,2\}$, there exists $j \in \{1,2\}$ such that the group $A_k$ is contained in $T_j^{(i)}$. Thus, there exist $i,j \in \{1,2\}$ such that $A_k \cap \Gamma \subseteq K_i^{(1)} \cap K_{j}^{(2)}$. In particular, $K_i^{(1)} \cap K_{j}^{(2)}$ is infinite. This shows $(P_{\mathrm{comp}})~(2)$. 

Finally, since $n \geq 5$, the $W_{n-2}$-stars $S_1$ and $S_2$ have a common refinement which is a $W_2$-star (take any $W_2$-star which refines $S$). Since the group of twists of a $W_2$-star contains a subgroup isomorphic to $\ZZ^{n-2}$ by Proposition~\ref{Prop Levitt stab}~$(3)$, this shows $(P_{\mathrm{comp}})~(3)$.

Conversely, suppose that $(H_1,H_2)$ satisfies $(P_{\mathrm{comp}})$. For $i \in \{1,2\}$, let $K_1^{(i)} \times K_2^{(i)}$ be the direct product of the groups of twists in $\Gamma$ about the two edges of $\mathcal{S}_i$. Then for every $i \in \{1,2\}$, the group $\big(H_i \cap K_1^{(i)}\big) \times \big(H_i \cap K_2^{(i)}\big)$ satisfies $(P_{W_{n-2}})~(1)$ by Proposition~\ref{Prop Stab Wn-2 satisfies PWn-2}. Hence Property~$(P_{\mathrm{comp}})~(2)$ implies that there exists $i,j \in \{1,2\}$ such that $$\left(H_1 \cap K_i^{(1)} \right)\cap \left(H_2 \cap K_j^{(2)}\right)$$ is infinite. For $i \in \{1,2\}$, let $S_1^{(i)}$ and $S_2^{(i)}$ be the two distinct $W_{n-1}$-stars on which $S_i$ collapses. By Proposition~\ref{Prop PWn-2 compatibility}, since $H_1 \cap H_2$ fixes pointwise the set $\{\mathcal{S}_1^{(1)},\mathcal{S}_2^{(1)},\mathcal{S}_1^{(2)},\mathcal{S}_2^{(2)}\}$, and since $H_1 \cap H_2$ contains a subgroup isomorphic to $\ZZ^{n-2}$ by $(P_{\mathrm{comp}})~(3)$, the $W_{n-1}$-stars $S_1^{(1)}$, $S_2^{(1)}$, $S_1^{(2)}$ and $S_2^{(2)}$ are pairwise compatible. Hence $S_1$ and $S_2$ have a common refinement $S$ which is either a $W_{n-3}$-star or a $W_{n-4}$-star. Since the groups of twists of $\mathcal{S}_1$ and $\mathcal{S}_2$ have infinite intersection, the refinement $S$ cannot be a $W_{n-4}$-star since otherwise the $W_{n-1}$-stars $S_1^{(1)}$, $S_2^{(1)}$, $S_1^{(2)}$ and $S_2^{(2)}$ would be pairwise nonequivalent and hence their groups of twists would have trivial intersection. Thus $S$ is a $W_{n-3}$-star. This concludes the proof.
\hfill\qedsymbol

\subsection{Conclusion}

In this last section, we complete the proof of our main theorem.

\begin{theo}\label{Theo Comm Wn complete statement}
Let $n \geq 5$ and let $\Gamma$ be a finite index subgroup of $C_n$. Then any isomorphism $f \colon H_1 \to H_2$ between two finite index subgroups of $\Gamma$ is given by conjugation by an element of $\Out(W_n)$ and the natural map: $$\Out(W_n) \to \Comm(\Out(W_n))$$ is an isomorphism. 
\end{theo}

\dem Suppose that $\mathcal{S}$ and $\mathcal{S}'$ are two distinct equivalence classes of $W_{n-2}$-stars. Then $\Stab_{\Gamma}(\mathcal{S})$ and $\Stab_{\Gamma}(\mathcal{S}')$ are not commensurable by Lemma~\ref{Lem normal subgroup of twist Wn-2 fixes 2 Wn-1}. Proposition~\ref{Prop characterization Wn-2 stars} shows that the collection $\mathcal{I}$ of all commensurability classes of $\Gamma$-stabilizers of equivalence classes of $W_{n-2}$-stars is $\Comm(\Gamma)$-invariant. Proposition~\ref{Prop Pcomp and compatibility} shows that the collection $\mathcal{J}$ of all pairs $([\Stab_{\Gamma}(\mathcal{S})],[\Stab_{\Gamma}(\mathcal{S}')])$ is also $\Comm(\Gamma)$-invariant. Since the natural homomorphism $\Out(W_n) \to \Aut(X_n)$ is an isomorphism by Theorem~\ref{Theo rigidity graph Wn-2}, the conclusion follows from Proposition~\ref{Prop isomorphism commensuration} and the fact that $\Comm(\Gamma)$ is isomorphic to $\Comm(\Out(W_n))$ since $\Gamma$ has finite index in $\Out(W_n)$.
\hfill\qedsymbol

\appendix

\section{Rigidity of the graph of $W_{n-1}$-stars}

The \emph{graph of $W_{n-1}$-stars}, denoted by $Y_n$, is the graph whose vertices are the $W_n$-equivariant homeomorphism classes of $W_{n-1}$-stars, where two equivalence classes $\mathcal{S}$ and $\mathcal{S}'$ are joined by an edge if there exist $S \in \mathcal{S}$ and $S' \in \mathcal{S}'$ such that $S$ and $S'$ are compatible. This graph arises naturally in the study of $\Out(W_n)$ as it is isomorphic to the full subgraph of the free splitting graph $\overline{K}_n$ of $W_n$ whose vertices are equivalence classes of $W_k$-stars, with $k$ varying in $\{0,\ldots,n-1\}$. As $\Aut(W_n)$ acts on $\overline{K}_n$ by precomposition of the marking, we have an induced action of $\Aut(W_n)$ on $Y_n$. As $\Inn(W_n)$ acts trivially on $Y_n$, the action of $\Aut(W_n)$ induces an action of $\Out(W_n)$. We denote by $\Aut(Y_n)$ the group of graph automorphisms of $Y_n$. In this section we prove the following theorem.

\begin{theo}\label{Theo rigidity graph Wn-1}
Let $n \geq 4$. The natural homomorphism $$\Out(W_n) \to \Aut(Y_n)$$ is an isomorphism.
\end{theo}

In order to prove this theorem, we take advantage of the action of $\Out(W_n)$ on the graph of $\{0\}$-stars and $F$-stars $L_n$. The strategy in order to prove Theorem~\ref{Theo rigidity graph Wn-1} is to construct an injective homomorphism \mbox{$\Phi \colon \Aut(Y_n) \to \Aut(L_n)$} such that every automorphism in the image preserves the set of $\{0\}$-stars and the set of $F$-stars.

The homomorphism $\Phi \colon \Aut(Y_n) \to \Aut(L_n)$ is defined as follows. Let $f \in \Aut(Y_n)$. Let $\mathcal{S}$ be the equivalence class of a $\{0\}$-star and let $S$ be a representative of $\mathcal{S}$. By Theorem~\ref{Theo scott swarup}, there exist exactly $n$ $W_{n-1}$-stars $S_1,\ldots,S_n$ refined by $S$. Moreover, these $W_{n-1}$-stars are pairwise compatible. For $i \in \{1,\ldots,n\}$, let $\mathcal{S}_i$ be the equivalence class of $S_i$. Since $f$ is an automorphism of $Y_n$, $f(\mathcal{S}_1),\ldots,f(\mathcal{S}_n)$ are pairwise adjacent in $Y_n$. Let $S_1',\ldots,S_n'$ be representatives of respectively $f(\mathcal{S}_1),\ldots,f(\mathcal{S}_n)$ that are pairwise compatible. Then Theorem~\ref{Theo scott swarup} implies that there exists a unique common refinement $S'$ of $S_1',\ldots,S_n'$ with exactly $n$ edges. Since, for every $i \in \{1,\ldots,n\}$, the splitting $S_i'$ is a $W_{n-1}$-star, the splitting $S'$ is necessarily a $\{0\}$-star. Let $\mathcal{S}'$ be the equivalence class of $S'$. We then define $\Phi(f)(\mathcal{S})=\mathcal{S}'$. If $\mathcal{T}$ is an $F$-star, we define $\Phi(f)(\mathcal{T})$ similarly.

\begin{lem}\label{Lem injectivity morphism Xn to Ln}
Let $n \geq 4$. Let $f \in \Aut(Y_n)$. Let $\Phi(f)$ be as above.

\medskip

\noindent{$(1)$ } The map $\Phi(f) \colon VL_n \to VL_n$ induces a graph automorphism $\widetilde{\Phi}(f) \colon L_n \to L_n$.

\medskip

\noindent{$(2)$ } If $\widetilde{\Phi}(f)=\mathrm{id}_{L_n}$, then $f=\mathrm{id}_{Y_n}$.
\end{lem}

\dem We prove the first statement. As $\Phi(f) \circ \Phi(f^{-1})=\Phi(f \circ f^{-1})=\mathrm{id}$, we see that $\Phi(f)$ is a bijection. Let $\mathcal{S}$ be the equivalence class of a $\{0\}$-star and let $\mathcal{T}$ be the equivalence class of an $F$-star. Suppose that $\mathcal{S}$ and $\mathcal{T}$ are adjacent in $L_n$. We prove that $\Phi(f)(\mathcal{S})$ and $\Phi(f)(\mathcal{T})$ are adjacent in $L_n$. Applying the same result to $f^{-1}$, this will prove that $\mathcal{S}$ and $\mathcal{T}$ are adjacent in $L_n$ if and only if $\Phi(f)(\mathcal{S})$ and $\Phi(f)(\mathcal{T})$ are adjacent in $L_n$, and this will conclude the proof. Let $S$ and $T$ be representatives of respectively $\mathcal{S}$ and $\mathcal{T}$. Let $S_1,\ldots,S_n$ be the $n$ $W_{n-1}$-stars refined by $S$, and let $T_1,\ldots,T_{n-1}$ be the $n-1$ $W_{n-1}$-stars refined by $T$. As $S$ refines $T$, and as $S$ refines exactly $n$ $W_{n-1}$-stars by Theorem~\ref{Theo scott swarup}, up to reordering, we can suppose that, for every $i \in \{1,\ldots,n-1\}$, we have $S_i=T_i$. For $i \in \{1,\ldots,n\}$, let $\mathcal{S}_i$ be the equivalence class of $S_i$, and let $S_i'$ be a representative of $\Phi(f)(\mathcal{S}_i)$ such that for distinct $i,j \in \{1,\ldots,n\}$, $S_i$ and $S_j$ are compatible. Then, by Theorem~\ref{Theo scott swarup}, a representative $T'$ of $\Phi(f)(\mathcal{T})$ is the unique (up to $W_n$-equivariant homomophism) $F$-star such that, for every $j \in \{1,\ldots,n-1\}$, $T'$ is compatible with $S_j'$. Moreover, a representative $S'$ of $\Phi(f)(\mathcal{S})$ is the unique $\{0\}$-star such that, for every $i \in \{1,\ldots,n\}$, $S'$ is compatible with $S_i'$. For $i \in \{1,\ldots,n\}$, let $x_i$ be the preimage by the marking of $W_n \backslash S_i'$ (well defined up to global conjugation) of the generator of the vertex group isomorphic to $F$ (which exists since $S_i'$ is a $W_{n-1}$-star). Then the preimages by the marking of $W_n \backslash T'$ of the generators of the groups associated with the $n-1$ leaves of the underlying graph of $W_n \backslash T'$ are $x_1,\ldots,x_{n-1}$ and the preimage by the marking of $W_n \backslash T'$ of the generator of the group associated with the center of the underlying graph of $W_n \backslash T'$ is $x_n$. Moreover, the preimages by the marking of $W_n \backslash S'$ of the generators of the groups associated with the $n$ leaves of the underlying graph of $W_n \backslash S'$ are $x_1,\ldots,x_{n}$. Let $v_n$ be the leaf of the underlying graph of $W_n \backslash S'$ such that the preimage by the marking of $W_n \backslash S'$ of the generator of the group associated with $v_n$ is $x_n$. Then $T'$ is obtained from $S'$ by contracting the edge adjacent to $v_n$. Thus $\Phi(f)(\mathcal{S})$ and $\Phi(f)(\mathcal{T})$ are adjacent in $L_n$.

\bigskip

The proof of the second statement is identical to the proof of \cite[Lemma~5.4]{Guerch2020symmetries}. We add the proof for completeness as the statement of \cite[Lemma~5.4]{Guerch2020symmetries} is about automorphisms of $\overline{K}_n$. Let $\mathcal{S} \in VY_n$ and let $S$ be a representative of $\mathcal{S}$. We prove that $f(\mathcal{S})=\mathcal{S}$. Let 
$$W_n=\left\langle x_1,\ldots, x_{n-1} \right\rangle \ast \left\langle x_n \right\rangle$$ 
be the free factor decomposition of $W_n$ induced by $S$. Let $S'$ be a representative of $f(\mathcal{S})$. Let $\mathcal{X}$ be the equivalence class of the $F$-star $X$ represented in Figure~\ref{proof X} on the left. 

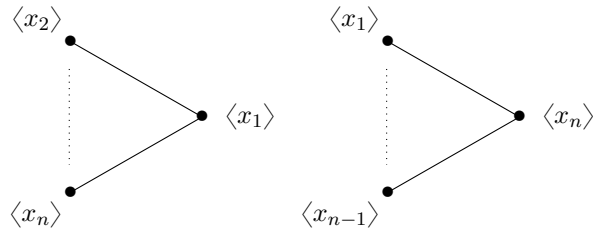
\begin{figure}[ht]
\centering
\begin{tikzpicture}[scale=2]
\draw (0:0) -- (150:1);
\draw (0:0) -- (210:1);
\draw (0:0) node {$\bullet$};
\draw (210:1) node {$\bullet$};
\draw (150:1) node {$\bullet$};
\draw[dotted] (160:0.93) -- (200:0.93);
\draw (0:0) node[right, scale=0.9] {$\;\;\left\langle x_1 \right\rangle$};

\draw (150:1) node[above left, scale=0.9] {$\left\langle x_2 \right\rangle$};
\draw (210:1) node[below left, scale=0.9] {$\left\langle x_n \right\rangle$};
\end{tikzpicture}
\begin{tikzpicture}[scale=2]
\draw (0:0) -- (150:1);
\draw (0:0) -- (210:1);
\draw (0:0) node {$\bullet$};
\draw (210:1) node {$\bullet$};
\draw (150:1) node {$\bullet$};
\draw[dotted] (160:0.93) -- (200:0.93);
\draw (0:0) node[right, scale=0.9] {$\;\;\left\langle x_n \right\rangle$};

\draw (150:1) node[above left, scale=0.9] {$\left\langle x_1 \right\rangle$};
\draw (210:1) node[below left, scale=0.9] {$\left\langle x_{n-1} \right\rangle$};
\end{tikzpicture}
\caption{The $F$-stars $X$ and $X'$ of the proof of Lemma~\ref{Lem injectivity morphism Xn to Ln}.}\label{proof X}
\end{figure}

Since $\Phi(f)(\mathcal{X})=\mathcal{X}$, the free splitting $S'$ is a $W_{n-1}$-star obtained from $X$ by collapsing $n-1$ edges. But if $T$ is a $W_{n-1}$-star obtained from $X$ by collapsing $n-1$ edges, then there exists $i \in \{1,\ldots,n\}$ such that the free factor decomposition of $W_n$ induced by $T$ is 
$$W_n=\left\langle x_1,\ldots, \widehat{x}_i,\ldots, x_{n} \right\rangle \ast \left\langle x_i \right\rangle.
$$
For $i \in \{1,\ldots,n\}$, we will denote by $T_i$ the $W_{n-1}$-star with associated free factor decomposition $\left\langle x_1,\ldots, \widehat{x}_i,\ldots, x_{n} \right\rangle \ast \left\langle x_i \right\rangle$, and by $\mathcal{T}_i$ its equivalence class. For $i \neq n$, the free splitting $T_i$ is a collapse of the $F$-star $X'$ depicted in Figure~\ref{proof X} on the right, whereas $S$ is not a collapse of $X'$.

Let $\mathcal{X}'$ be the equivalence class of $X'$. Since $\Phi(f)(\mathcal{X}')=\mathcal{X}'$, there does not exist a representative of $f(\mathcal{S})$ that is obtained from a representative of $\mathcal{X}'$ by collapsing a forest. Thus, for all $i \neq n$, we have $f(\mathcal{S}) \neq \mathcal{T}_i$. Therefore, as $\mathcal{S}=\mathcal{T}_n$, we conclude that $f(\mathcal{S})=\mathcal{S}$.
\hfill\qedsymbol

\bigskip

\noindent{\bf Proof of Theorem~\ref{Theo rigidity graph Wn-1}. } Let $n \geq 4$. We first prove injectivity. The homomorphism $\Out(W_n) \to \Aut(L_n)$ is injective by Theorem~\ref{Theo rigidity Ln}. Moreover, the homomorphism $\Out(W_n) \to \Aut(L_n)$ factors through $\Out(W_n) \to \Aut(Y_n) \to \Aut(L_n)$. Therefore we deduce the injectivity of $\Out(W_n) \to \Aut(Y_n)$. We now prove surjectivity. Let $f \in \Aut(Y_n)$. By Lemma~\ref{Lem injectivity morphism Xn to Ln}~$(1)$, we have a homomorphism $\Phi \colon \Aut(Y_n) \to \Aut(L_n)$ whose image consists in automorphisms preserving the set of $\{0\}$-stars and the set of $F$-stars. By Theorem~\ref{Theo rigidity Ln}, the automorphism $\Phi(f)$ is induced by an element $\gamma \in \Out(W_n)$. Since the homomorphism $\Aut(Y_n) \to \Aut(L_n)$ is injective by Lemma~\ref{Lem injectivity morphism Xn to Ln}~$(2)$, $f$ is induced by $\gamma$. This concludes the proof. 
\hfill\qedsymbol

\bibliographystyle{alphanum}
\bibliography{bibliographie}

\noindent \begin{tabular}{l}
Laboratoire de mathématique d'Orsay\\
UMR 8628 CNRS \\
Université Paris-Saclay\\
91405 ORSAY Cedex, FRANCE\\
{\it e-mail: yassine.guerch@universite-paris-saclay.fr}
\end{tabular}
\end{document}